\documentclass[a4paper, 10pt]{scrartcl}

\usepackage{doi}
\usepackage{bm} 
\usepackage{hyperref}
\usepackage[top=3cm, right = 2.5cm, bottom =3cm, left = 2.5cm]{geometry}
\usepackage{lmodern}

\usepackage{amsmath, amssymb, amsthm}

\newtheorem{lem}{Lemma}[section]
\newtheorem{thm}[lem]{Theorem}

\newtheorem{cor}[lem]{Corollary}

\theoremstyle{plain}

\theoremstyle{remark}
\newtheorem{rem}[lem]{Remark}
\numberwithin{equation}{section}

\usepackage{enumitem}
\usepackage{xcolor}
\usepackage{pgfplots}
\pgfplotsset{compat=newest}
\usetikzlibrary{plotmarks}
\usetikzlibrary{arrows.meta}
\usepgfplotslibrary{patchplots}
\usepackage{grffile}
\usepackage{graphicx}
\graphicspath{{./Figs/}}
\usepackage{subcaption}

\usepackage[]{csquotes}
\usepackage[authoryear, round]{natbib}
\newcommand\dm{\mathop{}\!\mathrm{d}}
\newcommand\sgn{\mathrm{sgn}}

\newcommand{\argmin}[1]{
	\mathrm{arg}\, \min_{#1}
}
\newcommand{\Lp}[2][\varOmega]{\ensuremath{L^{#2}\left(#1\right)}}

\newcommand{\hzero}[1][\varOmega]{\ensuremath{H^1_0\left(#1\right)}}
\newcommand{\sob}[3][\varOmega]{\ensuremath{W^{#2,#3}\left(#1\right)}}
\newcommand{\sobzero}[3][\varOmega]{\ensuremath{W^{#2,#3}_0\left(#1\right)}}

\begin{document}
\title{Gibbs Phenomena for \texorpdfstring{$\bm{L^q}$}{$L^q$}-Best Approximation in Finite Element Spaces - Some Examples}
\author{Paul Houston$^{a,}$\thanks{Paul.Houston@nottingham.ac.uk}, Sarah Roggendorf$^{a,}$\thanks{Corresponding author: Sarah.Roggendorf@nottingham.ac.uk} and Kristoffer G. van der Zee$^{a,}$\thanks{KG.vanderZee@nottingham.ac.uk} \bigskip\\
{\small $^a$School of Mathematical Sciences, The University of Nottingham},\\ {\small University Park, NG7 2RD, United Kingdom}}

\maketitle
\begin{abstract}
  Recent developments in the context of minimum residual finite element methods are paving the way for designing finite element methods in non-standard function spaces.
    This, in particular, permits the selection of a solution space in which the best approximation of the solution has desirable properties. One of the biggest challenges in designing finite element methods are non-physical oscillations near thin layers and jump discontinuities.
    In this article we investigate Gibbs phenomena in the context of $L^q$-best approximation of discontinuities in finite element spaces with $1\leq q<\infty$. Using carefully selected examples, we show that on certain meshes the Gibbs phenomenon can be eliminated in the limit as $q$ tends to $1$. The aim here is to show the potential of $L^1$ as a solution space in connection with suitably designed  meshes.
  \bigskip

\end{abstract}

\section{Introduction}

This article investigates the  Gibbs phenomenon in the context of the $L^q$-best approximation of discontinuous functions in finite element spaces by considering a few carefully selected simple examples that can be analysed in detail. The Gibbs phenomenon was originally discovered by Henry \cite{Wilbraham1848}, and described by Willard \cite{Gibbs1899} in the context of approximating jump discontinuities by partial sums of Fourier series. It also occurs in the best approximation of functions either by a trigonometric polynomial in the $L^1$-metric \citep{Moskona1995} or spline functions in the $L^2$-metric \citep{Richards1991}. The best approximation in finite element spaces consisting of piecewise polynomials is closely related to the last example. \cite{Saff1999} show  that in one dimension the best approximation of a jump discontinuity by polygonal lines leads to Gibbs phenomena for all $1<q<\infty$
 but vanishes as $q \rightarrow 1$; this is the starting point of our investigation.

 We consider several meshes in one and two dimensions and show that on certain meshes the over- and undershoots in the best approximation can be eliminated in the limit $q \rightarrow 1$. These results are extensions of \cite{Saff1999}. However, there exist meshes in both one and two dimensions that do not satisfy this property.   The aim of this article is therefore to illustrate which properties the underlying mesh must satisfy to ensure that the oscillations vanish in the $L^q$-best approximation of discontinuous functions.

This study of $L^q$-best approximations in finite element spaces is motivated by approximating solutions to partial differential
equations (PDEs) in subspaces of $L^1(\varOmega)$.  \cite{Guermond2004} points out
that there are only very few attempts at achieving this despite the fact that
first-order PDEs and their non-linear generalizations have been extensively studied
in $L^1(\varOmega)$. The existing numerical methods which seek  an approximation
directly in $L^1(\varOmega)$ include the ones outlined in the articles by  \cite{Lavery1988,Lavery1989,Lavery1991},
the reweighted least-squares method of \cite{Jiang1993,Jiang1998} and the methods outlined
in the series of articles by Guermond \emph{et al.}\ \citep{Guermond2004,Guermond2007,Guermond2008,Guermond2008/09,Guermond2009}.
More recently, a novel approach to designing finite element methods in
a very general Banach space setting has been introduced by \cite{Muga2017} and applied to the advection-reaction equation \citep{Muga2019} and to the convection-diffusion-reaction equation \citep{Houston2019b}.
This approach is based on the so-called discontinuous Petrov-Galerkin methods \citep[e.g.,][]{Demkowicz2014} and
extends the concept of optimal test norms and functions from Hilbert spaces to
more general Banach spaces. At least in an abstract sense, this approach outlines
how to design a numerical method that leads to a quasi-best approximation of the
solution in a space of choice, provided the continuous problem is well-posed in a suitable sense.
 In practice, there are hurdles to overcome
to design a practical method, but this is not the subject of this article. Nonetheless,
it opens up a new approach to designing numerical methods that raises the question of which norms
and spaces are favourable for the approximation of certain types of PDEs.

In the context of approximating solutions containing discontinuities and under resolved interior- and boundary layers, the numerical results for existing
$L^1$-methods suggest
such features can be
approximated as sharply as a given mesh permits without exhibiting spurious over- or
undershoots. This property clearly gives them an enormous advantage over traditional finite element
methods yielding approximations in subspaces of $L^2(\varOmega)$. Indeed, it is well-known
that even seemingly simple examples such as the
transport equation or convection-dominated diffusion equations require extra
care in the design of the method, with the standard Galerkin finite element method being
unstable, and alternative methods often requiring so-called \emph{stabilization} and/or \emph{shock-capturing} techniques \citep[e.g.,][]{John2007,John2007a,John2008,Roos2008}.

 \subsection{Notation}
 Throughout this article, we denote by $L^q(\varOmega)$, $1\le q<\infty$, the Lebesgue space of $q$-integrable
 functions on a bounded Lipschitz domain $\varOmega \subset \mathbb{R}^d$, $d\in\{1,2\}$;
 $L^\infty(\varOmega)$ is the Lebesgue space of functions on~$\varOmega$ with finite essential supremum; and $\sob{1}{q}$, $1\le q\le \infty$, is the Sobolev space of functions that are in $\Lp{q}$
such that their gradient is in $\Lp{q}^d$. Furthermore, $\sobzero{1}{q}\subset
\sob{1}{q}$ is the subspace of all functions with zero trace  on the boundary
$\partial\varOmega$. The corresponding norms are denoted by $\|\cdot \|_{\Lp{q}}$
and $\|\cdot\|_{\sob{1}{q}}$, respectively. For $q=2$, we furthermore use the usual
notation $H^1(\varOmega) := \sob{1}{2}$ and $\hzero:= \sobzero{1}{2}$.
 For $1\leq q\leq \infty$, we write $q'$ to denote the dual exponent such that
 $1/q+1/q' = 1$.  For any Banach space $V$,
 its dual space is denoted by $V'$. Furthermore, for $v \in V$ and $\varphi \in V'$, we have
 the duality pairing
 \begin{align*}
   \langle \varphi, v \rangle_{V', V} : = \varphi(v).
 \end{align*}
 The subdifferential of a function $f:V \rightarrow \mathbb{R}$ at a point $v \in V$ is denoted by $\partial f (v) \subset V'$. Furthermore, for $v,w \in V$, we write $\partial f(v)(w)$ to denote $\varphi(w)$ for an arbitrary $\varphi \in \partial f(v)$.

\subsection{Motivation}
{ To motivate the best approximation problem we analyse in this article, we consider the following simple convection-diffusion problem}: find $u$ such that
\begin{align}
  -\varepsilon u'' + u' = 0 \quad \text{ in } (0,1), && u(0) = 1, && u(1) = 0. \label{eq:1d_boundary_layer}
\end{align}
The analytical solution to this problem is given by
\begin{align*}
  u(x) = \frac{1- e^{-\frac{1-x}{\varepsilon}}}{1-e^{-\frac{1}{\varepsilon}}};
\end{align*}
in particular, there is a boundary layer near $x=1$ for small $\varepsilon$.
In two dimensions, we consider a rather straight forward extension of the one-dimensional example: find $u$ such that
\begin{align}
  -\varepsilon \Delta u + \partial_x u = 0 \quad \text{ in }(0,1)^2, && u(0, \cdot ) = 1, && u(1, \cdot) = 0, && \partial_{\bm n} u= 0 \text{ if } y = 0 \text{ or } y=1, \label{eq:2d_boundary_layer}
\end{align}
where $\bm n$ denotes the unit outward normal vector on the boundary of the domain.

We seek an approximation of the analytical solution in a finite dimensional space that
consists of  continuous piecewise linear polynomials defined on a given mesh. In one dimension, we are interested both in uniform and non-uniform meshes. In two dimensions we consider predominantly uniform and structured meshes, although we include one example of an unstructured mesh.

If $\varepsilon \ll 1$, then the second order term is completely dominated by the first-order term and away from the outflow boundary the solution is essentially given by the solution to the advection problem obtained by setting $\varepsilon$ to zero. For the above problems this means that $u \approx 1$ away from the outflow boundary. Due to the Dirichlet boundary conditions, a boundary layer forms near the outflow boundary. If the diameter of the elements near the boundary layer is large compared with $\varepsilon$, the layer is fully contained within these elements and, in the above problems, $u \approx 1$ in the rest of the domain.
 Numerically, this essentially means that we approximate the problems
\eqref{eq:1d_boundary_layer}/\eqref{eq:2d_boundary_layer} with $\varepsilon = 0$ while still keeping the boundary conditions
at both ends. Clearly, the  analytical solution for the above problems with $\varepsilon = 0$ and  the boundary
conditions only imposed on the inflow part of the boundary is $u \equiv 1$.
This motivates us to  consider the best approximations of $u \equiv 1$
by linear finite element functions satisfying the boundary conditions given in \eqref{eq:1d_boundary_layer} and \eqref{eq:2d_boundary_layer}, respectively.

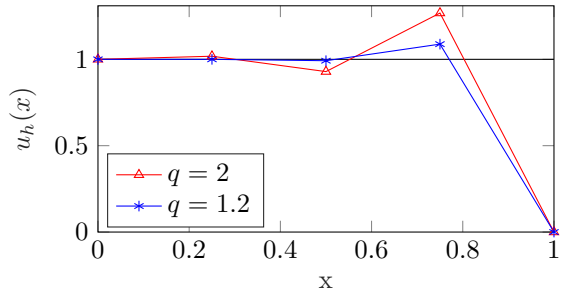
\begin{figure}
	\centering
	\begin{tikzpicture}[scale = 1.0]

		\begin{axis}[%
		width=6cm,
		height=3cm,
		at={(1.011in,0.642in)},
		scale only axis,
		xmin=0,
		xmax=1,
		xlabel style={font=\color{white!15!black}},
		xlabel={x},
		ymin=0,
		ymax=1.31,
		ylabel style={font=\color{white!15!black}},
		ylabel={$u_h(x)$},
		axis background/.style={fill=white},
		legend style={at={(0.02,0.02)}, anchor=south west, legend cell align=left, align=left, draw=white!15!black}
		]

		\addplot [color=black, mark=none, forget plot]
			table[row sep=crcr]{%
			0 1.0\\
			0.25 1.000\\
			0.5 1.0\\
			0.75 1.0\\
			1 1.0\\
		};
		\addplot [color=red, mark=triangle]
		  table[row sep=crcr]{%
		0 1.0\\
		0.25 1.018\\
		0.5 0.929\\
		0.75 1.268\\
		1 0 \\
		};
		\addlegendentry{$q=2$}

		\addplot [color=blue, mark=asterisk]
		  table[row sep=crcr]{%
		  0 1.0\\
			0.25 1.000\\
			0.5 0.992\\
			0.75 1.088\\
			1 0\\
		};
		\addlegendentry{$q=1.2$}
		\end{axis}
\end{tikzpicture}\caption{$L^q$-best approximation to $u \equiv 1$ by a piecewise linear function $u_h$ satisfying $u_h(0)=1$ \\and $u_h(1)=0$ on a uniform mesh consisting of four elements with $q=2$ and $q=1.2$.}\label{fig:motivating_example}
\end{figure}

 Fig.\ \ref{fig:motivating_example} shows the $L^q$ best approximation of $u\equiv 1$ by a piecewise linear function $u_h$ satisfying $u_h(0)=1$ and $u_h(1)=0$ on a uniform mesh consisting of four elements with $q=2$ and $q=1.2$. We can see that in both cases over- and undershoots are present in the approximation, but that the magnitude of these oscillations is significantly smaller for $q=1.2$. This example illustrates the phenomenon of reducing oscillations in the approximation as $q \rightarrow 1$ that we shall  investigate in this article.

\subsection{Problem Statement}
We consider a subdivision $\varOmega_h$ of the domain $\varOmega = (0,1)^d$, $d = 1,2$ into $n$ disjoint open simplicial elements (i.e., subintervals when $d=1$ and triangles when $d=2$) $\kappa_i$, $i=1, \dots, n$, such that $\bar{\varOmega} = \bigcup_{i=1}^n \bar{\kappa}_i$ and define  $U_h$ to be the standard finite element space consisting of continuous piecewise linear polynomials on the mesh $\varOmega_h$.
Let $u \equiv 1$ and consider the following (constrained) best approximation problem:
\begin{subequations}
\begin{align}
  u_h = \argmin{v_h \in U_h} \|u -v_h\|_{L^q((0,1)^d)}
\end{align}
subject to
\begin{align}
  \begin{aligned}
    u_h(0) &= 1, \quad u_h(1) &= 0 && \text{ if }  d=1,\\
    u_h(0, \cdot) &=1, \quad u_h(1,\cdot) &= 0 && \text{ if } d=2.
  \end{aligned}
\end{align}
\label{eq:best_approx}
\end{subequations}
Note that the constraint can be removed by using a Dirichlet lift argument as commonly employed in the context of finite element methods and restricting the space $U_h$ to the functions that are zero on the part of the boundary where boundary conditions are employed.

In one dimension, instead of $u(x)\equiv 1$, we also consider the $L^q$-best approximation of $u(x) = \sgn(x)$ on
$(-1,1)$ by a continuous piecewise linear function $u_h$ satisfying $-u_h(-1) = u_h(1)=1$. We use this example to establish the link between our work and \cite{Saff1999}.

There is a related body of literature studying the $L^2$-projection onto finite element spaces, such as \cite{Bank2014,Douglas1975,Crouzeix1987}. These works are mostly concerned with the stability of the projection operator in subspaces (e.g., $L^q(\varOmega)$, $W^{1,q}(\varOmega)$, $H^1_0(\varOmega)$).

\subsection{Summary of Results}
The main result of this article consists of the precise analysis of very simple examples that illustrate the behaviour of $L^q$-best approximations of discontinuities by continuous piecewise linear polynomials on coarse meshes. We have furthermore included some numerical examples that confirm the theoretical analysis and illustrate how the observed behaviour in simple model examples applies to more general scenarios. In particular, we demonstrate that the over- and undershoots observed in $L^q$-best approximations for $1<q<\infty$ decrease as $q \rightarrow 1$. Whether these oscillations disappear entirely depends on the mesh used to define the underlying finite dimensional approximation space. In one dimension, Gibbs phenomena can be eliminated on uniform meshes both for a boundary discontinuity and a jump discontinuity present in  the interior of the domain. For non-uniform meshes it depends on the relative sizes of the elements. In two dimensions, we show that there exist uniform and structured meshes for which Gibbs phenomena are not eliminated. But, we also include examples of meshes in two dimensions on which the over- and undershoots vanish as $q\rightarrow 1$. Furthermore, we will illustrate that there exist infinitely many $L^1$-best approximations in certain cases which is due to the fact that $L^1(\varOmega)$ is not strictly convex.

\begin{figure}
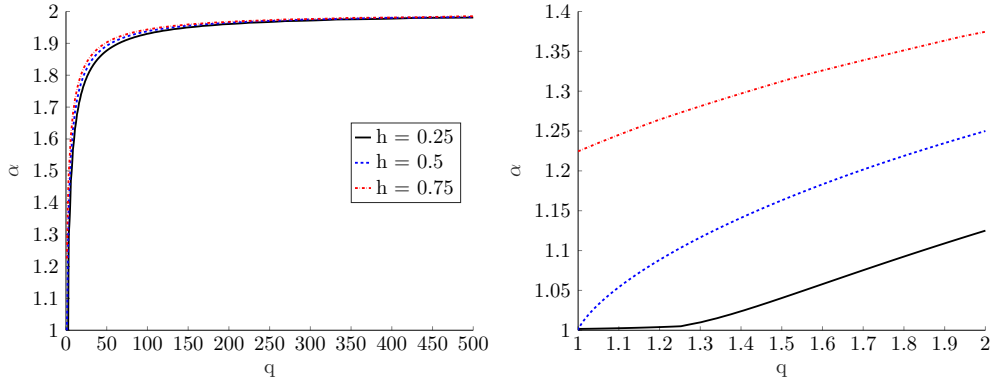

 \centering
 \input{overshoot_1d_1_500}
 \input{overshoot_1d_1_2}
 \caption{Values for $\alpha$ for different ranges of $q$ and three different choices of
 $h$.}\label{fig:alpha_one_d}
\end{figure}

The first example we consider is the approximation problem \eqref{eq:best_approx}
with $d = 1$.
The following theorem precisely characterises the $L^q$-best approximation for all $1 \leq q < \infty$ for any
two-element mesh on $(0, 1)$. Note that in this case,  the approximation problem only has one degree of freedom  due to the boundary conditions. Furthermore, we prove for an $N$-element mesh
that there exists an $L^1$-best approximation with no over- or undershoot if either a
grading-type mesh condition is satisfied or a stronger, simple element-size condition.
 The precise result is given below.
\begin{thm}[$L^q$-best approximation of a boundary discontinuity]\label{thm:lq_1d}
  \leavevmode
  \begin{enumerate}
   \item Consider the mesh given by the subdivision of $(0,1)$ into the two intervals $(0,1-h)$ and $(1-h, 1)$ with $h\in (0,1)$. For $1\leq q < \infty$, the solution of the approximation problem \eqref{eq:best_approx} with $u\equiv 1$ and $d =1$
    is given by a continuous piecewise linear polynomial $u_h$ that satisfies the boundary conditions and $u_h(1-h) =\alpha$, where $\alpha$ is defined as follows
  \begin{enumerate}
    \item If $q = 1$,
    \begin{align}
      \alpha  = \left\{ \begin{aligned} 1 && \text{ if } h \leq 0.5 \\
                        \sqrt{2h} && \text{ if } h > 0.5.\end{aligned} \right. \label{eq:alpha_l1_1d}
    \end{align}
    \item If $1<q<\infty$, then $\alpha > 1$ and
    \begin{align*}
      0 = -(1-h)\alpha^2 q(\alpha-1)^{q-1}-h(\alpha q+1)(\alpha-1)^q+h.
    \end{align*}
  \end{enumerate}
	\item  \label{lem:sufficient_1D_stronger}
	  Let the mesh be given by a subdivision of the interval $(0,1)$ into $N\geq 2$ intervals $(x_{i-1}, x_{i})$, $i = 1,\dots, N$, with $0 = x_0 < x_2 < \dots <x_{N-1} < x_{N} = 1$. The length $h_i$ of the $i$th subinterval is given
	  by $h_i = x_i - x_{i-1}$, $i = 1, \dots N$.
	  Define
	  \begin{align}
	    \vartheta_N&:=0,\\
	    \vartheta_i^2&:=\frac{1}{2}\left(1-\left(2(1-\vartheta_{i+1})^2-1\right)\frac{h_{i+1}}{h_i}\right), i= N-1, \dots, 1,\label{eq:def_vartheta}\\
	    M&:= \max \left(\{0\}\cup\left\{ i \in \left\{1, \dots , N-1\right\}\,:\, \vartheta_i \geq 1-\frac{1}{\sqrt{2}}\right\}\right).
	  \end{align}
	  Then
	\begin{align}
	  h_i \geq (2(1-\vartheta_{i+1})^2-1)h_{i+1}, \qquad \text{ for  } i = M, M+1, \dots, N-1,\label{eq:sufficient_cond_1d}
	\end{align}
	   is a sufficient condition for the existence of an $L^1$-best approximation $u_h$ of $u \equiv 1$ with $u_h(0) = 1$ and $u_h(1) = 0$ satisfying
	  $u_h(x_i) = 1$ for all $i = 1, \dots, N-1$. Hence, $u_h$ contains no over- or undershoots.
  \item Let the mesh be given by a subdivision of the interval $(0,1)$ into $N\geq 2$ intervals $(x_{i-1}, x_{i})$, $i = 1,\dots, N$, with $0 = x_0 < x_2 < \dots <x_{N-1} < x_{N} = 1$. The length $h_i$ of the $i$th subinterval is given
  by $h_i = x_i - x_{i-1}$, $i = 1, \dots N$. Then  $h_N \leq \min_{i = 1, \dots, N-1} h_i$ is a sufficient condition for the existence of an $L^1$-best approximation $u_h$ of $u \equiv 1$ with $u_h(0) = 1$ and $u_h(1) = 0$ satisfying
  $u_h(x_i) = 1$ for all $i = 1, \dots, N-1$.
\end{enumerate}
\end{thm}
\begin{rem}
Note that in the second part of Theorem \ref{thm:lq_1d}, condition \eqref{eq:sufficient_cond_1d} essentially states that elements cannot be too small compared to their neighbouring element closer to the discontinuity. Furthermore,  there are no conditions on the size of the elements contained in $(0, x_{M-1})$ if $M>0$. Moreover, it is always possible to ensure $M>0$  by selecting $h_M$ sufficiently large in comparison to $h_{M+1}$ such that $\vartheta_M> 1-1/\sqrt{2}$. This means that the mesh can be designed in such a way that it is allowed to be arbitrary away from the discontinuity without leading to oscillations. This observation is particularly useful if more than one discontinuity is to be approximated.
\end{rem}
\begin{rem}\label{rem:refinement}
  With very similar arguments as in the proof of the final part of Theorem \ref{thm:lq_1d}, it is easy to see that if $h_N>h_{N-1}$, but $h_{N-1} \leq h_i$ for all $i=1, \dots, N-2$, then every $L^1$-best approximation must contain over- or undershoots. Moreover, there exists an $L^1$-best approximation with overshoot only at the node $x_{N-1}$ and no further over- or undershoots, i.e., $u_h(x_i) = 1$ for
  $i = 1, \dots, N-2$ and
  \begin{align*}
    u_h(x_{N-1}) = \sqrt{\frac{2h_N}{h_N+h_{N-1}}}.
  \end{align*}
  The value at $u_h(x_{N-1})$ follows from the first part of the theorem by rescaling the interval.
\end{rem}
Fig.\ \ref{fig:alpha_one_d} shows $\alpha$ specified in Theorem \ref{thm:lq_1d} for two different ranges of $q$ and
three different choices of $h$. The plot
shows that $\alpha<2$ for all $1\leq q < \infty$ and that $\alpha$ decreases as $q\rightarrow 1$ for
all three choices of $h$. Furthermore, we can see that the behaviour as $q \rightarrow \infty$ is
very similar for all choices of $h$, but that there are clear differences  as $q \rightarrow 1$.
For $h = 0.25$ and $h = 0.5$, $\alpha$ approaches $1$ as $q \rightarrow 1$, hence the overshoot vanishes as $q\rightarrow 1$, whereas for
$h = 0.75$ it approaches $\sqrt{2h} \approx 1.2247$, hence the overshoot does not vanish. This is consistent with the results
obtained for the $L^1$-best approximation, cf., \eqref{eq:alpha_l1_1d}.

 In Section \ref{sec:examples_1d} we include examples of two three-element meshes violating the sufficient condition in part three of Theorem \ref{thm:lq_1d} such that one of the meshes satisfies \eqref{eq:sufficient_cond_1d}, whereas the other mesh violates this condition as well. We will demonstrate that for the latter mesh the overshoot does indeed not vanish entirely as $q\rightarrow 1$.

The second best approximation problem  we  analyse is the best approximation of $u(x)=\sgn(x)$ on $(-1,1)$ on a mesh consisting of exactly four elements that is symmetric with respect to $x = 0$. The main difference to the result in part one of Theorem \ref{thm:lq_1d} is that there exists a whole family of best approximations if $q = 1$. For $q>1$, we observe the same behaviour as before.

\begin{thm}[$L^q$-best approximation of a jump discontinuity]\label{thm:discontinuity}
  Consider the mesh given by the subdivision of $(-1,1)$ into the four intervals
  $(-1,-h)$, $(-h, 0)$, $(0,h)$ and $(h,1)$ with $h \in (0,1)$.
  For $1 \leq q < \infty$, the $L^q$-best approximation  of $u = \sgn(x)$ on $(-1,1)$ by a continuous piecewise linear function $u_h$ on the above mesh such that $-u_h(-1) = u_h(1) = 1$ can be characterised as follows.
  \begin{enumerate}
    \item If $h \leq 0.5$, there exists an $L^1$-best approximation for any $\beta \in [-1,1]$ such that $u_h(0) = \beta $ and $-u_h(-h) = u_h(h) = 1$. Conversely, any $L^1$-best approximation satisfies $-u_h(-h) = u_h(h) = 1$ and $u_h(0) \in [-1,1]$.
    \item If $h>0.5$, then for any $\beta \in [-1,1]$, $u_h$ satisfies
    \begin{subequations}
    \begin{align*}
      u_h(-h) &= \alpha := -\sqrt{2h}-\beta(\sqrt{2h}-1),\\
      u_h(0) & = \beta,\\
      u_h(h) & = \gamma := \sqrt{2h}-\beta(\sqrt{2h}-1),
    \end{align*}
  \end{subequations}
  defining an $L^1$-best approximation. Conversely, any $L^1$-best approximation satisfies $u_h(0) \in [-1,1]$ and is of the above form.
    \item The unique $L^q$-best approximation with $1<q<\infty$ is given by
    $-u_h(-h) = u_h(h)= \alpha$ and $u_h(0) = 0$, where $\alpha$ satisfies
    \begin{align*}
      0 = -(1-h)\alpha^2 q(\alpha-1)^{q-1}-h(\alpha q+1)(\alpha-1)^q+h
    \end{align*}
    and $\alpha >1$.
    \item In the limit $q \rightarrow 1$ the $L^q$-best approximation converges to the $L^1$-best approximation with $u_h(0) = 0$ for any $h\in (0,1)$. The corresponding $L^1$-best approximation is anti-symmetric and
    satisfies
    \begin{align*}
      -u_h(-h) = u_h(h) = \left\{ \begin{aligned}
      1 && \text{ if } h\leq 0.5,\\
      \sqrt{2h} && \text{ if } h > 0.5.
    \end{aligned}\right.
    \end{align*}
  \end{enumerate}
\end{thm}

We again observe that the presence of over- and undershoots in the $L^1$-best approximation depends on the choice of mesh. Furthermore, there exists a whole family of $L^1$-best approximations in this case which is possible since $L^1$ is not strictly convex and therefore minimizers are not necessarily unique. We recover uniqueness if we define the minimizer as the limit as $q \rightarrow 1$ of the $L^q$-minimizer. Moreover, it follows from the proof of Theorem \ref{thm:discontinuity} that the $L^1$-best approximation is unique if the subdivision of the interval is no longer symmetric, as we will see in Section \ref{sec:1d_disc}.

In order to see how this result relates to the work in \cite{Saff1999}, it first has to be noted that there are two major differences between our investigation and \cite{Saff1999}:
\begin{enumerate}
 \item The interval in \cite{Saff1999} is subdivided into $2n$ subintervals of equal length.
 In contrast to this, we only consider the special case that $(-1,1)$ is subdivided into $4$ subintervals and instead allow the subdivision to be non-uniform but still symmetric with respect to the center of the interval.
 \item  We consider bounded domains with fixed boundary conditions,  which are relevant to finite element approximations, whereas the investigation
 in \cite{Saff1999} considers the limit $n\rightarrow \infty$ for the interval $[-nh,nh]$ (ergo
 essentially an infinite domain) with no boundary conditions.
\end{enumerate}
In \cite{Saff1999} it is shown that for a uniform subdivision of the interval $[-nh,nh]$, the
over- and undershoots disappear as $n\rightarrow \infty$ and $q\rightarrow 1$.
The last point in Theorem \ref{thm:discontinuity} shows that, on a fixed mesh, we recover the result that the over- and undershoots disappear as $q \rightarrow 1$ for $h \leq 0.5$, which includes the case of a uniform mesh. However, if  $h>0.5$,  the over- and undershoots do not disappear as $q \rightarrow 1$.

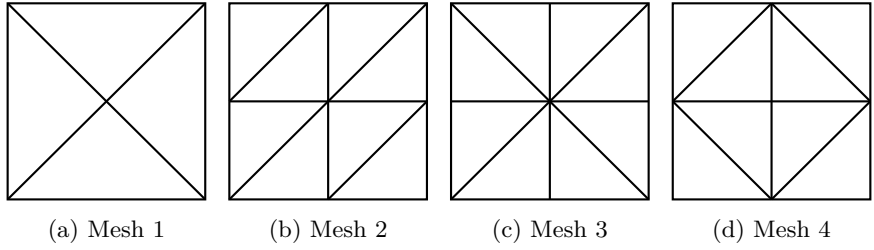
\begin{figure}[t]
  \centering
  \begin{subfigure}{2.8cm}
      \centering
    \begin{tikzpicture}[thick,scale=0.65, every node/.style={scale=0.9}]
    \draw [thick] (0,0) rectangle (4,4);
    \draw [thick](0,0) -- (4,4);
    \draw [thick](4,0) -- (0,4);
    \end{tikzpicture}
    \caption{Mesh 1}
  \end{subfigure}
\begin{subfigure}{2.8cm}
    \centering
\begin{tikzpicture}[thick,scale=0.65, every node/.style={scale=0.9}]
\draw [thick] (0,0) rectangle (4,4);
\draw [thick] (0,2) -- (4,2);
\draw [thick] (2,0) -- (2,4);
\draw [thick] (0,0) -- (4,4);
\draw [thick] (2,0) -- (4,2);
\draw [thick] (0,2) -- (2,4);
\end{tikzpicture}
\caption{Mesh 2}
\end{subfigure}
\begin{subfigure}{2.8cm}
    \centering
\begin{tikzpicture}[thick,scale=0.65, every node/.style={scale=0.9}]
\draw [thick] (0,0) rectangle (4,4);
\draw [thick] (0,2) -- (4,2);
\draw [thick] (2,0) -- (2,4);
\draw [thick] (0,0) -- (4,4);
\draw [thick] (0,4) -- (4,0);
\end{tikzpicture}
\caption{Mesh 3}
\end{subfigure}
\begin{subfigure}{2.8cm}
    \centering
\begin{tikzpicture}[thick,scale=0.65, every node/.style={scale=0.9}]
\draw [thick] (0,0) rectangle (4,4);
\draw [thick] (0,2) -- (4,2);
\draw [thick] (2,0) -- (2,4);
\draw [thick] (2,0) -- (4,2);
\draw [thick] (0,2) -- (2,4);
\draw [thick] (2,0) -- (0,2);
\draw [thick] (2,4) -- (4,2);
\end{tikzpicture}
\caption{Mesh 4}
\end{subfigure}
\caption{Four different meshes on $(0,1)^2$.}
\label{fig:meshes_2d_thm}
\end{figure}

The final theoretical result concerns the solution to \eqref{eq:best_approx} with $d = 2$ on the four meshes shown in Fig.\ \ref{fig:meshes_2d_thm}.  Note that the discrete space $U_h$ has only one degree of freedom on Mesh 1, corresponding to the value at the midpoint, and $U_h$ has three degrees of freedom on the other meshes, corresponding to the values at the three nodes on the line $x = 0.5$. For the first mesh, we analyse the $L^q$-best approximation for all $1 \leq q< \infty$ and show that the solution contains an overshoot that does not disappear as $q \rightarrow 1$.
For Mesh 2, we show that any $L^1$-best approximation must contain over- or undershoots and characterise an $L^1$-best approximation. Furthermore, we prove that there exists an $L^1$-best approximation on Meshes 3 and 4 without over- or undershoots. Finally, we demonstrate numerically in Section \ref{sec:2d_l1_rest_numerics} that the $L^q$-best approximation on Meshes 2, 3 and 4 indeed approaches the $L^1$-best approximation characterised in the theorem below.

\begin{thm}[$L^q$-best Approximation of a Boundary Discontinuity in Two Dimensions]\label{thm:best_approx_2d}
  \leavevmode
  \begin{enumerate}
    \item The unique solution to \eqref{eq:best_approx} on Mesh 1 for $q=1$ is given defined by $u_h(0.5,0.5) = \alpha$, where $\alpha$ satisfies
    \begin{align*}
      \alpha >1 &&\text{ and } &&0 = 2 \alpha^3- 5\alpha +2,
    \end{align*}
		hence $\alpha \approx 1.3200$.
    \item The unique solution to \eqref{eq:best_approx} on Mesh 1 for $1<q<\infty$ is defined by $u_h(0.5,0.5) = \alpha$, where $\alpha$ satisfies
    \begin{align*}
      \alpha>1 &&\text{ and } && 0 = (\alpha-1)^{q-1}\left[4\alpha^3 q +4(1-q)\alpha^2+(q-6)\alpha+2\right]-\alpha(q+4)+2.
    \end{align*}
		\item If $q>1$, the $L^q$-best approximation to \eqref{eq:best_approx} contains over- or undershoots on all four meshes.
    \item If $q=1$, there exists a solution to \eqref{eq:best_approx} on Mesh 2 such that $u_h(0.5, 1) = u_h(0.5,0.5) = 1$ and $u_h(0.5,0) = \alpha$, where $\alpha$ satisfies
    \begin{align*}
      \alpha > 1 && \text{ and }&& 0 = -3\alpha^3 + 8 \alpha -4,
    \end{align*}
		hence $\alpha \approx 1.2723$.
    Furthermore, $u_h(0.5, 1) = u_h(0.5,0.5) = u_h(0.5,0) = 1$ does not define an $L^1$-best approximation.
    \item If $q=1$, there exists a solution to \eqref{eq:best_approx} on Meshes 3 and  4 such that $u_h(0.5, 1) = u_h(0.5,0.5) = u_h(0.5,0) = 1$.
  \end{enumerate}
\end{thm}
Theorem \ref{thm:best_approx_2d} shows that, on Meshes 1 and 2, the $L^q$-best approximations exhibit an overshoot for all $q$, including $q=1$, while on Meshes 3 and 4 the $L^1$-best approximation does not contain any over- or undershoots.

\begin{figure}
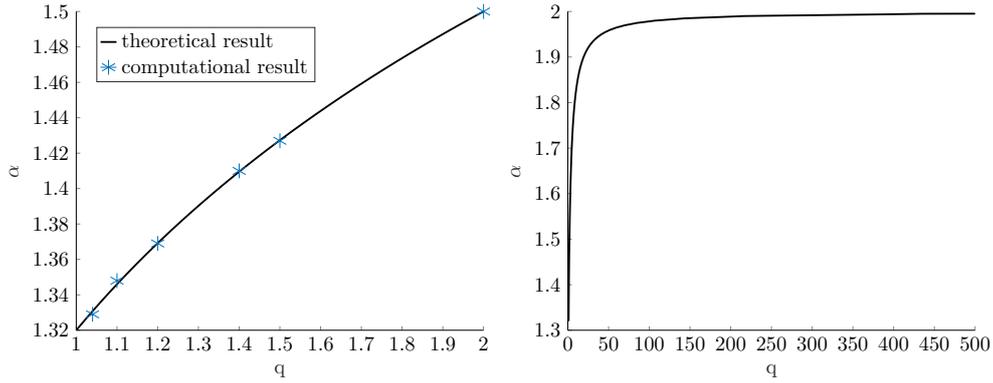

  \centering
  \input{crisscross_q1_2}
  \input{crisscross_q1_500}
  \caption{Values for $\alpha$ for different ranges of $q$ on Mesh 1.}\label{fig:alpha_crisscross}
\end{figure}

Fig.\ \ref{fig:alpha_crisscross} shows the parameter $\alpha$ defining the $L^q$-best approximation on Mesh 1 for two different ranges of $q$. The
plot shows that $\alpha <2$ for all $q$ and that $\alpha$ decreases as $q \rightarrow 1$,
where it approaches $1.32$. This is  consistent with the result in Theorem \ref{thm:best_approx_2d} obtained for
the $L^1$-best approximation. To confirm the theoretical results, we have also
determined the $L^q$-best approximation numerically by implementing \eqref{eq:Lqbest}
using FEniCS \citep{Alnes2015}. The solution to the resulting non-linear system can be approximated using a Newton iteration if $q$ is sufficiently close to $2$. Note that this solver is not robust in $q$ and stalls or diverges for $q$ close to $1$ and for $q \gg 2$. The left plot in Fig.\ \ref{fig:alpha_crisscross} shows numerically
determined approximations of $\alpha$ for selected values of $q$ which confirm
the theoretical results.

We also include further numerical experiments in Section \ref{sec:examples} illustrating that the observations remain the same if $u$ is a more general smooth function and that the over- and undershoots cannot be eliminated by refining the mesh.

\subsection{Outline of the Paper}
The remainder of this article is organised as follows: in Section \ref{sec:char_best} we describe a characterisation of the $L^q$-best approximation of a function in a finite dimensional subspace  that we will use to prove our theoretical results;
Sections \ref{sec:1d_layer}, \ref{sec:1d_disc} and \ref{sec:2d_layer} contain the proofs of Theorems \ref{thm:lq_1d}, \ref{thm:discontinuity} and \ref{thm:best_approx_2d}, respectively.
We conclude with several numerical examples in Section \ref{sec:examples} illustrating the effect of mesh refinement in one and two dimensions and showing the behaviour of the $L^q$-best approximation as $q \rightarrow 1$  in one dimension, as well as on structured and unstructured meshes in two dimensions.

\section{Characterisation of Best \texorpdfstring{$\bm{L^q}$}{$L^q$}-Approximation}\label{sec:char_best}
In this section we describe a characterisation of best-approximation in Banach spaces
and more specifically the Lebesgue spaces $L^{q}(\varOmega)$, $1\leq q < \infty$. This
 characterisation  will be used in the remainder of this article
to determine the best $L^q$-approximation for specific examples.

If $U$ is a Banach space and $f$ a function $f:U\rightarrow \mathbb{R}$, the
subdifferential $\partial f(u)$ of $f$ at a point $u \in U$ is defined as the set
\begin{align*}
  \partial f(u) : = \left\{u^{\prime} \in U^{\prime} \,:\, f(w)-f(u) \geq \langle u^{\prime}, w-u \rangle_{U^{\prime},U}, \, \forall w \in U \right\}
\end{align*}
If $f$ is G\^ateaux differentiable, the subdifferential is single
valued and agrees with the G\^ateaux derivative.
We now quote the following theorem, cf., \cite[Theorem~1.1]{Singer1970}.
\begin{thm}[Characterisation of best approximation]\label{thm:best_approximation}
   Let $U$ be a Banach space, $U_h \subset U$
  a closed subspace and $u\in U$. The following statements are equivalent:
\begin{enumerate}
  \item $\displaystyle u_h = \mathrm{arg}\, \min_{w_h \in U_h}\|u-w_h\|_U$.
  \item There exists a functional $r^{\prime} \in \partial \left( \| \cdot \|_U\right)(u-u_h)$ which annihilates
  $U_h$, i.e., \[\langle r^{\prime}, w_h\rangle_{U^{\prime}, U} = 0 \qquad \text{ for all }w_h \in U_h.\]
\end{enumerate}
\end{thm}
\begin{rem}
  The subdifferential $\partial \left(\|\cdot \|_U\right)(\cdot)$ can be characterised
  as follows, cf., e.g., \cite[Chapter~1,~Proposition~3.4]{Cioranescu1990}. For any $w \in U$,
  \begin{align}
    \partial \left(\|\cdot \|_U\right) (w) := \left\{
    \begin{aligned}
      \{ w^{\prime} \in U^{\prime} \,: \, \langle w^{\prime},  w \rangle_{U^{\prime}, U}  = \|w\|_U,\, \|w^{\prime}\|_{U^{\prime}} = 1\} &&\text{if } w \neq 0,\\
      \{ w^{\prime} \in U^{\prime} \,: \, \|w^{\prime}\|_{U^{\prime}} = 1\}&& \text{if } w = 0.
    \end{aligned}
    \right.\label{eq:subdiff_norm}
  \end{align}
  This characterisation allows us to translate the above formulation of
Theorem \ref{thm:best_approximation} directly into the formulation
  found in \cite{Singer1970}. In \cite{Muga2017} the same theorem is stated
   in terms of the so-called duality mapping, which can also be easily translated
   into the above formulation.
\end{rem}
First we will use Theorem \ref{thm:best_approximation} to characterise
best approximants in subspaces of $L^q(\varOmega)$, $1<q< \infty$. To this end, we determine the subdifferential  $\partial \left(\| \cdot\|_{L^q(\varOmega)}\right)(w)$ for an arbitrary $w \in L^q(\varOmega)$ and $1 < q< \infty$. Note that in this case
the norm is G\^ateaux differentiable; indeed, we can compute for $w\not\equiv 0$:
\begin{align*}
  \partial \left(\|\cdot \|_{L^q(\varOmega)}\right) (w)(v) = \left.\frac{\dm}{\dm t}\left(\int_{\varOmega} |w+tv|^q \dm \bm x\right)^{\frac{1}{q}}\right|_{t=0}
  = \|w\|_{L^q({\varOmega})}^{1-q}\int_{\varOmega} \sgn(w)|w|^{q-1}v \dm \bm x,
\end{align*}
where
\begin{align*}
  \mathrm{sgn}(w(x)) = \left\{
  \begin{aligned}
    -1 && \text{ if } w(\bm x) < 0,\\
    1 && \text{ if } w(\bm x) >0, \\
    0 && \text{ if } w(\bm x)= 0,
  \end{aligned}\right.
\end{align*}
hence  $\partial \left(\| \cdot\|_{L^q(\varOmega)}\right)(w) = \|w\|_{\Lp{q}} \sgn (w)|w|^{q-1}$ by the canonical identification of an element in the dual space of $\Lp{q}$ with a function in $\Lp{q'}$, where $1 = 1/q+1/q'$.

The following Corollary is an immediate consequence of this by setting $w = u-u_h$.
\begin{cor}\label{cor:best_lq}
Let $U:=\Lp{q}$ and $U_h \subset U$ a closed subspace. The function $u_h \in U_h$ is an $L^q$-best approximation of $u$ if and only if
\begin{align}
  \begin{aligned}
		\int_{\varOmega}\sgn(u-u_h)|u-u_h|^{q-1} v_h \dm \bm x = 0 \qquad \forall v_h \in U_h.
\end{aligned}\label{eq:best_lq}
\end{align}
\end{cor}

Next we will use \eqref{eq:subdiff_norm} to characterise
best approximations in subspaces of $L^1(\varOmega)$.  Note that in this
case the subdifferential $\partial \left(\|\cdot \|_{L^1(\varOmega)}\right) (w)$ is in general not single valued for an arbitrary $w \in \Lp{1}$ .
From \eqref{eq:subdiff_norm}, we deduce that
\begin{align*}
  \partial \left(\|\cdot \|_{L^1(\varOmega)}\right) (w)(v) = \int_{\varOmega}\psi v \mathrm{d}\bm x,
\end{align*}
where $\psi \in L^{\infty}(\varOmega)$ with the following properties
\begin{enumerate}
  \item $\displaystyle\|\psi\|_{L^{\infty}(\varOmega)} = 1$.
  \item $\displaystyle \int_{\varOmega}\psi w \mathrm{d}\bm x = \|w\|_{L^1(\varOmega)}$.
\end{enumerate}
It is easy to see that any $\psi$ such that $\psi = \sgn(w)$ if $w \neq 0$ and $|\psi| \leq 1$ almost everywhere satisfies the above conditions.
Conversely, the first property implies $|\psi(x)| \leq 1$ almost everywhere and
the second property implies that $\psi(x) = 1$ almost everywhere on
 $\{u(x) > 0 \}$ and $\psi(x) = -1$ almost everywhere on
 $\{u(x) < 0 \}$ since
 \begin{align*}
  \|w\|_{L^1(\varOmega)} =\int_{\varOmega} |w| \dm \bm x
  = \int_{\varOmega}\psi w \dm \bm x = \int_{\varOmega \cap \{w(\bm x) > 0 \}}\psi |w| \dm \bm x
  - \int_{\varOmega \cap \{w(\bm x) < 0 \}}\psi |w| \dm \bm x.
 \end{align*}
It is important to note, that the only condition on $\psi$ on the
set $\{w(x) = 0\}$ is that $|\psi|\leq 1$ almost everywhere.
The following Corollary characterising $L^1$-best approximations is a direct consequence of this by setting $w = u-u_h$.

\begin{cor}\label{cor:best_approx_l1}
Let $U:=\Lp{q}$ and $U_h \subset U$ a closed subspace. The function $u_h \in U_h$ is an $L^1$-best approximation of $u$ if and only if there exists a function
  $\psi_0 \in L^{\infty}(\varOmega \cap \{u(\bm x) = u_h(\bm x)\})$, $|\psi_0|\leq 1$, almost everywhere, such that for all $v_h \in U_h$
\begin{align}
  0 = \int_{\varOmega\cap \{u(\bm x) > u_h(\bm x) \}}v_h \mathrm{d}\bm x - \int_{\varOmega\cap \{u(\bm x) < u_h(\bm x) \}}v_h \mathrm{d}\bm x
  +\int_{\varOmega\cap \{u(\bm x) = u_h(\bm x) \}}\psi_0 v_h \mathrm{d}\bm x, \label{eq:best_approx_l1}
\end{align}
or, equivalently, for all $v_h \in U_h$,
\begin{align*}
	0 = \int_{\varOmega} \psi v_h \dm \bm x,
\end{align*}
where $\psi = \sgn(u-u_h)$ on $\{u(\bm x) \neq u_h(\bm x)\}$ and $\psi = \psi_0$ on $\{u(\bm x) = u_h(\bm x)\}$.
\end{cor}
Note that in the case that $u$ and $u_h$ only agree on a set of measure zero,
the choice of $\psi_0\in[-1,1]$ becomes irrelevant.

\section{Best Approximation of a Boundary Discontinuity in One Dimension}\label{sec:1d_layer}

In this section we consider the best approximation problem \eqref{eq:best_approx} in one dimension and provide a proof of Theorem \ref{thm:lq_1d}. We split this into three parts: Sections \ref{sec:1d_l1_bd} and \ref{sec:1d_lq_bd} contain the proof of the first part of the theorem, where the former addresses the case $q=1$ and the latter the case $1<q<\infty$;
Section \ref{sec:sufficient_gen_mesh} contains the proof of the second  and third part of the theorem.

In the first part of Theorem \ref{thm:lq_1d}, we consider the mesh consisting of the two subintervals $(0,1-h)$ and $(1-h, 1)$. The best approximation $u_h$ of $u \equiv 1$ by a continuous piecewise linear function satisfying the boundary conditions $u_h(0) = 1$ and $u_h(1) = 0$ is determined entirely by the value it takes at the point $x=1-h$.

Thereby, we can write $u_h = \varphi_0 + \alpha \varphi_1$, where $\alpha$ is to be determined and
 \begin{align*}
   \varphi_0 &= \left\{
   \begin{aligned}
     \frac{(1-h) -x}{1-h} && \text{ in } [0,1-h],\\
     0 && \text{ else,}
   \end{aligned}\right. &
   \varphi_1 &= \left\{
   \begin{aligned}
     \frac{x}{1-h} && \text{ in } [0,1-h],\\
     \frac{1-x}{h} && \text{ in } [1-h,1].
   \end{aligned}\right.
 \end{align*}

 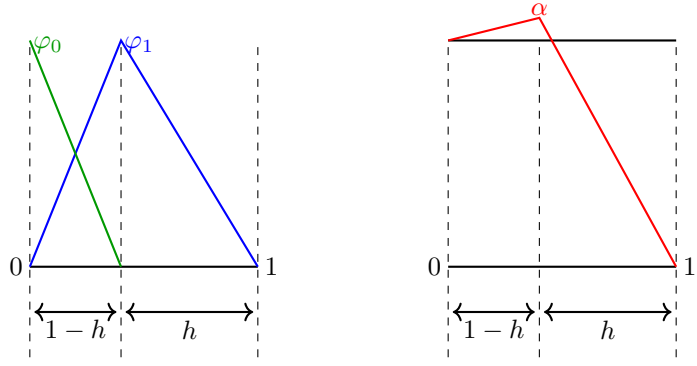
\begin{figure}
   \centering
 \begin{tikzpicture}[scale = 0.6]
 \draw [thick](0,0)--(5,0);
 \draw [blue,  thick] (0,0)--(2,5)--(5,0);
 \draw [green!60!black,  thick] (0,5)--(2,0);
 \draw [<->,thick](0.1,-1)--(1.9,-1);
 \draw [<->,thick](2.1,-1)--(4.9,-1);

 \node at (1,-1.4) {$1-h$};
 \node at (3.5,-1.4) {$h$};

 \node at (-0.3,0){$0$};
 \node at (5.3,0){$1$};
 \node[green!60!black] at (0.4,4.9) {$\varphi_0$};
 \node[blue] at (2.4,4.9) {$\varphi_1$};

 \draw [dashed] (0,-2)--(0,5);
 \draw [dashed] (2,-2)--(2,5);
 \draw [dashed] (5,-2)--(5,5);

 \end{tikzpicture}
 \hspace{1.5cm}
 \begin{tikzpicture}[scale = 0.6]
 \draw [thick](0,0)--(5,0);
 \draw [thick] (0,5)--(5,5);
 \draw [red, thick] (0,5)--(2,5.5)--(5,0);
 \draw [<->,thick](0.1,-1)--(1.9,-1);
 \draw [<->,thick](2.1,-1)--(4.9,-1);

 \node at (1,-1.4) {$1-h$};
 \node at (3.5,-1.4) {$h$};

 \node at (-0.3,0){$0$};
 \node at (5.3,0){$1$};
 \node[red] at (2,5.7) {$\alpha$};

 \draw [dashed] (0,-2)--(0,5);
 \draw [dashed] (2,-2)--(2,5);
 \draw [dashed] (5,-2)--(5,5);

 \end{tikzpicture}
 \caption{Left: Hat functions $\varphi_0$ and $\varphi_1$. Right: Approximation $u_h$ with $\alpha >1$.}
 \label{fig:1d_layer_basis}
 \end{figure}
 Fig.\ \ref{fig:1d_layer_basis} shows the two functions $\varphi_0$ and $\varphi_1$
 as well as an approximation $u_h$ of $u \equiv 1$ with $\alpha >1$.
 To eliminate the constraint by introducing a Dirichlet lift in the best approximation problem \eqref{eq:best_approx}, we could define the subspace $U_h$ as the span of $\varphi_1$  and redefine $u = 1 - \varphi_0$ and $u_h = \alpha \varphi_1$. Note, however, that $u-u_h$ remains the same.
 The main consequence of this observation is, that \eqref{eq:best_lq} and \eqref{eq:best_approx_l1} do not have to be satisfied for $w_h = \varphi_0$ due to the boundary condition constraint.

\subsection{$L^1$-Best Approximation}\label{sec:1d_l1_bd}
In this section we give a proof of the first part of Theorem \ref{thm:lq_1d} for the case $q=1$. More precisely, we show that
  the $L^1$-best approximation of $u \equiv 1$ on $(0,1)$ by a continuous piecewise linear function $u_h$ satisfying the boundary conditions $u_h(0)=0$ and $u_h(1)=1$ is given by
  $u_h = \varphi_0 + \alpha \varphi_1$, where
  \begin{align}
    \alpha  = \left\{ \begin{aligned} 1 && \text{ if } h \leq 0.5 ,\\
                      \sqrt{2h} && \text{ if } h > 0.5.\end{aligned} \right.
											\label{eq:alpha_l1_1d_proof}
  \end{align}

\begin{proof}
Using the characterisation of the best approximation given in
Corollary \eqref{cor:best_approx_l1},
 to determine the best approximation, we  distinguish between two cases:
 \begin{enumerate}
   \item The set $\{x \in (0,1) \,: \,(u-u_h)(x) = 0\}$ has measure zero. (For continuous piecewise linear functions this
   set has to consist of a finite number of points).  This means  that $\psi = \mathrm{sgn}(u-u_h)$ everywhere except on a set of measure zero and is thus uniquely defined almost everywhere.
   \item The set $\{x \in (0,1) \,: \,(u-u_h)(x) = 0\}$ has positive measure, i.e., the set contains
   an interval of positive length. This means that $\psi$ is not uniquely defined on a set with positive measure.
 \end{enumerate}

 Starting with the second case, we observe that this can only be true if $\alpha = 1$ and
 thus $u = u_h$ in $(0,1-h)$.
 In this case $\psi = \mathrm{sgn} (u-u_h) = 1$ in $(1-h,1)$ and we compute
 \begin{align*}
   \int_{1-h}^1\psi \varphi_1\mathrm{d}x = \int_{1-h}^1\frac{1-x}{h}\mathrm{d}x = \frac{h}{2}.
 \end{align*}
Assuming $-h/(1-h) \in [-1,1]$, we can choose $ \psi_0 = -h/(1-h)$ in $(0,1-h)$ and obtain
 \begin{align*}
 \int_{0}^{1-h} \psi \varphi_1\mathrm{d}x = -\frac{h}{1-h}\underbrace{\int_{0}^{1-h}(x+1)\mathrm{d}x}_{=(1-h)/2} = -\frac{h}{1-h}\frac{1-h}{2} = -\frac{h}{2}.
\end{align*}
Obviously, adding both integrals yields zero for this choice of $\psi_0$.
Therefore, by Corollary \ref{cor:best_approx_l1},  $\alpha =1 $ is an $L^1$-best approximation provided that
$-h/(1-h) \in [-1,1]$. This is the case if and only if $h \leq 1-h \Leftrightarrow h\leq 1/2$, hence proving \eqref{eq:alpha_l1_1d_proof} for $h\leq 1/2$.
If on the other hand $ h> 1/2 \Leftrightarrow 1-h <h$, $\alpha = 1$ does not yield an $L^1$-best approximation. Indeed,
\begin{align*}
\int_0^1 \psi\varphi_1 \dm x  \geq \frac{h}{2}-\frac{1-h}{2}>0,
\end{align*}
since $\psi = 1$ in $(0,1-h)$ and $\psi \geq -1$ in $(1-h,1)$.

Next we consider the first case; this implies that $\alpha \neq 1$ and $\psi = \sgn(u-u_h)$ almost everywhere in $(0,1)$. If $\alpha <1$, we have
$u-u_h>0$  everywhere in $(0,1)$ and since $\varphi_0 > 0$ everywhere in $(0,1)$, we also have
\begin{align*}
  \int_{0}^1\sgn(u-u_h) \varphi_1 \mathrm{d}x > 0.
\end{align*}
Therefore, $u_h = \varphi_0+\alpha \varphi_1$ cannot be an $L^1$-best approximation if
$\alpha < 1$. This leaves the case $\alpha > 1$. Now, $u-u_h < 0$ in $(0, \frac{\alpha-h}{\alpha})$ and $u-u_h > 0$
in $(\frac{\alpha-h}{\alpha},1)$; we compute
\begin{align*}
  \begin{aligned}
  \int_{0}^1 \mathrm{sgn}(u-u_h) \varphi_1\mathrm{d}x
  &= -\int_{0}^{1-h} \frac{x}{1-h} \mathrm{d}x - \int_{1-h}^{\frac{\alpha-h}{\alpha}}\frac{1-x}{h}\mathrm{d}x
  +\int_{\frac{\alpha-h}{\alpha}}^1\frac{1-x}{h}\mathrm{d}x\\
  &= -\frac{1-h}{2} -\frac{h(\alpha^2-1)}{2\alpha^2}+\frac{h}{2\alpha^2} = \frac{2h-\alpha^2}{2\alpha^2}.
\end{aligned}
\end{align*}
This integral becomes $0$ for $\alpha = \sqrt{2h}$. Note that this only yields an $L^1$-best
approximation if $h>1/2$.
Indeed, if $h\leq1/2$, then $\alpha = \sqrt{2h} \leq 1$, but we have assumed $\alpha>1$.

We have therefore shown that $\alpha =1$ is the only $L^1$ best approximation if $h\leq 1/2$ and
$\alpha = \sqrt{2h}$ is the only $L^1$-best approximation if $h>1/2$.
\end{proof}

\subsection{$L^q$-Best Approximation}\label{sec:1d_lq_bd}
 In this section we give a proof of the first part of Theorem \ref{thm:lq_1d} for the case $1<q<\infty$. More precisely, we show that
    the $L^q$-best approximation of $u \equiv 1$ on $(0,1)$ by a continuous piecewise linear function $u_h$ satisfying the boundary conditions $u_h(0)=0$ and $u_h(1)=1$ is given by
    $u_h = \varphi_0 + \alpha \varphi_1$, where $\alpha > 1$ and
    \begin{align*}
      0 = -(1-h)\alpha^2 q(\alpha-1)^{q-1}-h(\alpha q+1)(\alpha-1)^q+h.
    \end{align*}

\begin{proof}
Corollary \ref{cor:best_lq} implies that we seek $\alpha$ such that
\begin{align*}
  \int_{0}^1 \sgn(u-u_h) |u-u_h|^{q-1} \varphi_1 \dm \bm x = 0.
\end{align*}
Again, we have to split the integral on each element
into the parts where $u-u_h>0$, $u-u_h <0$ and $u-u_h = 0$. We consider three cases
\begin{align*}
  \text{ (a) } \alpha < 1,&&
  \text{ (b) } \alpha = 1,&&
  \text{ (c) } \alpha >1.
\end{align*}
If $\alpha < 1$, we have $u-u_h>0$ everywhere in $(0,1)$ and thus
both $\sgn(u-u_h)|u-u_h|^{q-1}>0$ and $\varphi_1>0$ in $(0,1)$. Therefore, we have
\begin{align*}
  \int_{0}^1\sgn(u-u_h) |u-u_h|^{q-1} \varphi_1\dm  x >0,
\end{align*}
hence $\alpha <1$ is not possible.
If $\alpha = 1$, we have $u-u_h = 0$ in $(0,1-h)$ and $u-u_h > 0$ in $(1-h,1)$. Thus,
\begin{align*}
  \int_{0}^1 \sgn(u-u_h) |u-u_h|^{q-1} \varphi_1 \dm  x
  =\int_{1-h}^1 \sgn(u-u_h) |u-u_h|^{q-1}\varphi_1\dm x > 0,
\end{align*}
hence $\alpha =1$ is not possible.
We can therefore assume $\alpha >1$. In this case $u-u_h < 0$ in $\left(0,\frac{\alpha-h}{\alpha}\right)$ and
 $u-u_h> 0$ in $\left(\frac{\alpha-h}{\alpha}, 1 \right)$. We compute
\begin{align*}
  \begin{aligned}
  \int_{0}^1&\sgn{(u-u_h)}|u-u_h|^{q-1}\varphi_1 \dm x \\&= -\int_{0}^{1-h}\left(\frac{(\alpha -1)x}{1-h}\right)^{q-1}\frac{x}{1-h} \dm x \\
&- \int_{1-h}^{\frac{\alpha-h}{\alpha}}\left(\frac{\alpha-h}{h} - \frac{\alpha x}{h}\right)^{q-1} \frac{1-x}{h} \dm x
+ \int_{\frac{\alpha-h}{\alpha}}^1 \left(\frac{\alpha x}{h}-\frac{\alpha-h}{h}\right)^{q-1} \frac{1-x}{h} \dm x\\
  &= -\frac{(1-h)(\alpha-1)^{q-1}}{q+1}-\frac{h(\alpha q +1)(\alpha-1)^q}{\alpha^2 q (q+1)}
  + \frac{h}{\alpha^2 q(q+1)}\\
  &= \frac{-(1-h)\alpha^2 q(\alpha-1)^{q-1}-h(\alpha q+1)(\alpha-1)^q+h}{\alpha^2 q (q+1)}.
\end{aligned}
\end{align*}
Hence, the $L^q$-best approximation can be determined by  finding $\alpha >1$ satisfying
\begin{align*}
  0 = -(1-h)\alpha^2 q(\alpha-1)^{q-1}-h(\alpha q+1)(\alpha-1)^q+h,
\end{align*}
existence of which is guaranteed since the $L^q$-best approximation always exists.
\end{proof}

\subsection{Sufficient Conditions on General Meshes}\label{sec:sufficient_gen_mesh}

 In this section we provide a proof of the second and third parts of Theorem \ref{thm:lq_1d}.
  To this end, let the mesh be given by a subdivision of the interval $(0,1)$ into $N\geq 2$ subintervals $(x_{i-1}, x_{i})$, $i = 1,\dots, N$, with $0 = x_0 < x_2 < \dots <x_{N-1} < x_{N} = 1$. The length $h_i$ of the $i$th subinterval is given
  by $h_i = x_i - x_{i-1}$, $i = 1, \dots N$.
	In order to prove the second part of Theorem \ref{thm:lq_1d}, we show that the following conditions are sufficient for the existence of an $L^1$-best approximation $u_h$ of $u \equiv 1$ with $u_h(0) = 1$ and $u_h(1) = 0$ satisfying
 $u_h(x_i) = 1$ for all $i = 1, \dots, N-1$:
 \begin{align}
   h_i \geq (2(1-\vartheta_{i+1})^2-1)h_{i+1} \qquad \text{ for  } i = M, M+1, \dots, N-1,\label{eq:sufficient_cond_1d_proof}
 \end{align}
 where
  \begin{align*}
    \vartheta_N&:=0,\\
    \vartheta_i^2&:=\frac{1}{2}\left(1-\left(2(1-\vartheta_{i+1})^2-1\right)\frac{h_{i+1}}{h_i}\right), i= N-1, \dots, 1,\\
    M&:= \max \left(\{0\}\cup\left\{ i \in \left\{1, \dots , N-1\right\}\,:\, \vartheta_i \geq 1-\frac{1}{\sqrt{2}}\right\}\right).
  \end{align*}
We then show that the much simpler condition  $h_N \leq \min_{i = 1, \dots, N-1} h_i$ implies \eqref{eq:sufficient_cond_1d_proof} which proves the third part of Theorem \ref{thm:lq_1d}.

\begin{proof}
  Define $u_h$ such that $u_h(0) = 1$, $u_h(1) = 0$ and
 $u_h(x_i) = 1$ for all $i = 1, \dots, N-1$. Furthermore, denote by $\varphi_i$ the hat function that is $1$ at $x_i$ and $0$ at
 $x_j$, $j\neq i$, for $i = 1, \dots, N-1$.

For $\alpha \in (0,1]$, define $\psi_{\alpha}(x)$ as follows:
\begin{align*}
  \psi_{\alpha}(x) = \left\{ \begin{aligned}
      &(-1)^{N-i+1}&& x \in (x_{i-1},x_{i-1}+\vartheta_i h_{i}), &&\text{ for all } i =M+1, \dots N,\\
      &(-1)^{N-i} && x \in (x_{i-1}+\vartheta_i h_{i}, x_{i}) &&\text{ for all } i =M+1, \dots N,\\
      &(-1)^{N-M+1}\alpha && x \in (x_{M-1},x_{M-1}+\tilde{\vartheta}_M h_{M})&&\text{ if } M>0,\\
      &(-1)^{N-M} && x \in (x_{M-1}+\tilde{\vartheta}_M h_{M}, x_{M}) &&\text{ if } M>0,\\
      &0 && \text{otherwise },\\
\end{aligned}
  \right.
\end{align*}
where
\begin{align*}
  \tilde{\vartheta}_M^2 = \frac{2}{1+\alpha}\vartheta_M^2 =  \frac{1}{\alpha+1}\left(1-\left(2(1-\vartheta_{M+1})^2-1\right)\frac{h_{M+1}}{h_M}\right).
\end{align*}
We claim that there exists $\alpha \in (0,1]$ such that
\begin{align}
  \int_0^1 \psi_{\alpha}(x)\varphi_i(x) dx = 0 \qquad \forall i = 1, \dots N-1,\label{eq:psi_alpha_l1_cond}
\end{align}
if \eqref{eq:sufficient_cond_1d_proof} is satisfied.
First we show that $\psi_{\alpha}$ is well defined. For this we require $\vartheta_i$ to be well defined, i.e., we require $\vartheta_i\in [0,1]$, for all $i = M, \dots, N-1$. This is trivially true
for $\vartheta_N$. Otherwise,  for $i \geq M$,
\begin{align*}
  \vartheta_i^2 \geq 0 \Leftrightarrow h_i \geq (2(1-\vartheta_{i+1})^2-1)h_{i+1}.
\end{align*}
Furthermore, $\vartheta_{i+1} < 1-1/\sqrt{2}$ for $i\geq M$ by definition of $M$ and therefore
\begin{align*}
  (1-\vartheta_{i+1})^2 > \frac{1}{2} \Rightarrow \left(2(1-\vartheta_{i+1})^2-1\right) >0 \Rightarrow \vartheta_i^2 < \frac{1}{2} \Rightarrow  \vartheta_i < 1.
\end{align*}
Next, note that $\vartheta_N = 0$ implies $\psi_{\alpha}(x) = 1  = \mathrm{sgn}(u-u_h)$ on $(x_{N-1}, x_N)$ and that $\|\psi_{\alpha}\|_{L^{\infty}((0,1))}=1$.
Thus, proving \eqref{eq:psi_alpha_l1_cond} immediately implies that
$u_h$ is an $L^1$-best approximation of $u \equiv 1$.

For all $i = 1, \dots, N-1$ it holds that
\begin{align}
  \int_0^1\psi_{\alpha}(x)\varphi_i ( x)dx = \int_{x_{i-1}}^{x_{i}}\psi_{\alpha}(x)\varphi_i( x) dx +\int_{x_{i}}^{x_{i+1}}\psi_{\alpha}(x)\varphi_i( x) dx.\label{eq:phi_i_zero}
\end{align}
With this in mind, we now consider  \eqref{eq:psi_alpha_l1_cond} for $i>M$. In this case
\begin{align}
  \int_{x_{i-1}}^{x_{i}}\psi_{\alpha}(x)\varphi_i ( x)dx &= (-1)^{N-i+1}\frac{h_i}{2}\left(2\vartheta_i^2-1\right),\\
  \int_{x_{i}}^{x_{i+1}}\psi_{\alpha}(x)\varphi_i( x) dx &= (-1)^{N-i+2}\frac{h_{i+1}}{2}\left(1-2(1-\vartheta_{i+1})^2\right),\label{eq:phii_hi+1}
\end{align}
Hence \eqref{eq:phi_i_zero} becomes zero if and only if
\begin{align*}
0=  \frac{h_i}{2}\left(2\vartheta_i^2-1\right) -\frac{h_{i+1}}{2}\left(1-2(1-\vartheta_{i+1})^2\right) \Leftrightarrow \vartheta_i^2&=\frac{1}{2}\left(1-\left(2(1-\vartheta_{i+1})^2-1\right)\frac{h_{i+1}}{h_i}\right),
\end{align*}
which is the definition of $\vartheta_i$.

Next, we consider $i=M$. We have already established $\vartheta_M^2 \leq 1/2$. With $\alpha>0$, this implies $\tilde{\vartheta}_M < 1$. Furthermore, \eqref{eq:phii_hi+1} still holds with $i=M$, whereas we obtain
\begin{align*}
\int_{x_{M-1}}^{x_{M}}\psi_{\alpha}\varphi_i dx  =(-1)^{N-M+1}((\alpha + 1)\vartheta^2_M - 1)\frac{h_{M}}{2}.
\end{align*}
Hence \eqref{eq:phi_i_zero} becomes zero for $i = M$ if and only if
\begin{align*}
&0=  \frac{h_M}{2}\left((1+\alpha)\tilde{\vartheta}_M^2-1\right) -\frac{h_{M+1}}{2}\left(2(1-\vartheta_{M+1})^2-1\right)\\ \Leftrightarrow \qquad&\tilde{\vartheta}_M^2=\frac{1}{1+\alpha}\left(1-\left(2(1-\vartheta_{M+1})^2-1\right)\frac{h_{M+1}}{h_M}\right),
\end{align*}
which is the definition of $\tilde{\vartheta}_M$.
Finally, we have to show that there exists $\alpha \in (0,1]$ such that
\begin{align*}
  0 &= \int_{0}^1\psi_{\alpha}\varphi_{M-1} dx =\underbrace{\int_{x_{M-2}}^{x_{M-1}} \psi_{\alpha}\varphi_{M-1}}_{= 0, \text{ since }\psi_{\alpha}(x) = 0 \text{ for } x< x_{M-2}} +\int_{x_{M-1}}^{x_{M}}\psi \varphi_{M-1} dx \\
  &= (-1)^{N-M+1}(\alpha-(\alpha + 1)(1-\tilde{\vartheta}_M)^2) \frac{h_{M}}{2}.
\end{align*}
The last expression becomes zero if and only if
\begin{align*}
  (1-\tilde{\vartheta}_M)^2 = \frac{\alpha}{1+\alpha} \Leftrightarrow \tilde{\vartheta}_M= \frac{1}{\sqrt{1+\alpha}}\left( \sqrt{1+\alpha}-\sqrt{\alpha}\right).
\end{align*}
Hence, we need $\alpha$ such that
\begin{align}
  \left( \sqrt{1+\alpha}-\sqrt{\alpha}\right) = \sqrt{\left(1-\left(2(1-\vartheta_{M+1})^2-1\right)\frac{h_{M+1}}{h_M}\right)}.\label{eq:alpha_theta}
\end{align}
For $\alpha >0$, $g(\alpha) =\sqrt{1+\alpha}-\sqrt{\alpha}$, is a strictly decreasing function of $\alpha$ and thus bijectively maps $(0,1]$ onto $[\sqrt{2}-1, 1)$. The equation \eqref{eq:alpha_theta} therefore has a unique solution $\alpha \in (0,1]$ if and only if
\begin{align*}
  \sqrt{\left(1-\left(2(1-\vartheta_{M+1})^2-1\right)\frac{h_{M+1}}{h_M}\right)} \in [\sqrt{2}-1, 1)
   \quad\Leftrightarrow \quad  & \vartheta_M \in \left[1-\frac{1}{\sqrt{2}}, \frac{1}{\sqrt{2}}\right).
\end{align*}
By the definition of $M$, we have $\vartheta_{M} >1-\frac{1}{\sqrt{2}}$ and  $\vartheta_{M+1} < 1-\frac{1}{\sqrt{2}} \Rightarrow \vartheta_M < \frac{1}{\sqrt{2}}$.
This shows that \eqref{eq:sufficient_cond_1d_proof} is indeed a sufficient condition for the existence of an $L^1$-best approximation such that $u_h(x_i) = 1$ for all $i=1, \dots, N-1$, which finishes the proof of part two of Theorem \ref{thm:lq_1d}.

  Next, we show that $h_N \leq \min_{i = 1, \dots, N-1} h_i$ implies that \eqref{eq:sufficient_cond_1d_proof} is satisfied. Let $\vartheta_i$ for $i =1, \dots N-1$ be defined as in Lemma \ref{lem:sufficient_1D_stronger}.
We first show that if $\vartheta_i$ is real and $\vartheta_i < 1-1/\sqrt{2}$ for some $k \geq 1$ and all $i \geq k$, then the following holds:
\begin{align}
  \frac{h_i}{2} \left(2(1-\vartheta_i)^2 -1\right)\leq \frac{h_{i+1}}{2}\left(2(1-\vartheta_{i+1})^2 -1\right) \quad \text{ for all } i \geq k.\label{eq:decreasing_integral}
\end{align}
Indeed, by the definition of $\vartheta_i$, we have
\begin{align*}
  \frac{h_{i+1}}{2}\left(2(1-\vartheta_{i+1})^2 -1\right) = \frac{h_i}{2}(1-2\vartheta_i).
\end{align*}
Hence, \eqref{eq:decreasing_integral} is equivalent to
\begin{align*}
  \frac{h_i}{2} \left(2(1-\vartheta_i)^2 -1\right)\leq \frac{h_i}{2}(1-2\vartheta_i^2) \Leftrightarrow \vartheta_i(\vartheta_i - 1) \leq 0,
\end{align*}
which is true since $\vartheta_i \in (0,1)$ for $i \geq k$.
Next, note that if we apply \eqref{eq:decreasing_integral} recursively, we obtain
\begin{align}
  \frac{h_{k}}{2}\left(2(1-\vartheta_{k})^2 - 1\right)\leq \frac{h_N}{2}.\label{eq:bound_hn}
\end{align}

In order to prove that \eqref{eq:sufficient_cond_1d_proof} is satisfied,
first note that \eqref{eq:sufficient_cond_1d_proof} reduces to $h_{N-1} \geq h_N$ for $i=N-1$ since $\vartheta_N = 0$.
Now assume for the sake of contradiction that for some $j>0$ and all $i =j,  j+1, \dots N-1$, $\vartheta_i < 1-1/\sqrt{2}$ and assume \eqref{eq:sufficient_cond_1d_proof} holds for $i = j+1, \dots N-1$ but not for $i=j$.  In this case, we have
\begin{align*}
  h_N \leq h_j < \left(2(1-\vartheta_j)^2 -1\right)h_{j+1} \stackrel{\text{ by }\eqref{eq:bound_hn}}{\leq} h_N,
\end{align*}
which is a contradiction.
\end{proof}
\begin{rem}
  There is a more direct proof of the third part of Theorem \ref{thm:lq_1d} by defining $\psi(x) \equiv (-1)^jh_N/h_{N-j}$ on $(x_{N-j-1}, x_{N-j})$ for $j = 1, \dots, N-1$ and showing that
  \begin{align*}
    \int_{0}^1\psi\varphi_i dx = 0, \qquad \text{ for } i=1,\dots, N-1.
  \end{align*}
  However, the proof given above shows that the condition $h_N \leq \min_{i = 1, \dots, N-1} h_i$  is always stronger than \eqref{eq:sufficient_cond_1d_proof}.
\end{rem}

\section{Over- and Undershoots at Jump Discontinuities}\label{sec:1d_disc}

In this section we consider the $L^q$-best approximation of $u(x) = \sgn(x)$ in $(-1,1)$ as an example of  a jump discontinuity
in the interior of the domain and provide a proof of Theorem \ref{thm:discontinuity}. We split this into three parts: in Section \ref{sec:h_l_05_disc}, we consider the first part of Theorem \ref{thm:discontinuity}, i.e., we consider the case where the $L^1$-best approximation does not exhibit Gibbs phenomena.
In Section \ref{sec:h_g_05_disc}, we consider the case where the $L^1$-best approximation does exhibit Gibbs phenomena (part two of Theorem \ref{thm:discontinuity}); finally, in Section \ref{sec:q_disc} we consider the
$L^q$-best approximation for $1<q<\infty$ and the limit as $q \rightarrow 1$ (parts three and four of Theorem \ref{thm:discontinuity}).

\begin{figure}
  \centering
\begin{tikzpicture}[scale = 0.6]
\draw [thick](0,0)--(10,0);
\draw [red!50!yellow,  thick] (0,0)--(2,5)--(5,0);
\draw [green!60!black,  thick] (2,0)--(5,5)--(8,0);
\draw [blue, thick] (5,0)--(8,5)--(10,0);
\draw [<->,thick](0.1,-1)--(1.9,-1);
\draw [<->,thick](2.1,-1)--(4.9,-1);
\draw [<->,thick](5.1,-1)--(7.9,-1);
\draw [<->,thick](8.1,-1)--(9.9,-1);
\node at (1,-1.4) {$1-h$};
\node at (3.5,-1.4) {$h$};
\node at (6.5,-1.4) {$h$};
\node at (9,-1.4) {$1-h$};
\node at (-0.8,0){$-1$};
\node at (10.8,0){$1$};
\node[red!50!yellow] at (2.4,4.9) {$\varphi_1$};
\node[green!60!black] at (5.4,4.9) {$\varphi_2$};
\node[blue] at (8.3,4.9) {$\varphi_3$};
\draw [dashed] (0,-2)--(0,5);
\draw [dashed] (2,-2)--(2,5);
\draw [dashed] (5,-2)--(5,5);
\draw [dashed] (8,-2)--(8,5);
\draw [dashed] (10,-2)--(10,5);

\end{tikzpicture}
\caption{Basis for $U_h$}
\label{fig:jump_discontinuity_basis}
\end{figure}
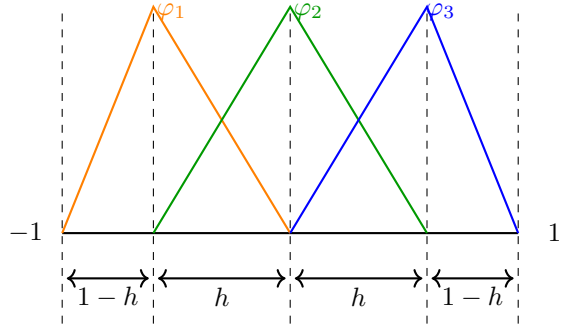

For this example, $\varOmega = (-1,1)$ and $u(x) = \sgn(x) \in L^1(\varOmega)$.
We seek an $L^q$-best approximation of this function by a
continuous piecewise linear function on the mesh consisting of $(-1,-h)$, $(-h,0)$, $(0,h)$
and $(h,1)$. We fix the boundary conditions at $-1$ and $1$, i.e.,
$u_h(1) = u(1)=1$ and $u_h(-1) = u(-1)=-1$.
The finite dimensional approximation space $U_h$ is given by the span of the hat functions $\varphi_i$, $i=1,2,3$,
depicted in Fig.\ \ref{fig:jump_discontinuity_basis}.

The condition for $u_h$ to be
an $L^1$ best-approximation in Corollary \ref{cor:best_approx_l1} can be written as follows: there exists $\psi:(-1,1)
\rightarrow [-1,1]$ such that
\begin{align}
  \int_{-1}^1\psi \varphi_i \dm x = 0, \qquad \text{ for } i=1,2,3,\label{eq:best_approx_gibbs}
\end{align}
and $\psi(x) = \sgn((u-u_h)(x))$ on  $\{x \in (-1,1) \,:\, u(x)\neq u_h(x)\}$.

\subsection{$L^1$-Best Approximation without Over- or Undershoots}\label{sec:h_l_05_disc}
In this section we provide a proof of the first part of Theorem \ref{thm:discontinuity}.
More precisely, if $h\leq 0.5$,  a continuous piecewise linear function $u_h$ on the mesh shown in Fig.\ \ref{fig:jump_discontinuity_basis} such that $-u_h(-1) = u_h(1) = 1$ is an $L^1$-best approximation of $u(x) = \sgn(x)$ if and only if   $u_h(0) = \beta $, with $\beta \in [-1,1]$ arbitrary,  and $-u_h(-h) = u_h(h) = 1$.

\begin{proof}
We first show that $u_h$ must satisfy $-u_h(-h) = u_h(h) = 1$. For the sake of contradiction assume $u_h(-h) < -1$. In this case $u-u_h >0$ on $(-1,-h)$ and $u-u_h >0$ in $(-h, -h+\delta)$ with $\delta >0$. Thus, $\psi = \sgn(u-u_h)$ almost everywhere and we obtain
\begin{align*}
  \int_{-1}^{1}\sgn(u-u_h) \varphi_1 \dm x = \underbrace{\int_{-1}^{-h}\varphi_1 \dm x}_{=(h-1)/2}+\int_{-h}^{0} \sgn(u-u_h) \varphi_1 \dm x >
    \frac{h-1}{2}- \underbrace{\int_{-h}^0 \varphi_1 \dm x}_{= h/2} \geq 0,
\end{align*}
since $h\leq 0.5 \Rightarrow h-1 \geq h$. This is a contradiction since the condition
\eqref{eq:best_approx_gibbs} is violated. Note that we can use the same argument with
opposite sign for $u_h(-h)>-1$. Hence, we have proven that  $u_h(-h)=-1$. Due to the symmetry of the problem, it can easily be seen that the argument for $u_h(h) = 1$ is identical.

If $-u_h(-h) = u_h(h) =1$, then $u_h(0) \not\in [-1,1]$ implies $\psi = \sgn(u-u_h) = -\sgn(u_h(0))$ on $(-h,h)$ and hence
\begin{align*}
  \int_{-1}^{1} \psi\varphi_2 \dm x =   \int_{-h}^{h} \psi\varphi_2 \dm x
  = -\sgn(u_h(0))h \neq 0.
\end{align*}
Again, \eqref{eq:best_approx_gibbs} is violated which implies $u_h(0) \in [-1,1]$.

Next we establish, that $-u_h(-h) = u_h(h)=1$ and $u_h(0)=\beta$ is indeed an $L^1$-best approximation of $\sgn(x)$ for any $\beta \in[-1,1]$. To this end, we distinguish the following three cases:
\begin{align*}
  \text{(a) } u_h(0) = 1, && \text{(b) }u_h(0) \in (-1,1), && \text{(c) } u_h(0) = -1.
\end{align*}
All three cases are shown in Fig.\ \ref{fig:jump_discontinuity}a-c, where
$u_h(0) =0$ was chosen as an example for the second case.
\begin{figure}
  \centering
\begin{subfigure}{3.2cm}
    \centering
\begin{tikzpicture}[scale = 0.3]
\draw [thick](0,0)--(10,0);
\draw [thick] (0,-5)--(5,-5);
\draw [thick] (5,5)--(10,5);
\draw [dashed] (0,-5.5)--(0,5.5);
\draw [dashed] (3,-5.5)--(3,5.5);
\draw [dashed] (5,-5.5)--(5,5.5);
\draw [dashed] (7,-5.5)--(7,5.5);
\draw [dashed] (10,-5.5)--(10,5.5);
\draw [red, very thick] (0,-5)--(3,-5)--(5,5)--(10,5);
\end{tikzpicture}\caption{$u_h(0) =1$}
\end{subfigure}
\begin{subfigure}{3.2cm}
    \centering
\begin{tikzpicture}[scale = 0.3]
\draw [thick](0,0)--(10,0);

\draw [thick] (0,-5)--(5,-5);
\draw [thick] (5,5)--(10,5);
\draw [dashed] (0,-5.5)--(0,5.5);
\draw [dashed] (3,-5.5)--(3,5.5);
\draw [dashed] (5,-5.5)--(5,5.5);
\draw [dashed] (7,-5.5)--(7,5.5);
\draw [dashed] (10,-5.5)--(10,5.5);
\draw [red, very thick] (0,-5)--(3,-5)--(7,5)--(10,5);
\end{tikzpicture}
\caption{$u_h(0) =0$}\label{fig:jump_l1_0}
\end{subfigure}
\hspace{0.5em}
\begin{subfigure}{3.2cm}
    \centering
\begin{tikzpicture}[scale = 0.3]
\draw [thick](0,0)--(10,0);

\draw [thick] (0,-5)--(5,-5);
\draw [thick] (5,5)--(10,5);
\draw [dashed] (0,-5.5)--(0,5.5);
\draw [dashed] (3,-5.5)--(3,5.5);
\draw [dashed] (5,-5.5)--(5,5.5);
\draw [dashed] (7,-5.5)--(7,5.5);
\draw [dashed] (10,-5.5)--(10,5.5);
\draw [red, very thick] (0,-5)--(5,-5)--(7,5)--(10,5);
\end{tikzpicture}
\caption{$u_h(0) =-1$}
\end{subfigure}
\begin{subfigure}{3.2cm}
    \centering
\begin{tikzpicture}[scale = 0.3]
\draw [thick](0,0)--(10,0);
\draw [thick] (0,-5)--(5,-5);
\draw [thick] (5,5)--(10,5);
\draw [dashed] (0,-5.5)--(0,5.5);
\draw [dashed] (2,-5.5)--(2,5.5);
\draw [dashed] (5,-5.5)--(5,5.5);
\draw [dashed] (8,-5.5)--(8,5.5);
\draw [dashed] (10,-5.5)--(10,5.5);
\draw [red, very thick] (0,-5)--(2,-5.5);
\draw [red, very thick] (2,-5.5)--(8,5.5);
\draw [red, very thick] (8,5.5)--(10, 5);
\end{tikzpicture}
\caption{$h>1/2$}\label{fig:jump_l1_1_1}
\end{subfigure}
\caption{$L^1$-best approximation of a jump discontinuity}
\label{fig:jump_discontinuity}
\end{figure}
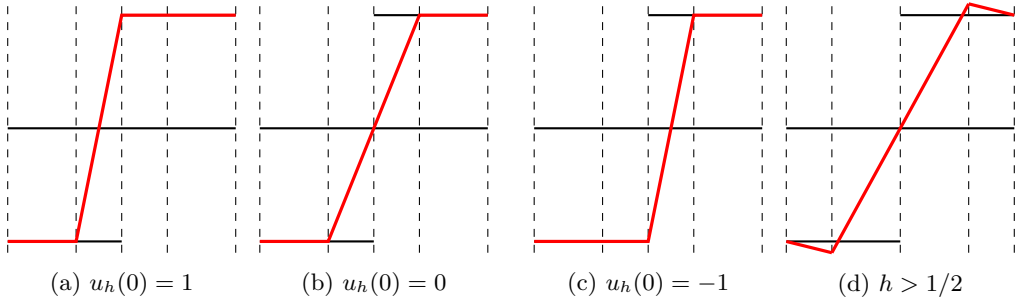
In the first case, $u - u_h = 0$ in $(-1,-h)$ and $u-u_h <0$ in $(-h,0)$. We compute
\begin{align*}
  \int_{-1}^{1}\psi \varphi_1 \dm x = \int_{-1}^{-h}\psi_0\varphi_1 \dm x
    - \underbrace{\int_{-h}^0 \varphi_1 \dm x.}_{= h/2}
\end{align*}
This integral is equal to zero for $\psi_0 \equiv h/(1-h)$ on $(-1,-h)$ which is a valid choice provided $h\leq1-h$.
Next, we compute
\begin{align*}
  \int_{-1}^{1} \psi\varphi_2 \dm x = -\underbrace{\int_{-h}^0 \varphi_2 \dm x}_{h/2}
  + \int_{0}^h \psi_0\varphi_2 \dm x.
\end{align*}
Note that $u-u_h = 0$ in $[0,h]$; it is easy to see that this integral becomes zero if and only if $\psi_0 \equiv 1$ in
$[0,h]$.
The remaining integral then becomes
\begin{align*}
  \int_{-1}^{1} \psi\varphi_3 \dm x = \underbrace{\int_{0}^h\varphi_3 \dm x}_{h/2}
  + \int_{h}^1\psi_0\varphi_3 \dm x.
\end{align*}
This integral is zero for $\psi_0 \equiv -h/(1-h)$ in $[h,1]$ which is again a valid choice if $h\leq 1-h$. Therefore, we have that  case (a) defines an $L^1$-best approximation
if $h\leq 1-h \Leftrightarrow h \leq 1/2$.

Due to the symmetry of the problem, it is easy to see that the third case, i.e.,
$u_h(0)=-1$ also defines an $L^1$ best approximation if  $h \leq 1-h$.
This leaves the second case. As in the first case, $u-u_h = 0$ in $[-1,-h]$ and
$u-u_h < 0$ in $[-h,0]$. This implies that we have to require $h\leq 1-h$ in
order to find $\psi_0$ such that the integral involving $\varphi_1$ becomes zero.
Furthermore, note that -- again due to the symmetry of the problem -- the same
applies to the integral involving $\varphi_3$. Finally, since
$u-u_h < 0$ in $[-h,0]$ and $u-u_h > 0$ in $[0,h]$, we can compute
\begin{align*}
  \int_{-1}^{1}\psi\varphi_2 \dm x = -\int_{-h}^0 \varphi_2\dm x+ \int_{0}^h
  \varphi_2 \dm x = 0.
\end{align*}
Thus, any choice of $u_h(0)\in[-1,1]$ defines an $L^1$-best approximation if $h\leq 1/2$. This completes the proof.
\end{proof}
\begin{rem}
Note that we have shown that there is a whole family of $L^1$-best approximations
with no over- or undershoots
for this particular example if $h\leq 1/2$. The situation is quite different if we instead consider a non-symmetric subdivision of the
interval $(-1,1)$ into $(-1,-h_1)$, $(-h_1,0)$, $(0,h_2)$ and $(h_2,1)$ with
$h_1 \neq h_2$. The integral involving $\varphi_2$ then implies
that the case $-1<u_h(0)<1$ does not yield an $L^1$-best approximation; the case $u_h(0)=1$ is an $L^1$-best
approximation if and only if $h_1 < h_2 \leq 1/2$, and the case
$u_h(0)=-1$ is an $L^1$-best approximation if and only if $h_2 < h_1 \leq 1/2$.
\end{rem}
\subsection{$L^1$-Best Approximation with Over- and Undershoots}\label{sec:h_g_05_disc}
In this section we provide a proof of the second part of Theorem \ref{thm:discontinuity}. More precisely, if $h> 0.5$, we show that a continuous piecewise linear function $u_h$ on the mesh shown in Fig.\ \ref{fig:jump_discontinuity_basis} such that $-u_h(-1) = u_h(1) = 1$ is an $L^1$-best approximation of $u(x) = \sgn(x)$ if and only if
  \begin{subequations}
  \begin{align}
    u_h(-h) &= \alpha := -\sqrt{2h} -\beta(\sqrt{2h} -1)\label{eq:alpha_sgn}\\
    u_h(0) & = \beta,\\
    u_h(h) & = \gamma := \sqrt{2h} -\beta(\sqrt{2h} -1)\label{eq:gamma_sgn},
  \end{align}\label{eq:uh_h_g_05}
\end{subequations}
with $\beta \in [-1,1]$ arbitrary.

\begin{proof}
	We start by showing that there exists $\psi$ as in Corollary \ref{cor:best_approx_l1}, such that $u_h$ as defined in \eqref{eq:uh_h_g_05} satisfies \eqref{eq:best_approx_l1} with $v_h = \varphi_1$.
In order to determine $\psi(x)$, we have to find the points where $u_h$ intersects $u$. Note that $u_h$ cannot intersect
$u$ in $(-1,-h)$ or $(h,1)$ unless it is identical with $u$ in all of $(-1,-h)$ or $(h,1)$, respectively. In order to determine the intersections in $(-h,0)$ and
$(0,h)$, we write $u_h$ as defined in \eqref{eq:uh_h_g_05} as follows:
\begin{align*}
  u_h(x) = \left\{\begin{aligned}
        \frac{\beta - \alpha}{h}x+\beta, && x \in(-h,0),\\
        \frac{\gamma - \beta}{h}x + \beta, && x\in (0,h).
\end{aligned}\right.
\end{align*}
We start with finding the intersection in $(-h,0)$;
note that $\beta-\alpha = 0 \Leftrightarrow \beta = -1$ and $u-u_h = 0$ everywhere in $(-1,0)$ in this case. If we assume $\beta \neq -1$, we obtain in $(-h,0)$ and
\begin{align*}
  0 = u-u_h = -1-\frac{\beta - \alpha}{h}x-\beta \Leftrightarrow  x = \frac{h(1+\beta)}{\alpha - \beta} = -\frac{h}{\sqrt{2h}} =: -h \vartheta.
\end{align*}
If now $\beta > -1$, we have $\alpha < -1$ and thus $\sgn(u-u_h) = 1$ in $(-1,-h\vartheta)$ and $\sgn(u-u_h) = -1$ in $(-h\vartheta,0)$. Hence $\psi = \sgn(u-u_h)$ in $(-1,0)$ and we compute
\begin{align}
  \begin{aligned}
  \int_{-1}^1\psi\varphi_1 \dm x &= \int_{-1}^{-h}\varphi_1 \dm x +
  \int_{-h}^{-\vartheta h} \varphi_1 \dm x - \int_{-\vartheta h}^0 \varphi_1 \dm x\\ &=
  \frac{1-h}{2}+\left(\frac{h}{2}-\frac{\vartheta^2h}{2}\right)- \frac{\vartheta^2h}{2}
  = \frac{1-2\vartheta^2h }{2} =0.
\end{aligned}\label{eq:phi_1}
\end{align}
Note that the above integral is zero if and only if $\vartheta = 1/\sqrt{2h}$.
Conversely, if $\beta > -1$ we know that $\sgn(u-u_h) = -1$ in $(-h\tilde{\vartheta}, 0)$ for some $\tilde{\vartheta} \in (0,1]$. Since $1-h < h$, replacing $\vartheta$ by $\tilde{\vartheta}$ in \eqref{eq:phi_1} implies
$\tilde{\vartheta} = 1/\sqrt{2h} = \vartheta$, which shows that \eqref{eq:best_approx_l1} uniquely defines the intersection point as $-\vartheta h$ if $\alpha < -1$.
 Note, that $\alpha <-1$ is necessary if $\beta >-1$, because $u-u_h$ has to either change sign or become zero for $u_h$ to satisfy \eqref{eq:best_approx}.
The point of intersection of $u$ and $u_h$ at $-\vartheta h$ and the value of $\alpha = u_h(-h)$ uniquely defines $u_h(0) = \beta$ such that \eqref{eq:alpha_sgn} holds. Hence, it is also a necessary condition that $u_h$ satisfies \eqref{eq:alpha_sgn} in order for $u_h$ to be an $L^1$-best approximation if $\beta >-1$.
Note that in this case the computation for $\beta < -1$ is completely analogous with opposite signs and therefore  we have not yet established that $\beta \geq -1$ is a necessary condition.

Similarly, $\gamma - \beta = 0 \Leftrightarrow \beta = 1$ and $u-u_h = 0$ everywhere in $(0,1)$ in this case. If we assume $\beta \neq 1$, we obtain in $(0,h)$
\begin{align*}
  0 = u-u_h = 1 -\frac{\gamma - \beta}{h}x - \beta \Leftrightarrow x = \frac{h(\beta -1)}{\gamma - \beta} = h\vartheta.
\end{align*}

If now $\beta < 1$, we have $\gamma > 1$ and thus $\sgn(u-u_h)=1$ in $(0, h\vartheta)$ and $\sgn(u-u_h) = -1$ in $(h\vartheta, 1)$.
Due to the symmetry of the problem, the computation is up to the sign the same as for $\beta \neq -1$ and $\varphi_1$ and we obtain
\begin{align*}
  \int_{-1}^{1}\psi \varphi_3 \dm x = 0.
\end{align*}
Hence $u_h$ as defined in \eqref{eq:uh_h_g_05} satisfies \eqref{eq:best_approx_l1} with $v_h = \varphi_3$.
Conversely, $u_h(0) = \beta$ implies $u_h(h) = \gamma$ as defined in \eqref{eq:uh_h_g_05} by an analogous argument to the proof of the implication $u_h(0) = \beta \Rightarrow u_h(-h) = \alpha$.
Thus, $u_h(h) = \gamma$ as defined in \eqref{eq:gamma_sgn} is a necessary condition for $u_h$ to be an $L^1$-best approximation if $\beta <1$.
Again note that the computation for $\beta > 1$ is completely analogous with opposite signs  and therefore  we have not yet established that $\beta \leq 1$ is a necessary condition.

To complete the cases $\beta \neq 1$ and $\beta \neq -1$, only the integral involving $\varphi_2$ remains. If $\beta \in (-1,1)$, we have $\gamma >1$ and $\alpha <-1$.
 Therefore, $\sgn(u-u_h) = -1$ in $(-\vartheta h, 0)$ and $(\vartheta h, h)$ and $\sgn(u-u_h) = 1$ in $(-h, -\vartheta h)$ and $(0,\vartheta h)$. Using the symmetry of $\varphi_2$, we obtain
\begin{align*}
  \int_{-1}^{1}\psi\varphi_2 \dm x = \underbrace{\int_{-h}^{-\vartheta h} \varphi_2 \dm x - \int_{\vartheta h}^h \varphi_2 \dm x}_{= 0} + \underbrace{\int_{0}^{\vartheta h} \varphi_2 \dm x - \int_{-\vartheta h}^0 \varphi_2 \dm x}_{=0} = 0.
\end{align*}
Hence, $u_h$ as defined in \eqref{eq:uh_h_g_05} also satisfies \eqref{eq:best_approx_l1} with $v_h = \varphi_2$ and we have established that \eqref{eq:uh_h_g_05} defines an $L^1$-best approximation for $\beta >1$ and $\beta<-1$.
Conversely, we have also shown that $u_h$ must satisfy \eqref{eq:uh_h_g_05} in this case.

The only two remaining cases are $\beta =1$ and $\beta = -1$. If $\beta =1$, then $\gamma = 1$ and $\alpha < -1$. Thus, $\sgn(u-u_h)$ on $(-1,0)$ is the same as in the case $\beta \in (-1,1)$ and $u-u_h = 0$ in $(0,1)$.
We now simply have to determine a valid choice for $\psi(x)$ on $(0,1)$ such that all integrals in \eqref{eq:best_approx_gibbs} are zero. One possible choice is trivially given by simply choosing the same as in the case $\beta \in (-1,1)$.
Note that we have already established that $\alpha$ has to be of the form \eqref{eq:alpha_sgn} for any $\beta > -1$ including $\beta = 1$.
Furthermore, $\gamma =1$ is necessary if $\beta = 1$ since otherwise either $u-u_h>0$ in $(0,1)$ (if $\gamma <1$) or $u-u_h<0$ in $(0,1)$ (if $\gamma>1$) and \eqref{eq:best_approx_l1} is violated with $v_h =\varphi_3$.
Due to the symmetry of the problem, the case $\beta = -1$ is analogous.

We finish the proof by showing that $|\beta| \leq 1$ is necessary. If $\beta < -1$ and hence in particular $\beta <1$, we have already established that  $\gamma >1$.  Furthermore, $\alpha >-1$ again follows from the fact that $u-u_h$ has to change sign within $(-1,0)$. Therefore, $\sgn(u-u_h) = -1$ in $(-h, -\vartheta h)$ and $(\vartheta h, h)$ and $\sgn(u-u_h) = 1$ in $(-\vartheta h, \vartheta h)$. Thus, we obtain

\begin{align*}
  \int_{-1}^{1}\psi\varphi_2 \dm x &= -\int_{-h}^{-\vartheta h} \varphi_2 \dm x +\int_{-\vartheta h}^0 \varphi_2 \dm x - \int_{0}^{\vartheta h} \varphi_2 \dm x
  + \int_{\vartheta h}^h \varphi_2 \dm x   \\
  &= 2\left(\int_{\vartheta h}^h \varphi_2 \dm x - \int_{0}^{\vartheta h} \varphi_2 \dm x \right) = h\left(2(1-\vartheta)^2-1\right) = 0\\
  & \Leftrightarrow (1-\vartheta)^2 = \frac{1}{2} \Leftrightarrow \vartheta = 1- \frac{1}{\sqrt{2}}.
\end{align*}
On the other hand, we require $\vartheta = 1/\sqrt{2 h}$ for the other integrals to become zero. So,
\begin{align*}
  \frac{1}{\sqrt{2h}} = 1-\frac{1}{\sqrt{2}} \Leftrightarrow h = \left( \frac{1}{\sqrt{2}-1}\right)^2 > 1.
\end{align*}

This is a contradiction, since $h \in (0,1)$. For $\beta >1$, the sign of $u-u_h$ on $(-h, h)$ is exactly opposite compared to the case $\beta <-1$; therefore, it is easy to see that $\beta >1$ also leads to a contradiction.
\end{proof}
\subsection{$L^q$-Best Approximation}\label{sec:q_disc}
In this section, we prove the third and fourth part of Theorem \ref{thm:discontinuity}.
More precisely, we show that a continuous piecewise linear function $u_h$ on the mesh shown in Fig.\ \ref{fig:jump_discontinuity_basis} such that $-u_h(-1) = u_h(1) = 1$ is an $L^q$-best approximation of $u(x) = \sgn(x)$ for $1<q<\infty$ if and only if
    $-u_h(-h) = u_h(h)= \alpha$ and $u_h = 0$, where $\alpha$ satisfies
    \begin{align*}
      0 = -(1-h)\alpha^2 q(\alpha-1)^{q-1}-h(\alpha q+1)(\alpha-1)^q+h
    \end{align*}
    and $\alpha >1$.
Furthermore, we show that in the limit $q \rightarrow 1$ the $L^q$-best approximation converges to the $L^1$-best approximation as defined in \eqref{eq:uh_h_g_05} with $\beta = 0$, for any $h \in (0,1)$, i.e.,  the corresponding $L^1$-best approximation is anti-symmetric and
    satisfies
    \begin{align*}
      -u_h(-h) = u_h(h) = \left\{ \begin{aligned}
      1 && \text{ if } h\leq 0.5,\\
      \sqrt{2h} && \text{ if } h > 0.5.
    \end{aligned}\right.
    \end{align*}

\begin{proof}
 We  use the characterisation of the $L^q$-best approximation in Corollary \ref{cor:best_lq}. Due to the uniqueness of the $L^q$-best approximation
for $1<q<\infty$ and the symmetry of the problem, we may assume that the $L^q$-best approximation is an odd function. This means that $u_h(0) = 0$ and $-u_h(-h) = u_h(h) = \alpha$ for some $\alpha \in \mathbb{R}$. It
is easy to see that
\begin{align*}
  \int_{-1}^1\sgn(u-u_h)|u-u_h|^{q-1} \varphi_2 \dm x=0
\end{align*}
for any choice  of $\alpha$ and that
\begin{align*}
  \int_{-1}^1\sgn(u-u_h)|u-u_h|^{q-1} \varphi_3 \dm x = 0 && \Leftrightarrow && \int_{-1}^1\sgn(u-u_h)|u-u_h|^{q-1} \varphi_1 \dm x = 0.
\end{align*}
To determine for which $\alpha$ the latter two  integrals become zero, note that this is the same situation as in the example presented in Sections
\ref{sec:1d_l1_bd} and \ref{sec:1d_lq_bd}, only mirrored. Therefore, we again obtain
that $\alpha$ satisfies
\begin{align*}
  0 = -(1-h)\alpha^2 q(\alpha-1)^{q-1}-h(\alpha q+1)(\alpha-1)^q+h.
\end{align*}
This completes the proof of the first part of the lemma.

Since, $u_h$ is an odd function, for any $1<q<\infty$, the limit as $q \rightarrow 1$ must be an odd function as well and must therefore be zero at $x = 0$. Since the $L^1$-best approximation is uniquely determined by the value it takes at zero according to the first two parts of Theorem \ref{thm:discontinuity}, this completes the proof. Therefore, in the limit we obtain
the solution in Fig.\ \ref{fig:jump_l1_0} if $h\leq 1/2$. The corresponding $L^1$-best
approximation for $h>1/2$ is shown in Fig.\ \ref{fig:jump_l1_1_1}.
\end{proof}

\section{Best Approximation of a Boundary Discontinuity in Two Dimensions}\label{sec:2d_layer}
\begin{figure}[t!]
  \centering
  \begin{subfigure}{0.24\textwidth}
    \begin{tikzpicture}[thick,scale=0.65, every node/.style={scale=0.9}]
    \draw [thick] (0,0) rectangle (4,4);
    \draw [thick](0,0) -- (4,4);
    \draw [thick](4,0) -- (0,4);
    \node at (2,3) {$\tau_2$};
    \node at (3,2) {$\tau_1$};
    \node at (2,1) {$\tau_0$};
    \node at (1,2) {$\tau_3$};
    \node at (-0.2,-0.3) {(0,0)};
    \node at (-0.2,4.3) {(0,1)};
    \node at (4.2,-0.3) {(1,0)};
    \node at (4.2,4.3) {(1,1)};
    \end{tikzpicture}
    \caption{Mesh 1}
  \end{subfigure}
\begin{subfigure}{0.24\textwidth}
\begin{tikzpicture}[thick,scale=0.65, every node/.style={scale=0.9}]
\draw [thick] (0,0) rectangle (4,4);
\draw [thick] (0,2) -- (4,2);
\draw [thick] (2,0) -- (2,4);
\draw [thick] (0,0) -- (4,4);
\draw [thick] (2,0) -- (4,2);
\draw [thick] (0,2) -- (2,4);
\node at (-0.2,-0.3) {1};
\node at (-0.2,4.3) {1};
\node at (4.2,-0.3) {0};
\node at (4.2,4.3) {0};
\node at (-0.2, 2) {1};
\node at (2, 4.3) {$\bm \gamma$};
\node at (2, -0.3) {$\bm \alpha$};
\node at (2.3, 1.7) {$\bm \beta$};
\node at (1.3,0.7){$\tau_0$};
\node at (0.7,1.3){$\tau_1$};
\node at (1.3, 2.7){$\tau_4$};
\node at (0.7, 3.3){$\tau_5$};
\node at (3.3,0.7){$\tau_2$};
\node at (2.7,1.3){$\tau_3$};
\node at (3.3, 2.7){$\tau_6$};
\node at (2.7, 3.3){$\tau_7$};
\end{tikzpicture}
\caption{Mesh 2}\label{fig:mesh2}
\end{subfigure}
\begin{subfigure}{0.24\textwidth}
\begin{tikzpicture}[thick,scale=0.65, every node/.style={scale=0.9}]
\draw [thick] (0,0) rectangle (4,4);
\draw [thick] (0,2) -- (4,2);
\draw [thick] (2,0) -- (2,4);
\draw [thick] (0,0) -- (4,4);
\draw [thick] (0,4) -- (4,0);
\node at (-0.2,-0.3) {1};
\node at (-0.2,4.3) {1};
\node at (4.2,-0.3) {0};
\node at (4.2,4.3) {0};
\node at (2, 4.3) {$\bm \gamma$};
\node at (2, -0.3) {$\bm \alpha$};
\node at (2.2, 1.4) {$\bm \beta$};
\node at (1.3,0.7){$\tau_0$};
\node at (0.7,1.3){$\tau_1$};
\node at (1.3, 3.3){$\tau_4$};
\node at (0.7, 2.7){$\tau_5$};
\node at (3.3,1.3){$\tau_2$};
\node at (2.7,0.7){$\tau_3$};
\node at (3.3, 2.7){$\tau_6$};
\node at (2.7, 3.3){$\tau_7$};
\end{tikzpicture}
\caption{Mesh 3}
\end{subfigure}
\begin{subfigure}{0.24\textwidth}
\begin{tikzpicture}[thick,scale=0.65, every node/.style={scale=0.9}]
\draw [thick] (0,0) rectangle (4,4);
\draw [thick] (0,2) -- (4,2);
\draw [thick] (2,0) -- (2,4);
\draw [thick] (2,0) -- (4,2);
\draw [thick] (0,2) -- (2,4);
\draw [thick] (2,0) -- (0,2);
\draw [thick] (2,4) -- (4,2);
\node at (-0.2,-0.3) {1};
\node at (-0.2,4.3) {1};
\node at (4.2,-0.3) {0};
\node at (4.2,4.3) {0};
\node at (2, 4.3) {$\bm \gamma$};
\node at (2, -0.3) {$\bm \alpha$};
\node at (2.3, 1.7) {$\bm \beta$};
\node at (1.3,1.3){$\tau_0$};
\node at (0.7,0.7){$\tau_1$};
\node at (1.3, 2.7){$\tau_4$};
\node at (0.7, 3.3){$\tau_5$};
\node at (3.3,0.7){$\tau_2$};
\node at (2.7,1.3){$\tau_3$};
\node at (3.3, 3.3){$\tau_6$};
\node at (2.7, 2.7){$\tau_7$};
\end{tikzpicture}
\caption{Mesh 4}
\end{subfigure}
\caption{Four different meshes on $(0,1)^2$.}
\label{fig:meshes_2d}
\end{figure}
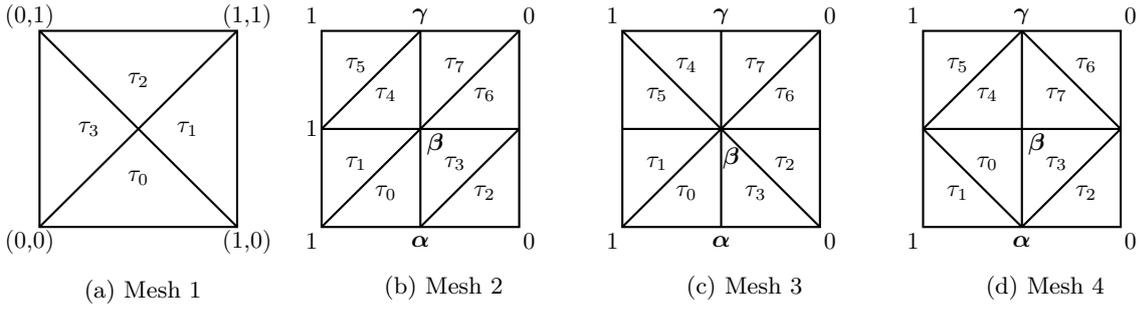

In this section we consider the best approximation problem \eqref{eq:best_approx} with $d=2$ and provide a proof of Theorem \ref{thm:best_approx_2d}. Hence, we consider the function $u\equiv 1$ on $(0,1)^2$. We consider the four meshes shown in Fig.\ \ref{fig:meshes_2d}
and determine the best approximation of $u$ by
a continuous function $u_h$ that is a linear polynomial on each of the triangles and takes the
following values in the four corners: $u_h(0,0)=u_h(0,1) = 1$ and $u_h(1,0)=u_h(1,1)=0$. For all meshes except the first one, we additionally fix the boundary conditions $u_h(0,0.5) = 1$ and $u_h(1,0.5) = 0$.

The free parameter of the best approximation problem for the first mesh is $\alpha = u_h(0.5,0.5)$;
there are three free parameters for each of the remaining meshes. For Meshes 2-4, we denote by $v_1$ the continuous piecewise
linear function that is $1$ at the node $(0.5,0)$ and $0$ at all other nodes;
by $v_2$ the continuous piecewise linear function that is $1$ at $(0.5,0.5)$ and $0$ at all
other nodes; and by $v_3$ the continuous piecewise linear function that is $1$ at
$(0.5,1)$ and zero at all other nodes. The coefficients defining the solution
$u_h$ are denoted as follows
\begin{align*}
  u(0.5,0) = \alpha,&& u(0.5,0.5) = \beta, && u(0.5,1) = \gamma.
\end{align*}

We split the proof of Theorem \ref{thm:best_approx_2d} into three parts: in Section \ref{sec:l1_2d}, we prove the first part of the theorem, i.e., we consider Mesh 1 with $q=1$; in Section \ref{sec:lq_2d} we continue with part two of the theorem and consider Mesh 1 with $1<q<\infty$; Section \ref{sec:2d_l1_rest} finally contains a proof of parts three, four and five of the theorem, i.e., we show that the $L^q$-best approximation contains over- or undershoots on all four meshes if $q>1$ and consider $q=1$ for Meshes 2, 3 and 4.

\subsection{$L^1$-Best Approximation}\label{sec:l1_2d}
In this section we prove the first part of Theorem \ref{thm:best_approx_2d}. More precisely, we show that if we consider Mesh 1, the $L^1$-best approximation is unique and  $u_h(0.5,0.5) = \alpha$, where $\alpha$ satisfies
\begin{align*}
	\alpha >1 &&\text{ and } &&0 = 2 \alpha^3- 5\alpha +2,
\end{align*}
hence $\alpha \approx 1.3200$.

\begin{proof}
 We again use the characterisation of the $L^1$-best approximation in Corollary \ref{cor:best_approx_l1}. The space $U_h$ is the span of the continuous function $v_h$ that is a linear polynomial on each element,
zero at the boundary of the domain and $1$ at the centroid $(0.5,0.5)$.
We will use the reference triangle $\hat{\tau}$ as depicted on the left in Fig.\ref{fig:crisscross_pos_neg}
 for all computations. To this end, we define the affine transformations
$\xi_i : \tau_i \rightarrow \hat{\tau}$, $i=0,1,2,3$ that are each composed of a scaling by $0.5$,
a rotation and a translation. On $\hat{\tau}$ we define the basis
functions $\hat{\varphi}_0$, $\hat{\varphi_1}$ and $\hat{\varphi}_2$ as
\begin{align*}
  \hat{\varphi}_0 = 1-x-y, \qquad \hat{\varphi}_1 = x, \qquad \hat{\varphi}_2 = y,
\end{align*}
{  in the coordinates of the reference element $\hat{\tau}$}.
Note that $\xi_i(v_h|_{\tau_i}) = \hat{\varphi}_0$ for all $i=0,1,2,3$.

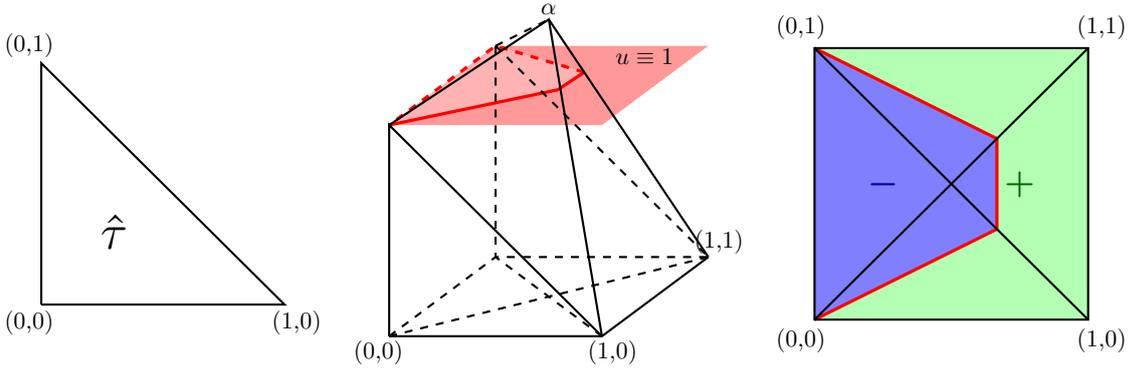
\begin{figure}
  \centering
  \begin{subfigure}{0.28\textwidth}
  \begin{tikzpicture}[thick,scale=0.8, every node/.style={scale=0.9}]
    \draw[thick] (0,0)--(4,0)--(0,4)--(0,0);
    \node at (-0.2,-0.3) {(0,0)};
    \node at (-0.2,4.3) {(0,1)};
    \node at (4.2,-0.3) {(1,0)};
    \node at  (1.2, 1.2){\huge$\hat{\tau}$};
  \end{tikzpicture}
\end{subfigure}
\begin{subfigure}{0.34\textwidth}
\begin{tikzpicture}[thick,scale=0.7, every node/.style={scale=0.9}]
\fill[red!40!white] (0,4)-- (4,4)--(6,5.5)--(2, 5.5) -- (0,4);
\fill[red!30!white] (0,4)-- (3.1666, 4.6666)--(3.6666,5)--(2,5.5) -- (0,4);
\draw [red, very thick] (0,4) -- (3.1666, 4.6666)--(3.6666,5);
\draw [red, very thick, dashed] (3.6666,5)--(2,5.5) -- (0,4);
\draw [thick](0,0) -- (4,0) -- (6,1.5);
\draw [dashed] (0,0) -- (2,1.5) -- (6,1.5);
\draw [dashed] (0,0) -- (6,1.5);
\draw [dashed] (4,0) -- (2,1.5);
\draw [thick] (0,0) -- (0,4) -- (4,0);
\draw [thick] (0,4) -- (3, 6) -- (4,0);
\draw [thick] (3,6) -- (6, 1.5);
\draw [dashed] (2, 5.5) -- (3,6);
\draw [dashed] (2, 5.5) -- (2, 1.5);
\draw [dashed] (2, 5.5) -- (6, 1.5);

\node at (-0.2,-0.3) {(0,0)};
\node at (6.2,1.8) {(1,1)};
\node at (4.2,-0.3) {(1,0)};
\node at (3,6.2){$\alpha$};
\node at (4.8, 5.3){$u\equiv 1$};
\end{tikzpicture}
\end{subfigure}
\begin{subfigure}{0.34\textwidth}
\begin{tikzpicture}[thick,scale=0.9, every node/.style={scale=0.9}]
\fill[green!30!white] (0,0) -- (4, 0)-- (4,4) -- (0,4) -- (0,0);
\fill[blue!50!white] (0,4)--(2.6666, 2.6666) -- (2.6666,1.3333) -- (0,0) -- (0,4);
\draw[very thick, red ] (0,4)--(2.6666, 2.6666) -- (2.6666,1.3333) -- (0,0);
\draw [thick] (0,0) rectangle (4,4);
\draw [thick](0,0) -- (4,4);
\draw [thick](4,0) -- (0,4);
\node[blue!60!black] at (1,2) {\Large$\bm-$};
\node[green!40!black] at (3,2) {\Large$\bm+$};
\node at (-0.2,-0.3) {(0,0)};
\node at (-0.2,4.3) {(0,1)};
\node at (4.2,-0.3) {(1,0)};
\node at (4.2,4.3) {(1,1)};
\end{tikzpicture}
\end{subfigure}
\caption{Left: Reference element $\hat{\tau}$. Center: The function $u_h$ with $\alpha >1$. The intersection with $u$ is marked with red lines.
Right: The mesh with the area where $(u-u_h)<0$ coloured in blue and the area where $(u-u_h)>0$ coloured in green.}
\label{fig:crisscross_pos_neg}
\end{figure}

We consider the following two cases:
\begin{enumerate}
  \item The set $\{(x,y) \in (0,1)^2 \,: \,(u-u_h)(x,y) = 0\}$ has measure zero.
  \item The set $\{(x,y) \in (0,1)^2 \,: \,(u-u_h)(x,y) = 0\}$ has positive measure.
\end{enumerate}
The second case is only possible if $\alpha =1$. In this case we have
$u_h = u$ in $\tau_3$ and $u_h < u$ in $\tau_i$, $i=0,1,2$.
Noting that the $\tau_i$ and $\hat{\tau}$ are similar and that the area of
$\hat{\tau}$ is precisely twice the area of any $\tau_i$, we obtain that
\begin{align*}
  \sum_{i=0}^2 \int_{\tau_i} \psi v_h \mathrm{d}\bm x = \sum_{i=0}^2 \int_{\tau_i} \mathrm{sgn}(u-u_h)v_h \mathrm{d}\bm x
  = \frac{3}{2}\int_{\hat{\tau}}\hat{\varphi}_0\mathrm{d}\bm x
  = \frac{3}{2}\int_0^1\int_0^{1-y}(1-x-y)\mathrm{d}x\mathrm{d}y = \frac{1}{4}.
\end{align*}
On the other hand,
\begin{align*}
  \int_{\tau_3}\psi v_h \mathrm{d}\bm x \geq -\int_{\tau_3} v_h \mathrm{d}\bm x  =- \frac{1}{2}\int_{\hat{\tau}} \hat{\varphi}_0\mathrm{d}\bm x = -\frac{1}{12},
\end{align*}
where $\psi(x,y)$  arbitrary on $\tau_3$ with $-1\leq \psi(x,y) \leq 1$. Hence, for any
choice of $\psi$, we have
\begin{align*}
  \sum_{i=0}^3 \int_{\tau_i} \psi v_h \mathrm{d}\bm x \geq \frac{1}{4}-\frac{1}{12} > 0.
\end{align*}
This shows that $\alpha = 1$ does not yield a best approximation on this mesh.

We can furthermore rule out the case $\alpha < 1$; indeed, in this case
$u-u_h > 0$ in the whole domain $(0,1)^2$ and since $v_h > 0$ in $(0,1)^2$, we have
\begin{align*}
\int_{(0,1)^2} \psi v_h \mathrm{d}\bm x=  \int_{(0,1)^2} \mathrm{sgn}(u-u_h)v_h \mathrm{d}\bm x >0.
\end{align*}
The only remaining case is $\alpha >1$; in this setting we have
$u-u_h < 0$ in $\tau_3$ and thus
\begin{align*}
\int_{\tau_3} \psi v_h \dm \bm{x}  =\int_{\tau_3} \mathrm{sgn}(u-u_h)v_h \dm \bm{x} = -\frac{1}{2}\int_{\hat{\tau}} \hat{\varphi}_0 \dm \bm x = -\frac{1}{12}.
\end{align*}
In the remaining three triangles, $u-u_h$ changes sign within the element.
This is illustrated in Fig. \ref{fig:crisscross_pos_neg}. Furthermore, $u-u_h=0$ on a set of measure zero, hence $\psi = \sgn(u-u_h)$ almost everywhere. In order
to determine the sections of each element where $u-u_h$ is positive and negative, respectively,
we will consider $\xi_i(u-u_h)$ for $i=0,1,2$.

 Due to the symmetry of the approximation problem, we
can assume
\begin{align*}
  \int_{\tau_0}\mathrm{sgn}(u-u_h)v_h \dm \bm x = \int_{\tau_2}\mathrm{sgn}(u-u_h)v_h \dm \bm x.
\end{align*}
We compute
\begin{align}
  \xi_0(u-u_h) = 1 - \alpha \hat{\varphi}_0 -\hat{\varphi}_1 = -(\alpha-1)+(\alpha-1)x+\alpha y.\label{eq:xi_0}
\end{align}
Thus,
\begin{align}
  \begin{aligned}
    \xi_0(u-u_h)>0 \text{ on } \left\{(x,y)\in \hat{\tau} \,:\, y > \frac{\alpha-1}{\alpha}(1-x) \right\},\\
    \xi_0(u-u_h)<0 \text{ on } \left\{(x,y)\in \hat{\tau} \,:\, y < \frac{\alpha-1}{\alpha}(1-x) \right\}
  \end{aligned} \label{eq:xi_0_sgn}
\end{align}
and
\begin{align*}
  \begin{aligned}
  \int_{\tau_0}\mathrm{sgn}(u-u_h)v_h \dm \bm x &= \frac{1}{2}\int_{\hat{\tau}}\hat{\varphi}_0 \dm \bm x
  -\int_0^{1} \int_0^{\frac{\alpha-1}{\alpha}(1-x)} \hat{\varphi}_0 \dm y \dm x\\
  & = \frac{1}{12} - \int_0^1 \int_0^{\frac{\alpha-1}{\alpha}(1-x)} (1-x-y)\dm y \dm x\\
  & = \frac{1}{12} - \frac{1}{6}\left(\frac{\alpha-1}{\alpha}\right)\left(\frac{\alpha+1}{\alpha}\right).
\end{aligned}
\end{align*}
Similarly,
\begin{align}
  \xi_1(u-u_h) = 1 - \alpha \hat{\varphi}_0 = -(\alpha-1)+\alpha x +\alpha y.\label{eq:xi_1}
\end{align}
Thus,
\begin{align}
  \begin{aligned}
    \xi_1(u-u_h)>0 \text{ on } \left\{(x,y)\in \hat{\tau} \,:\, y > \frac{\alpha-1}{\alpha}-x \right\},\\
    \xi_1(u-u_h)<0 \text{ on } \left\{(x,y)\in \hat{\tau}\,:\, y < \frac{\alpha-1}{\alpha}-x \right\}
  \end{aligned}\label{eq:xi_1_sgn}
\end{align}
and
\begin{align*}
  \begin{aligned}
  \int_{\tau_1}\mathrm{sgn}(u-u_h)v_h \dm \bm x &= \frac{1}{2}\int_{\hat{\tau}}\hat{\varphi}_0 \dm \bm x
  -\int_0^{\frac{\alpha-1}{\alpha}} \int_0^{\frac{\alpha-1}{\alpha}-x} \hat{\varphi}_0 \dm y \dm x\\
  & = \frac{1}{12} - \int_0^{\frac{\alpha-1}{\alpha}} \int_0^{\frac{\alpha-1}{\alpha}-x} (1-x-y)\dm y \dm x\\
  & = \frac{1}{12} - \frac{1}{6}\left(\frac{\alpha-1}{\alpha}\right)^2\left(\frac{\alpha+2}{\alpha}\right).
\end{aligned}
\end{align*}
Finally,
\begin{align*}
  \begin{aligned}
    \int_{(0,1)^2} \mathrm{sgn}(u-u_h)v_h \dm \bm x &= \sum_{i=0}^3\int_{\tau_i}\mathrm{sgn}(u-u_h)v_h \dm \bm x\\
    &=\frac{1}{6}\left[1-2\left(\frac{\alpha-1}{\alpha}\right)\left(\frac{\alpha+1}{\alpha}\right)-\left(\frac{\alpha-1}{\alpha}\right)^2\left(\frac{\alpha+2}{\alpha}\right)\right]\\
    &=-\frac{1}{6 \alpha^3}\left(2\alpha^3- 5\alpha +2\right).
  \end{aligned}
\end{align*}
The polynomial $2\alpha^3-5\alpha+2$ has three roots $\alpha_i$, $i=0,1,2$, where
\begin{align*}
  \alpha_0 \approx -1.7623, && \alpha_1 \approx 0.43232, && \alpha_2 \approx 1.3200.
\end{align*}
Only $\alpha_2$ satisfies the condition $\alpha > 1$ and therefore $\alpha \approx 1.3200$
yields the only $L^1$-best approximation.
\end{proof}

\begin{figure}
 \centering

  \begin{subfigure}{0.65\textwidth}
      \centering
\begin{tikzpicture}[thick,scale=0.9, every node/.style={scale=0.9}]
\draw [thick] (0,0) rectangle (4,4);
\draw [thick](0,0) -- (4,4);
\draw [thick](4,0) -- (0,4);
\draw [thick,->] (4.2,2)--(4.8,2);
\end{tikzpicture}
\begin{tikzpicture}[thick,scale=0.9, every node/.style={scale=0.9}]
  \draw [thick] (0,0) rectangle (4,4);
  \draw [thick](0,0) -- (4,4);
  \draw [thick](4,0) -- (0,4);
  \draw [thick](0,2)--(4,2);
  \draw [thick] (2,0)--(2,4);
  \draw [thick] (0,2)--(2,0);
  \draw [thick] (4,2)--(2,0);
  \draw [thick] (4,2)--(2,4);
  \draw [thick] (0,2)--(2,4);
\end{tikzpicture}
\caption{Refinement for which the overshoot in the $L^1$-best approximation remains constant.}
 \label{fig:crisscross_refinement}
\end{subfigure}
\begin{subfigure}{0.32\textwidth}
    \centering
    \begin{tikzpicture}[thick,scale=0.9, every node/.style={scale=0.9}]
    \fill[green!30!white] (0,0) -- (4, 0)-- (4,4) -- (0,4) -- (0,0);
    \fill[blue!50!white] (2,2)--(3.3333, 1.3333) -- (3.3333,0.6666) -- (2,0)--(2,2);
    \fill[blue!50!white](2,4)--(3.3333, 3.3333) -- (3.3333,2.6666) -- (2,2) --(2,4);
    \fill[blue!50!white](2,2)--(0.6666, 1.3333) -- (0.6666,0.6666) -- (2,0)--(0,0)--(0,2)--(2,2);
    \fill[blue!50!white](2,4)--(0.6666, 3.3333) -- (0.6666,2.6666) -- (2,2)--(0,2)--(0,4)--(2,4);
    \draw[very thick, red ] (2,2)--(3.3333, 1.3333) -- (3.3333,0.6666) -- (2,0);
    \draw[very thick, red ] (2,4)--(3.3333, 3.3333) -- (3.3333,2.6666) -- (2,2);
    \draw[very thick, red ] (2,2)--(0.6666, 1.3333) -- (0.6666,0.6666) -- (2,0);
    \draw[very thick, red ] (2,4)--(0.6666, 3.3333) -- (0.6666,2.6666) -- (2,2);
    \draw [thick] (0,0) rectangle (4,4);
    \draw [thick](0,0) -- (4,4);
    \draw [thick](4,0) -- (0,4);
    \draw [thick](0,2)--(4,2);
    \draw [thick] (2,0)--(2,4);
    \draw [thick] (0,2)--(2,0);
    \draw [thick] (4,2)--(2,0);
    \draw [thick] (4,2)--(2,4);
    \draw [thick] (0,2)--(2,4);

    \node[blue!60!black] at (2.5,3) {\Large$\bm{-1}$};
    \node[green!40!black] at (1.5,3) {\Large$\bm{+1}$};
    \end{tikzpicture}
    \caption{Possible choice for $\sgn (0)$}\label{fig:crisscross_sign}
\end{subfigure}
\caption{Uniform refinement of the mesh preserving the structure.}
\end{figure}
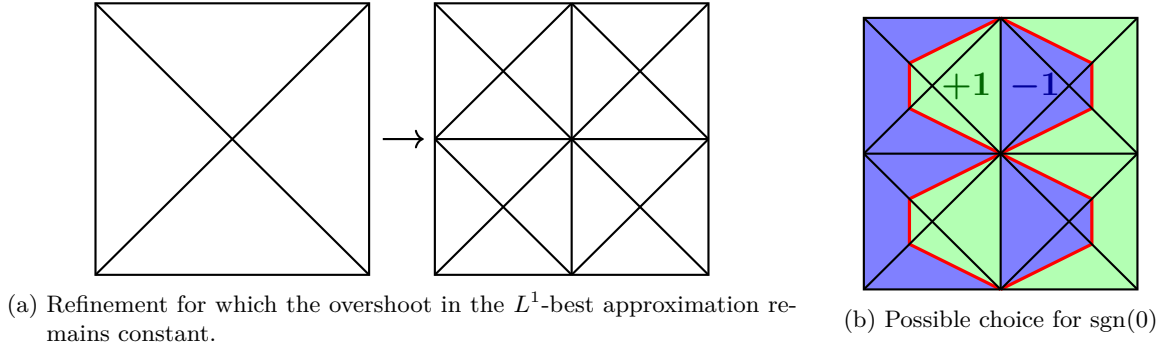
\begin{rem}[Uniform refinement]\label{rem:refinement_crisscross}
  If the mesh is refined uniformly, keeping the same structure as
  shown in Fig.\ \ref{fig:crisscross_refinement},   it is easy to see
  that an $L^1$-best approximation is given by $u(x_i,y_j) = \alpha $, with $\alpha$ as specified in Section \ref{sec:l1_2d}, if the node $(x_i, y_j)$ is connected with the boundary $x =1$, and $u(x_i, y_j)=1$ at the remaining interior nodes. Indeed, in
  this case we can choose $\psi = \psi_0$ on the set $\{x\,:\, u(x)-u_h(x) = 0\}$ as shown
  in Fig.\ \ref{fig:crisscross_sign}. This shows that the overshoot
  in the $L^1$-best approximation remains constant under this type
  of mesh refinement.
\end{rem}
\bigskip
\subsection{$L^q$-Best Approximation}\label{sec:lq_2d}
 In this section we prove the second part of Theorem \ref{thm:best_approx_2d}. More precisely, we show that the unique solution to \eqref{eq:best_approx} on Mesh 1 for $1<q<\infty$ is defined by $u_h(0.5,0.5) = \alpha$, where $\alpha$ satisfies
  \begin{align*}
    \alpha>1 &&\text{ and } && 0 = (\alpha-1)^{q-1}\left[4\alpha^3 q +4(1-q)\alpha^2+(q-6)\alpha+2\right]-\alpha(q+4)+2.
  \end{align*}
\begin{proof}
Corollary \ref{cor:best_lq} implies that we seek $\alpha$ such that
\begin{align}
  \int_{(0,1)^2} \sgn(u-u_h) |u-u_h|^{q-1} v_h \dm \bm x = 0.\label{eq:Lqbest}
\end{align}
Again, we have to split the integral on each element
into the parts where $u-u_h>0$, $u-u_h <0$ and $u-u_h = 0$.
We consider the three cases
\begin{align*}
  \text{ (a) } \alpha < 1,&&
  \text{ (b) } \alpha = 1,&&
  \text{ (c) } \alpha >1.
\end{align*}
If $\alpha < 1$, we have $u-u_h>0$ everywhere in $(0,1)^2$ and thus
both $\sgn(u-u_h)|u-u_h|^{q-1}>0$ and $v_h>0$ in $(0,1)^2$. Therefore, we have
\begin{align*}
  \int_{(0,1)^2} \sgn(u-u_h) |u-u_h|^{q-1} v_h \dm \bm x >0.
\end{align*}
If $\alpha = 1$, we have $u-u_h = 0$ in $\tau_3$ and $u-u_h > 0$ in $\tau_i$, $i=0,1,2$. Thus,
\begin{align*}
  \int_{(0,1)^2} \sgn(u-u_h) |u-u_h|^{q-1} v_h \dm \bm x
  = \sum_{i=1}^2 \int_{\tau_i } \sgn(u-u_h) |u-u_h|^{q-1} v_h \dm \bm x > 0.
\end{align*}
We can therefore again assume $\alpha >1$; in this case $u-u_h < 0$ in $\tau_3$ and
we compute
\begin{align*}
  \begin{aligned}
  \int_{\tau_3} \sgn(u-u_h) |u-u_h|^{q-1} v_h \dm \bm x
  &= -\frac{1}{2}\int_{\hat{\tau}}(\xi_3(u_h- u))^{q-1}\hat{\varphi}_0 \dm \bm x\\
  &= -\frac{1}{2}\int_{\hat{\tau}} (\alpha - 1)^{q-1}\hat{\varphi}_0^q \dm \bm x
  = -\frac{(\alpha-1)^{q-1}}{2(q+1)(q+2)}.
\end{aligned}
\end{align*}
Next, consider $\tau_0$ and $\tau_2$; using the symmetry of the problem and equations
\eqref{eq:xi_0} and \eqref{eq:xi_0_sgn}, we obtain
\begin{align*}
  \begin{aligned}
  \int_{\tau_2}\mathrm{sgn}&(u-u_h)v_h \dm \bm x+\int_{\tau_0}\mathrm{sgn}(u-u_h)v_h \dm \bm x\\
  &=2\int_{\tau_0}\mathrm{sgn}(u-u_h)v_h \dm \bm x\\
  &= \int_{\hat{\tau}}\sgn(\xi_0(u-u_h))|\xi_0(u-u_h)|^{q-1}\hat{\varphi}_0 \dm \bm x\\
  &=\int_0^1\int_{\frac{\alpha-1}{\alpha}(1-x)}^{1-x}(-(\alpha-1)+(\alpha-1)x+\alpha y)^{q-1}(1-x-y)\dm y \dm x\\
  &\quad-\int_0^1\int_{0}^{\frac{\alpha-1}{\alpha}(1-x)}((\alpha-1)-(\alpha-1)x-\alpha y)^{q-1}(1-x-y)\dm y \dm x\\
  & = \frac{1}{\alpha^2 q(q+1)(q+2)} - \frac{(\alpha-1)^q(\alpha q+1)}{\alpha^2q(q+1)(q+2)}.
\end{aligned}
\end{align*}
Finally, we consider $\tau_1$; using  equations
\eqref{eq:xi_1} and \eqref{eq:xi_1_sgn}, we obtain
\begin{align*}
  \begin{aligned}
  \int_{\tau_1}\mathrm{sgn}&(u-u_h)v_h \dm \bm x\\
  &= \frac{1}{2} \int_{\hat{\tau}}\sgn(\xi_1(u-u_h))|\xi_1(u-u_h)|^{q-1}\hat{\varphi}_0 \dm \bm x\\
  &= \frac{1}{2}\int_0^{\frac{\alpha-1}{\alpha}}\int_{\frac{\alpha-1}{\alpha}-x}^{1-x}(-(\alpha-1)+\alpha x+\alpha y)^{q-1}(1-x-y)\dm y \dm x\\
  & \quad +\frac{1}{2}\int_{\frac{\alpha-1}{\alpha}}^1\int_{0}^{1-x}(-(\alpha-1)+\alpha x+\alpha y)^{q-1}(1-x-y)\dm y \dm x\\
  &\quad-\frac{1}{2}\int_0^{\frac{\alpha-1}{\alpha}}\int_{0}^{\frac{\alpha-1}{\alpha}-x}((\alpha-1)- \alpha x-\alpha y)^{q-1}(1-x-y)\dm y \dm x\\
  & = \frac{\alpha-1}{2\alpha^3 q(q+1)} + \frac{1}{2\alpha^3(q+1)(q+2)}
  -\frac{(\alpha-1)^q(\alpha^2q+(2-q)\alpha - 2)}{2\alpha^3q(q+1)(q+2)}.
\end{aligned}
\end{align*}
Combining all integrals, we obtain
\begin{align*}
  \begin{aligned}
    \int_{(0,1)^2} \mathrm{sgn}&(u-u_h)v_h \dm \bm x = \sum_{i=0}^3\int_{\tau_i}\mathrm{sgn}(u-u_h)v_h \dm \bm x\\
    &= -\frac{(\alpha-1)^{q-1}\left[4\alpha^3 q +4(1-q)\alpha^2+(q-6)\alpha+2\right]-\alpha(q+4)+2}{2\alpha^3 q(q+1)(q+2)}.
  \end{aligned}
\end{align*}
Hence, the $L^q$ best approximation can be determined by finding $\alpha  >1$ satisfying
\begin{align*}
  0 = (\alpha-1)^{q-1}\left[4\alpha^3 q +4(1-q)\alpha^2+(q-6)\alpha+2\right]-\alpha(q+4)+2.
\end{align*}
\end{proof}

\subsection{$L^q$-Best Approximation on Meshes 2, 3 and 4} \label{sec:2d_l1_rest}
In this section we prove parts three, four and five of Theorem \ref{thm:best_approx_2d}. More precisely, we show the following:
  \begin{itemize}
		\item If $q>1$, the $L^q$-best approximation to \eqref{eq:best_approx} contains over- or undershoots on all four meshes.
  \item If $q=1$, there exists a solution to \eqref{eq:best_approx} on Mesh 2 such that $u_h(0.5, 1) = u_h(0.5,0.5) = 1$ and $u_h(0.5,0) = \alpha$, where $\alpha$ satisfies
    $\alpha > 1$  and $0= -3\alpha^3 + 8 \alpha -4$.

  Furthermore, $u_h(0.5, 1) = u_h(0.5,0.5) = u_h(0.5,0) = 1$ does not define an $L^1$-best approximation.
  \item If $q=1$, there exists a solution to \eqref{eq:best_approx} on Mesh 3 and  Mesh 4 such that $u_h(0.5, 1) = u_h(0.5,0.5) = u_h(0.5,0) = 1$.
\end{itemize}

\begin{proof}
To see that the first point (part three of Theorem \ref{thm:best_approx_2d}) is true, note that if there are no over- or undershoots, i.e., $u_h(0.5, 1) = u_h(0.5,0.5) = u_h(0.5,0) = 1$, we have $u-u_h = 0$ in $(0,0.5)\times(0,1)$ and $u-u_h>0$ in $(0.5,1)\times (0,1)$ and hence
\begin{align*}
	\int_{(0,1)^2}\sgn(u-u_h)|u-u_h|^{q-1} v_i \dm \bm x >0 \quad
	\text{ for all } i = 1,2,3 \text{ and } 1<q<\infty,
\end{align*}
which contradicts Corollary \ref{cor:best_lq}.

Furthermore, it is easy to see that for the second and third mesh,
$\alpha = \beta = \gamma = 1$ is an $L^1$-best approximation,  again using the characterisation in Corollary \ref{cor:best_approx_l1}.
 Indeed, we have $\psi = \sgn (u-u_h)$ whenever $u(\bm x)\neq u_h(\bm x)$ and we can choose $\psi = \psi_0 \equiv 1$ otherwise. In this case, for each of the three nodes $(0.5,0)$, $(0.5, 0.5)$ and $(0.5,1)$, we have that $\psi=-1$ on exactly half of the
connected elements and $\psi=1$ on the remaining connected elements.
Note that
\begin{align*}
  \int_{\tau_i}v_j \dm \bm x = \frac{1}{24} \qquad \forall i,j,
\end{align*}
and thus the terms with $\psi =-1$ and $\psi=1$
cancel each other. This proves the third point (part five of Theorem \ref{thm:best_approx_2d}).

The second point (part four of Theorem \ref{thm:best_approx_2d}) is more interesting. The above argument could now only
be applied to $v_2$ which is $1$ at the node $(0.5,0.5)$. First consider the node
$(0.5,1)$; the connected elements are $\tau_5$, $\tau_4$ and $\tau_7$, cf., Fig.\ \ref{fig:mesh2}.
If $\beta = \gamma = 1$, we have that $u-u_h = 0$ in $\tau_4$ and $\tau_5$ and
$u-u_h>0$ in $\tau_7$. If we choose $\psi \equiv -1$ on $\tau_4$ and $\psi\equiv 0$ on $\tau_5$, we obtain
\begin{align}
  \begin{aligned}
\int_{(0,1)^2}\psi v_3\dm \bm x
&= \int_{\tau_7}v_3 \dm \bm x-\int_{\tau_4}v_3\dm \bm x \
&= \frac{1}{24}-\frac{1}{24} = 0.\label{eq:v3}
\end{aligned}
\end{align}
Note that this is independent of $\alpha$. Considering $v_1$ shows that
there is no $L^1$-best approximation with $\alpha = \beta =1$. Indeed in this case, we have
for any $-1 \leq \psi \leq 1$,
\begin{align*}
\begin{aligned}
\int_{(0,1)^2}\psi v_1\dm \bm x
&\geq \int_{\tau_0}v_1 \dm \bm x+\int_{\tau_1}v_1\dm \bm x -\int_{\tau_3}v_1\dm \bm x\\
&= \frac{1}{24}+\frac{1}{24} -\frac{1}{24} = \frac{1}{24}>0.
\end{aligned}
\end{align*}
We will now show that there is an $L^1$-best approximation with
$\beta = \gamma =1$ and $\alpha>1$. We have already established that \eqref{eq:v3}
is independent of $\alpha$; this leaves the integrals involving
$v_2$ and $v_3$.  We have $u-u_h >0$ in $\tau_i$, $i =3,6,7$, and $u-u_h < 0$ in
$\tau_0$. Moreover, we have already fixed $\psi \equiv -1$ in $\tau_4$. If we now
furthermore choose $\psi \equiv -1$ in $\tau_1$, we obtain
\begin{align*}
  \int_{(0,1)^2}\psi v_2 \dm \bm x = \sum_{i =3,6,7}\int_{\tau_i}v_2 \dm \bm x -\sum_{i =0,1,4}\int_{\tau_i}v_2 \dm \bm x = 0;
\end{align*}
this is again independent of $\alpha$. Finally, we consider $v_3$ to determine
$\alpha$. Let again $\xi_i$ be the linear
transformation that maps $\tau_i$ onto $\hat{\tau}$;
we have $u-u_h>0$ on $\tau_0$. Furthermore,
\begin{align*}
  \xi_2(u-u_h) = 1-\alpha y, && \xi_3(u-u_h) = (1-\alpha) x +y.
\end{align*}
Thus,
\begin{align*}
  \begin{aligned}
    \xi_2(u-u_h)>0 \text{ on } \left\{(x,y)\in \hat{\tau} \,:\, y < \frac{1}{\alpha}\right\},\\
    \xi_2(u-u_h)<0 \text{ on } \left\{(x,y)\in \hat{\tau}\,:\, y > \frac{1}{\alpha}\right\}
  \end{aligned}
\end{align*}
and
\begin{align*}
  \begin{aligned}
    \xi_3(u-u_h)>0 \text{ on } \left\{(x,y)\in \hat{\tau} \,:\, y > (\alpha-1)x\right\},\\
    \xi_3(u-u_h)<0 \text{ on } \left\{(x,y)\in \hat{\tau}\,:\, y < (\alpha-1)x\right\}.
  \end{aligned}
\end{align*}
Therefore, we compute
\begin{align*}
  \begin{aligned}
\int_{\tau_2} \psi v_1\dm \bm x = \int_{\tau_2} \sgn{(u-u_h)}v_1\dm \bm x
  &= \frac{1}{4}\int_{\hat{\tau}}y\dm x \dm y - \frac{1}{2}\int_{\frac{1}{\alpha}}^1\int_0^{1-y}y \dm x \dm y\\
  &= \frac{-\alpha^3+6\alpha -4}{24\alpha^3},
\end{aligned}
\end{align*}
\begin{align*}
  \begin{aligned}
  \int_{\tau_3} \psi v_1\dm \bm x=\int_{\tau_3} \sgn{(u-u_h)}v_1\dm \bm x
  &= -\frac{1}{4}\int_{\hat{\tau}}x\dm x\dm y +\frac{1}{2} \int_{0}^{\frac{1}{\alpha}}\int_{(1-\alpha)x}^{1-x}x \dm x \dm y\\
  &= \frac{2-\alpha^2}{24\alpha^2}
\end{aligned}
\end{align*}
and
\begin{align*}
  \int_{\tau_0} \psi v_1\dm \bm x=\int_{\tau_0} \sgn{(u-u_h)}v_1\dm \bm x = -\frac{1}{24}.
\end{align*}
Putting all three integrals together yields
\begin{align*}
  \int_{(0,1)^2}\psi v_1\dm \bm x = \frac{1}{24\alpha^3}(-3\alpha^3+8\alpha-4).
\end{align*}
Hence, $\alpha > 1$ has to satisfy the equation
\begin{align*}
  0 = -3\alpha^3+8\alpha-4.
\end{align*}
Three roots of the above polynomial are $\alpha_0  \approx -1.8414$, $\alpha_1 \approx 0.56913$ and
$\alpha_2 \approx 1.2723$. Only the third root satisfies $\alpha > 1$ and therefore defines an
$L^1$-best approximation.

\end{proof}
\section{Numerical Examples}\label{sec:examples}
In this section we consider selected examples of meshes for which we have determined the solution of the best approximation problem \eqref{eq:best_approx} numerically by interpreting the condition \eqref{eq:best_lq} as a variational problem that can be implemented using standard finite element techniques. Here, we have used FEniCS \cite{Alnes2015} for the implementation.
In Section \ref{sec:refinement} we illustrate that the overshoot in the $L^q$-best approximation does not vanish if the mesh is refined and that these observations even apply if $u$ is a more general smooth function.

In Section \ref{sec:examples_2d} we illustrate that the $L^q$-approximation on the three meshes considered in Section \ref{sec:2d_l1_rest} (second half of Theorem \ref{thm:best_approx_2d}) converges to the $L^1$-best approximation characterised in the theorem. Furthermore, we show how the understanding of these special cases can be applied to predict the behaviour of the $L^q$-best approximation on a more general mesh.

\subsection{Refinement of the Mesh}\label{sec:refinement}
\subsubsection{Gibbs Phenomenon on Meshes in Two Dimensions}

We start with Mesh 1 depicted in Fig.\ \ref{fig:meshes_2d} and the refinement shown in Fig.\ \ref{fig:crisscross_refinement} that preserves the structure of the mesh. We have already shown in Remark \ref{rem:refinement_crisscross} that for $q=1$ there exists an $L^1$-best approximation such that the overshoot remains constant as we refine the mesh.
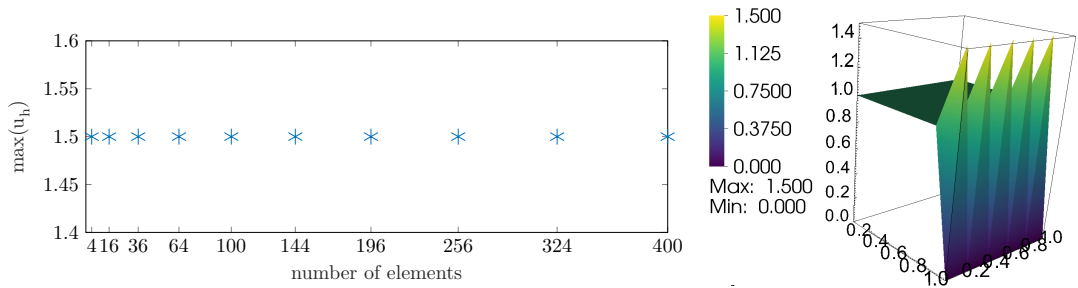
\begin{figure}[t]
  \centering
%
%
\definecolor{mycolor1}{rgb}{0.00000,0.44700,0.74100}%
\begin{tikzpicture}[thick, scale=0.5, every node/.style={scale=1.5}]

\begin{axis}[%
width=6.028in,
height=2in,
at={(1.011in,0.642in)},
scale only axis,
xmin=0,
xmax=400,
xlabel style={font=\color{white!15!black}},
xlabel={number of elements},
xtick = data,
ymin=1.4,
ymax=1.6,
ylabel style={font=\color{white!15!black}},
ylabel={$\text{max(u}_\text{h}\text{)}$},
axis background/.style={fill=white}
]
\addplot [color=mycolor1, draw=none, mark=asterisk, mark options={solid, mycolor1, scale = 3.0}, forget plot]
  table[row sep=crcr]{%
4	1.5\\
16	1.5\\
36	1.5\\
64	1.5\\
100	1.5\\
144	1.5\\
196	1.5\\
256	1.5\\
324	1.5\\
400	1.5\\
};
\end{axis}
\end{tikzpicture}%
	\includegraphics[width = 5cm]{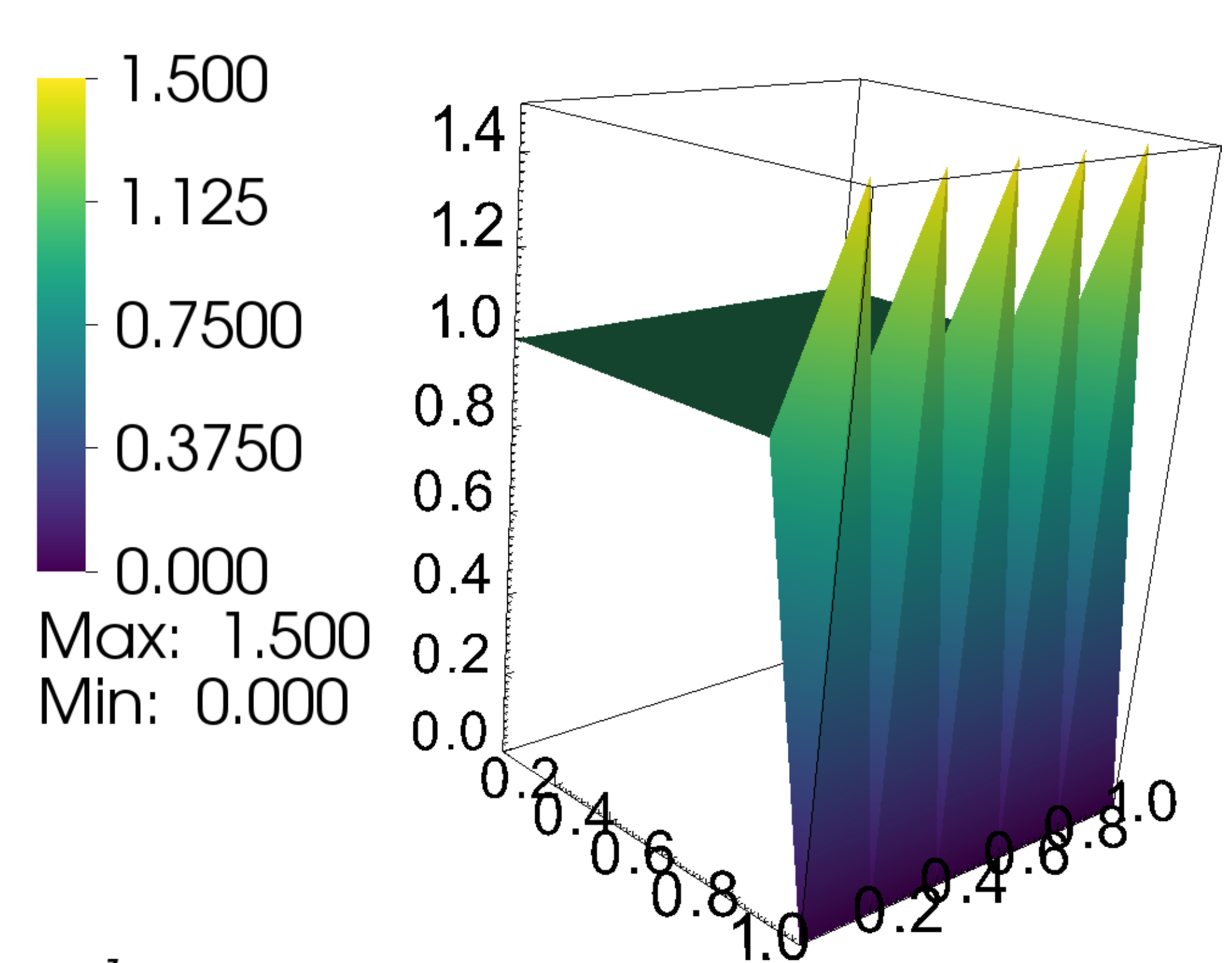}
  \caption{Left: $\max(u_h)$ for $q=2$ and several refinements as shown in
  Fig.\ \ref{fig:crisscross_refinement}.  Right: $u_h$ with 100 elements.}
  \label{fig:crisscross_q2}
\end{figure}
Indeed, Fig.\ \ref{fig:crisscross_q2} shows the maximum value of $u_h$ for this example with $q = 2$
 and for several refinements of the mesh, as well as the approximation $u_h$ for a mesh with this structure consisting of $100$ elements. We can
clearly see that the maximum value remains constant under this
type of refinement which suggests that the maximum overshoot
also remains constant for $q\neq 1$, as well as in the limit $q \rightarrow 1$.

\subsubsection{Gibbs Phenomenon on Meshes in One Dimension}
Next, we consider a one-dimensional example such that $u$ is not piecewise linear and compute the $L^q$-best approximation numerically.
Let
$u (x)= 1+0.1 \sin(2\pi x)$ on $(0,1)$ and consider the $L^q$-best approximation
$u_h$ with $u_h(0) = 1$ and $u_h(1) = 0$ on four different grids:
two uniform grids with 5 and 100 elements, respectively, and two meshes
where all elements are the same size except the last one which is twice
the size of the others. Again we consider a mesh with 5 elements and
one with 100 elements. Note that the latter two meshes violate the conditions in parts two and three of Theorem \ref{thm:lq_1d}, but satisfy the condition
in Remark \ref{rem:refinement}. We therefore expect the overshoot to vanish as $q \rightarrow 1$ in the first two cases and to decrease but still be present in the last two. Remark \ref{rem:refinement} and the observations for the previous example suggests that for $u\equiv 1$, we could expect the overshoot to be the same both when 5 and 100 elements are employed on both  the uniform and the non-uniform meshes.

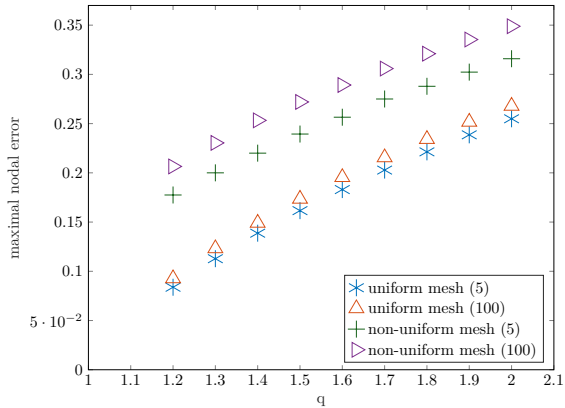
\begin{figure}[t]
  \centering
%
%
\definecolor{mycolor1}{rgb}{0.00000,0.44700,0.74100}%
\definecolor{mycolor2}{rgb}{0.85000,0.32500,0.09800}%
\definecolor{mycolor3}{RGB}{21, 97, 24}%
\definecolor{mycolor4}{rgb}{0.49400,0.18400,0.55600}%
\begin{tikzpicture}[thick, scale=0.4, every node/.style={scale=1.5}]

\begin{axis}[%
width=6.028in,
height=4.754in,
at={(1.011in,0.642in)},
scale only axis,
xmin=1.0,
xmax=2.1,
xlabel style={font=\color{white!15!black}},
xlabel={q},
ymin=0.0,
ymax=0.37,
ylabel style={font=\color{white!15!black}},
ylabel={maximal nodal error},
axis background/.style={fill=white},
legend style={at={(0.55,0.02)}, anchor=south west, legend cell align=left, align=left, draw=white!15!black, only marks}
]
\addplot [color=mycolor1, draw=none, mark=asterisk, mark options={solid, mycolor1, scale = 4.0}]
  table[row sep=crcr]{%
2	0.2549\\
1.9	0.2386\\
1.8	0.2213\\
1.7	0.2028\\
1.6	0.183\\
1.5	0.1616\\
1.4	0.1387\\
1.3	0.1128\\
1.2	0.0838\\
};
\addlegendentry{uniform mesh (5)}

\addplot [color=mycolor2, draw=none, mark=triangle, mark options={solid, mycolor2, scale = 4.0}]
  table[row sep=crcr]{%
2	0.2679\\
1.9	0.2517\\
1.8	0.2342\\
1.7	0.2156\\
1.6	0.1955\\
1.5	0.1735\\
1.4	0.1491\\
1.3	0.1232\\
1.2	0.0927\\
};
\addlegendentry{uniform mesh (100)}

\addplot [color=mycolor3, draw=none, mark=+, mark options={solid, mycolor3, scale = 4.0}]
  table[row sep=crcr]{%
2	0.3159\\
1.9	0.3023\\
1.8	0.2879\\
1.7	0.275\\
1.6	0.2565\\
1.5	0.2394\\
1.4	0.22\\
1.3	0.2\\
1.2	0.1775\\
};
\addlegendentry{non-uniform mesh (5)}

\addplot [color=mycolor4, draw=none, mark=triangle, mark options={solid, rotate=270, mycolor4, scale = 4.0}]
  table[row sep=crcr]{%
2	0.3489\\
1.9	0.3354\\
1.8	0.321\\
1.7	0.3059\\
1.6	0.2892\\
1.5	0.272\\
1.4	0.2533\\
1.3	0.2304\\
1.2	0.2065\\
};
\addlegendentry{non-uniform mesh (100)}

\end{axis}
\end{tikzpicture}%
  \caption{Maximal nodal error for different values of $q$ and four different meshes.}\label{fig:smooth_u}
\end{figure}

Fig.\ \ref{fig:smooth_u} shows the maximum error at the nodes in
all four cases for several values of $q$. We observe that
the overshoot indeed decreases as $q\rightarrow 1$. Furthermore, we see that the overshoot is very similar for the coarse and
fine meshes in both cases which confirms that the overshoot does not disappear under mesh refinement. However, the overshoot is not identical for $5$ and for $100$ elements in both cases which can be attributed to the fact that $u$ is not constant. Furthermore, note that the overshoot for the non-uniform mesh
is consistently larger than for the uniform mesh, which  suggests
that it does not disappear entirely as $q \rightarrow 1$. Note that on the non-uniform mesh when $u \equiv 1$ and $q=1$, the overshoot would be $\frac{2}{3}\sqrt{3}-1 \approx 0.15$; see Remark \ref{rem:refinement}.

\subsection{(Vanishing) Overshoot in One Dimension}\label{sec:examples_1d}

To illustrate the graded mesh condition in part two of Theorem \ref{thm:lq_1d}, we consider two three-element meshes on $(0,1)$. For the first one we choose $h_1 = 0.1$ and $h_2 = h_3 = 0.45$, i.e., the mesh consists of the subintervals $(0,0.1)$, $(0.1,0.55)$ and $(0.55, 1)$.
For the second one we choose $h_1 = 0.1$, $h_2 = 0.5$ and $h_3 = 0.4$, i.e., the mesh consisting of the subintervals $(0,0.1)$, $(0.1,0.6)$ and $(0.6, 1)$. We will check the condition \eqref{eq:sufficient_cond_1d} for both meshes; Indeed, we will see that for the first mesh the condition is violated but it is satisfied for the second mesh. In the latter case, we therefore know that there exists an $L^1$-best approximation without over- or undershoots.
In the former case, it is a priori unknown whether or not such an $L^1$-best approximation exists, since it is an open problem whether \eqref{eq:sufficient_cond_1d} is also a necessary condition.

In the first case, we obtain from \eqref{eq:def_vartheta} that $\vartheta_3 = \vartheta_2 = 0$ yielding the following sufficient conditions for the existence of an $L^1$-best approximation without over- or undershoots: $h_2 \geq h_3$ and $h_1 \geq h_2$. The second condition is violated.
In fact, it is easy to show that, if $h_2 =h_3$, the condition $h_1\geq h_2$ is necessary for the existence of an $L^1$-best approximation without over- or undershoots. Moreover, one can show that the $L^1$-best approximation is unique in this case by solving \eqref{eq:best_approx_l1} for the points where $u$ and $u_h$ intersect. The intersection points uniquely determine $u_h(0.1)\approx 0.9931$ and $u_h(0.55) \approx 1.0247$. For brevity, the details are omitted here.

For the second mesh, we again have $\vartheta_3 = 0$ and \eqref{eq:sufficient_cond_1d} with $i = 2$ becomes $h_2 \geq h_3$, which holds for $h_2 = 0.5$ and $h_3=0.4$. For $i = 2$, we obtain from \eqref{eq:def_vartheta} that
	$\vartheta_2^2 = 0.1$.
Hence,  $\vartheta_2 >1- \frac{1}{\sqrt{2}}$ and there is no condition on $h_1$ according to Theorem \ref{thm:lq_1d}. Therefore, there exists an $L^1$-best approximation without over- or undershoots.

\begin{figure}
	\centering
	\begin{tikzpicture}[scale = 1.0]

		\begin{axis}[%
		width=5cm,
		height=4cm,
		at={(1.011in,0.642in)},
		scale only axis,
		xmin=0,
		xmax=1,
		xlabel style={font=\color{white!15!black}},
		xlabel={x},
		ymin=0,
		ymax=1.31,
		ylabel style={font=\color{white!15!black}},
		ylabel={$u_h(x)$},
		axis background/.style={fill=white},
		legend style={at={(0.0,0.0)}, anchor=south west, legend cell align=left, align=left, draw=white!15!black}
		]

		\addplot [color=black, mark=none, forget plot]
			table[row sep=crcr]{%
			0 1.0\\
			0.25 1.000\\
			0.5 1.0\\
			0.75 1.0\\
			1 1.0\\
		};
		\addplot [color=red, mark=triangle]
			table[row sep=crcr]{%
		0 1.0\\
		0.1 0.886075949367\\
		0.55 1.27848101266\\
		1 0 \\
		};
		\addlegendentry{$q=2$, $h_2 = 0.45$}

		\addplot [color=blue, mark=asterisk]
			table[row sep=crcr]{%
			0 1.0\\
			0.1 0.9931\\
			0.55 1.0247\\
			1 0 \\
		};
		\addlegendentry{$q=1$, $h_2 = 0.45$}
		\addplot [color=red!50!yellow, mark=triangle, mark options = {rotate=270}]
			table[row sep=crcr]{%
		0 1.0\\
		0.1 0.895287958115\\
		0.6 1.25130890052\\
		1 0 \\
		};
		\addlegendentry{$q=2$, $h_2 = 0.5$}

		\addplot [color=green!60!black, mark=+]
			table[row sep=crcr]{%
			0 1.0\\
			0.1 1.0\\
			0.6 1.0\\
			1 0 \\
		};
		\addlegendentry{$q=1$, $h_2 = 0.5$}
		\end{axis}
\end{tikzpicture}
%
%
\definecolor{mycolor1}{rgb}{0.00000,0.44700,0.74100}%
\definecolor{mycolor2}{rgb}{0.85000,0.32500,0.09800}%
\definecolor{mycolor3}{RGB}{21, 97, 24}%
\definecolor{mycolor4}{rgb}{0.49400,0.18400,0.55600}%
\begin{tikzpicture}[thick, scale=1.0, every node/.style={scale=1.0}]

\begin{axis}[%
width=5cm,
height=4cm,
at={(1.011in,0.642in)},
scale only axis,
xmin=1.0,
xmax=2.1,
xlabel style={font=\color{white!15!black}},
xlabel={q},
ymin=0.0,
ymax=0.30,
ylabel style={font=\color{white!15!black}},
ylabel={maximal nodal error},
axis background/.style={fill=white},
legend style={at={(0.55,0.02)}, anchor=south west, legend cell align=left, align=left, draw=white!15!black, only marks}
]
\addplot [color=blue, draw=none, mark=asterisk, mark options={solid, blue, scale = 1.0}]
  table[row sep=crcr]{%
2	0.278481012658\\
1.9	0.261569221008\\
1.8	0.243023311461\\
1.7	0.223214657052\\
1.6	0.2054817613\\
1.5	0.195872547251\\
1.4	0.180713603317\\
1.3	0.159296647375\\
1.2	0.122984344462\\
1.0 0.0247\\
};
\addlegendentry{$h_2 = 0.45$}

\addplot [color=red, draw=none, mark=triangle, mark options={solid, red, scale = 1.0}]
  table[row sep=crcr]{%
2	0.251308900524\\
1.9	0.233823117755\\
1.8	0.21543733823\\
1.7	0.200058567538\\
1.6	0.184084910766\\
1.5	0.1653851901\\
1.4	0.141926486021\\
1.3	0.109658068939\\
1.2	0.0636077124789\\
1.0 0.0\\
};
\addlegendentry{$h_2 = 0.5$}

\end{axis}
\end{tikzpicture}%
	\caption[$L^q$-best approximations on two three-element meshes on $(0,1)$]{$L^q$-best approximations on two three-element meshes on $(0,1)$ with $h_1 = 0.1$ and two different choices for $h_2$. Left: best approximation with $q=2$ and $q=1$. Right: maximal nodal error for several values of $q$. }
	\label{fig:3el}
\end{figure}
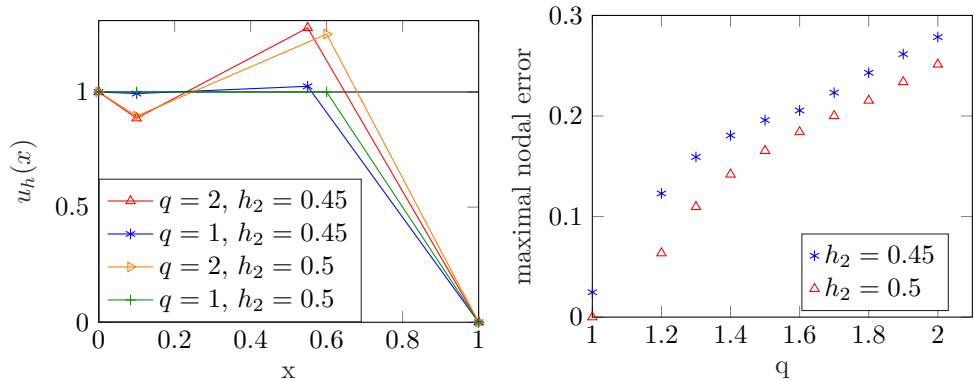
Fig.\ \ref{fig:3el} shows the $L^q$-best approximation on both meshes for $q = 2$ and $q=1$ on the left and the maximal nodal error on both meshes for several values of $q$ on the right. The approximations for $q>1$ were again obtained using the implementation of the best approximation problem in FEniCS. We can clearly see, that the maximal overshoot is always larger on the first mesh. In both cases it decreases as $q\rightarrow 1$, but the overshoot only vanishes completely on the second mesh. However, even on the first mesh the maximal overshoot is very small for $q = 1$. Note that, if $h_2$ and $h_1$ as chosen for the first mesh were swapped, the maximal overshoot for $q=1$ would be $u_h(0.55)-1 = 0.2792$ according to Remark \ref{rem:refinement} and thus significantly larger than the overshoot we can observe. This shows that  the effect of an element being too small and causing the $L^1$-best approximation to contain over- and undershoots is much weaker away from the discontinuity than near the discontinuity.

\subsection{(Vanishing) Overshoot in Two Dimensions}\label{sec:examples_2d}
\subsubsection{Overshoot on Meshes 2, 3 and 4 from Section \ref {sec:2d_l1_rest}}\label{sec:2d_l1_rest_numerics}

\begin{figure}
 \centering
\begin{subfigure}{4.7cm}
   \centering
\includegraphics[width = 4.7cm]{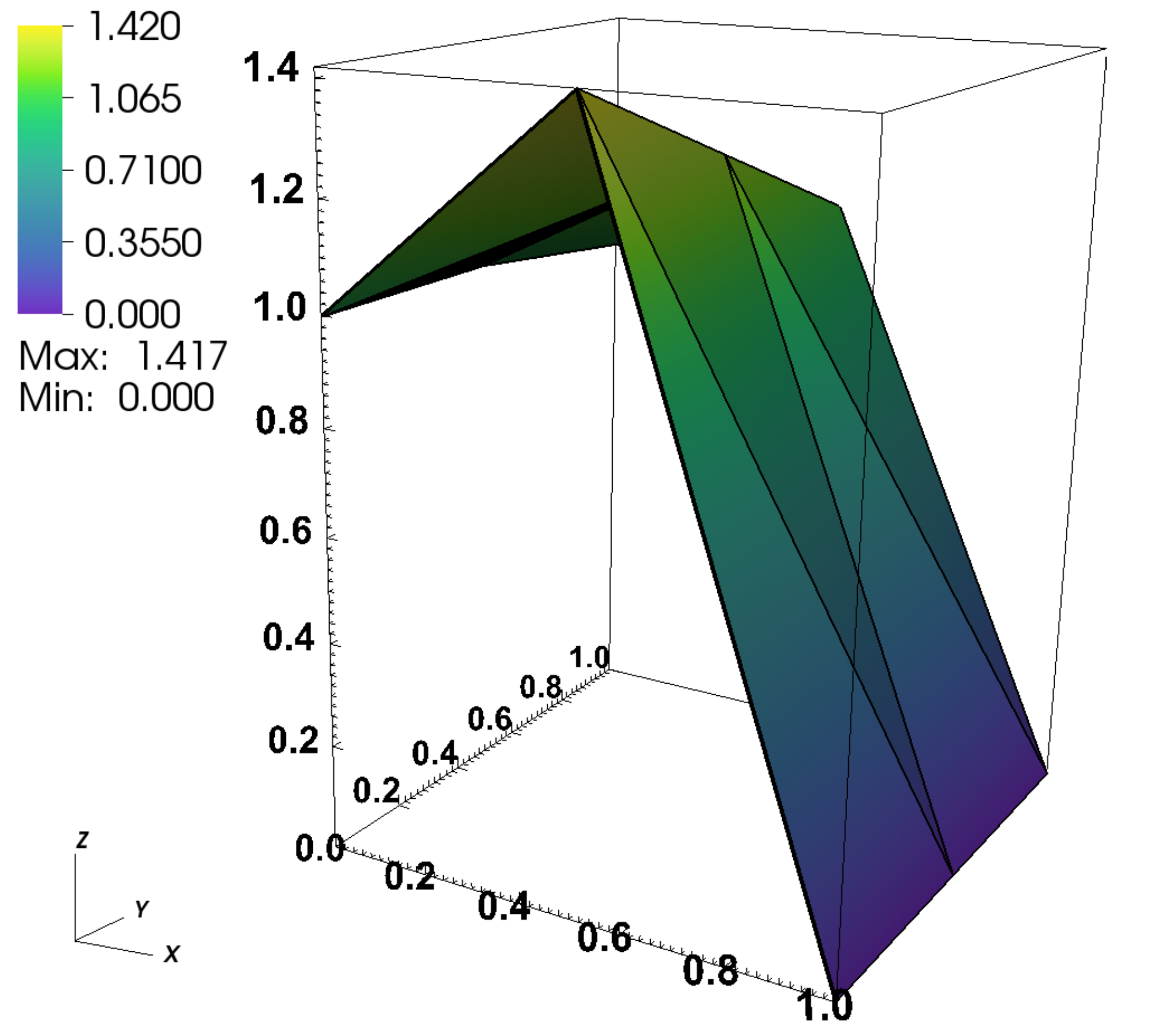}
\caption{Mesh 2, $q=2$}
\end{subfigure}
\begin{subfigure}{4.7cm}
   \centering
\includegraphics[width = 4.7cm]{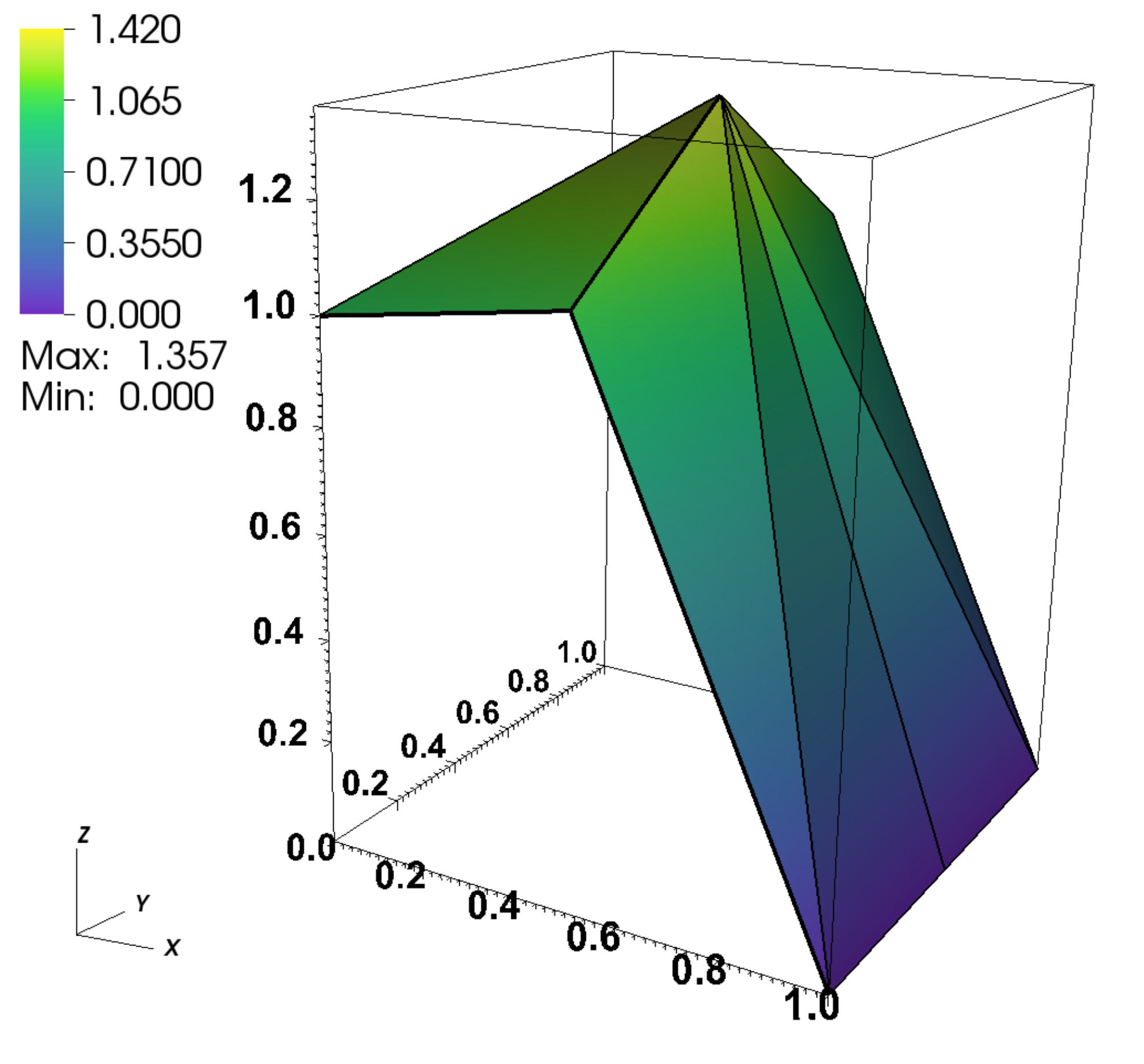}
\caption{Mesh 3, $q=2$}
\end{subfigure}
\begin{subfigure}{4.7cm}
   \centering
\includegraphics[width = 4.7cm]{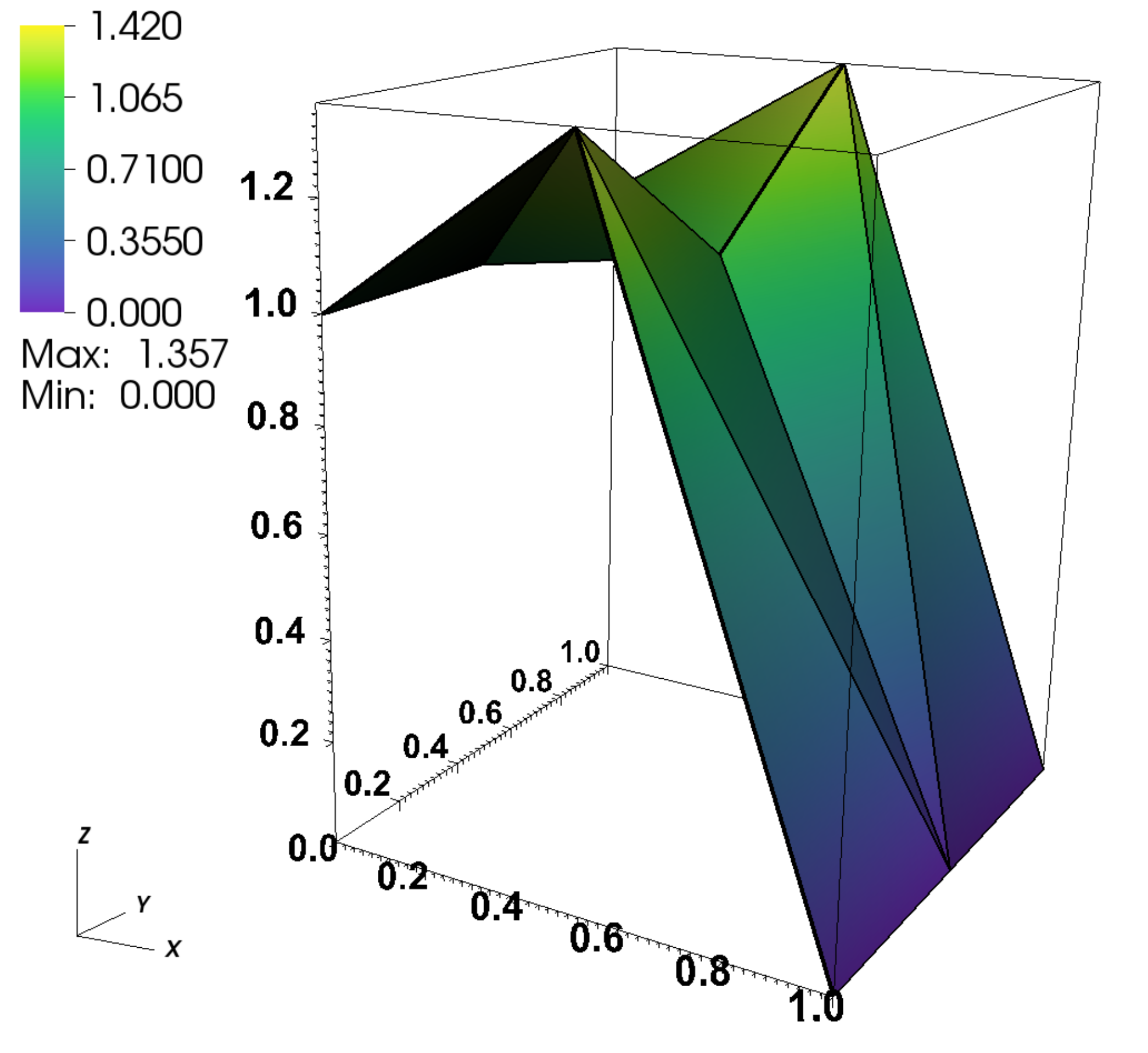}
\caption{Mesh 4, $q=2$}
\end{subfigure}
\begin{subfigure}{4.7cm}
   \centering
\includegraphics[width = 4.7cm]{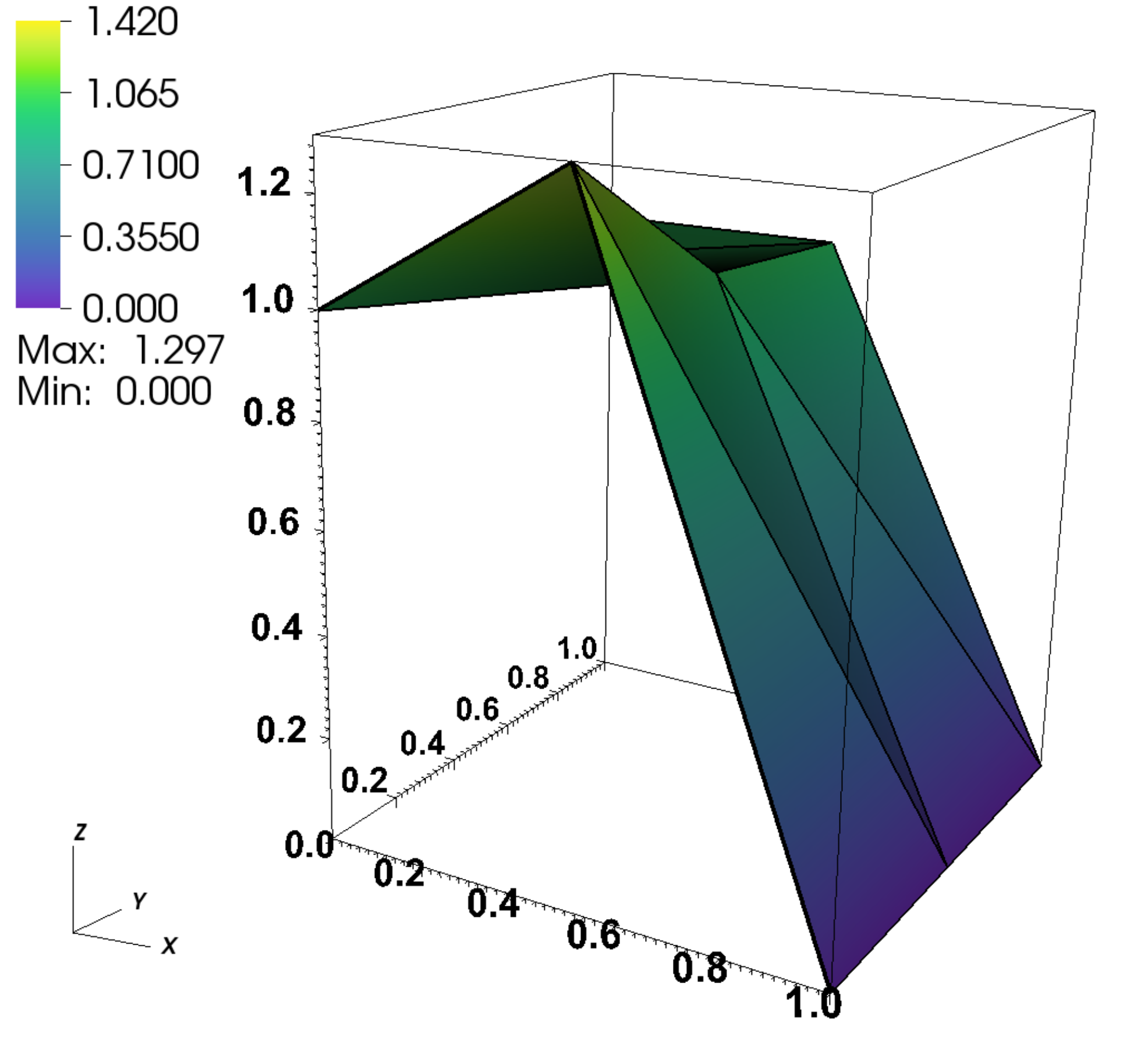}
\caption{Mesh 2, $q=1.1$}
\end{subfigure}
\begin{subfigure}{4.7cm}
   \centering
\includegraphics[width = 4.7cm]{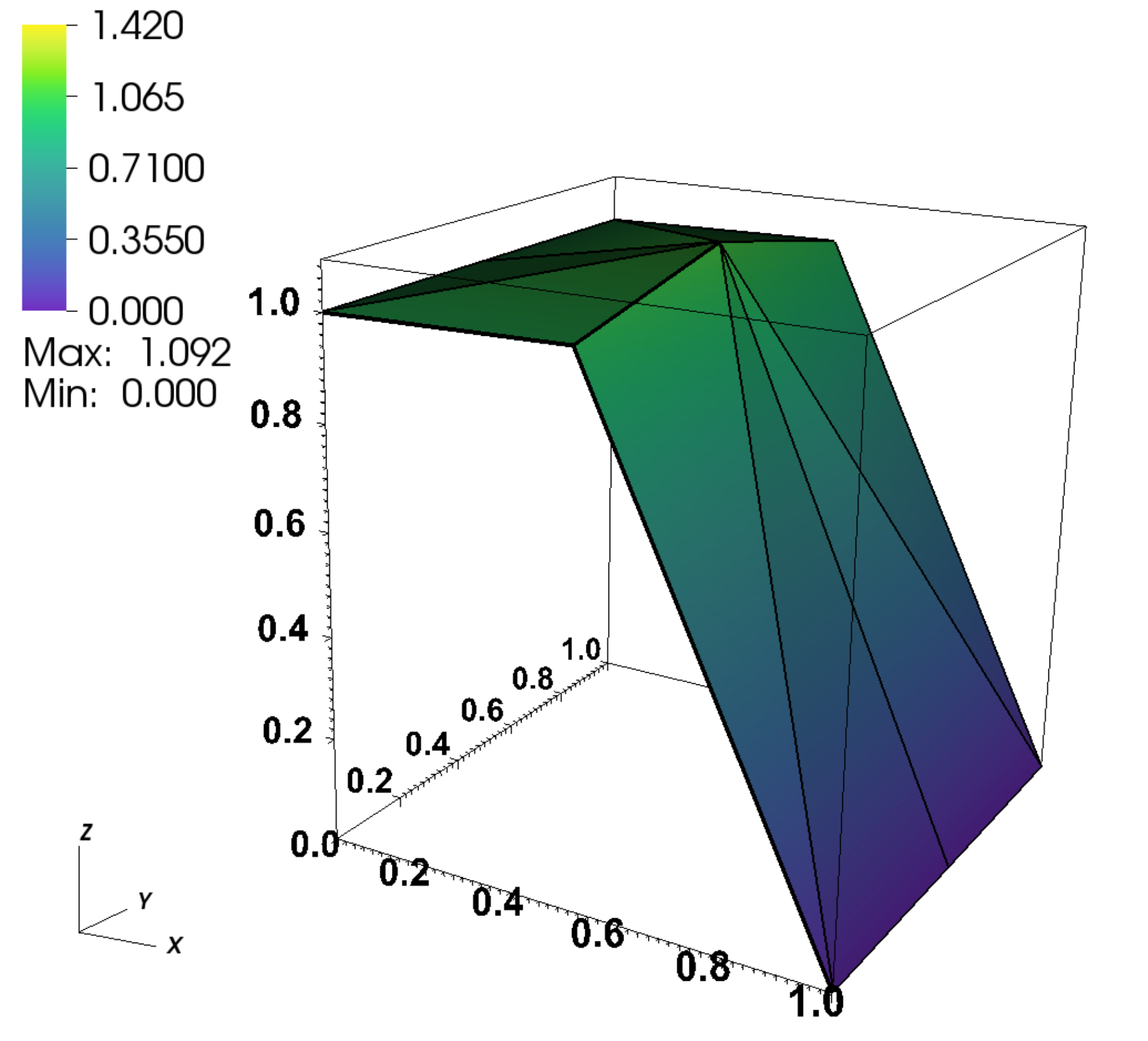}
\caption{Mesh 3, $q=1.1$}
\end{subfigure}
\begin{subfigure}{4.7cm}
   \centering
\includegraphics[width = 4.7cm]{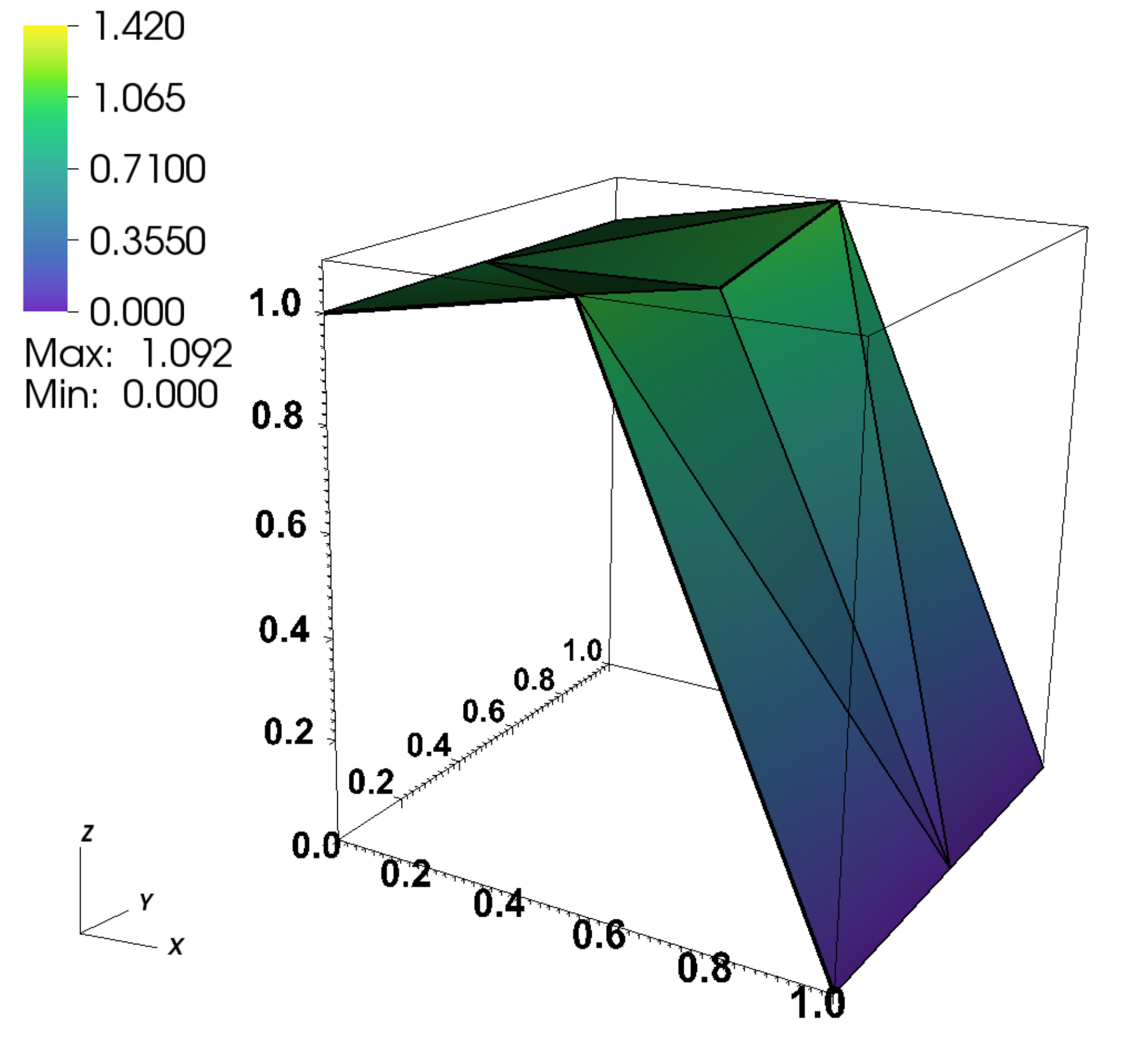}
\caption{Mesh 4, $q=1.1$}
\end{subfigure}
\caption{$L^q$-best approximation of a boundary discontinuity}
\label{fig:surface_plots}
\end{figure}

Fig.\ \ref{fig:surface_plots} shows the best approximations for $q = 2$ and $q=1.1$ for three of the
meshes we have considered in Section \ref{sec:2d_l1_rest}. Even just a comparison of these two cases for each of the meshes
illustrates clearly how the overshoot gradually vanishes on Mesh 3 and Mesh 4. On Mesh 2,
the overshoot vanishes away from the boundary $y = 0$; this is consistent with the $L^1$-best
approximation described above that only exhibits an overshoot at the node $(0.5,0)$ and
no overshoot at all other nodes.

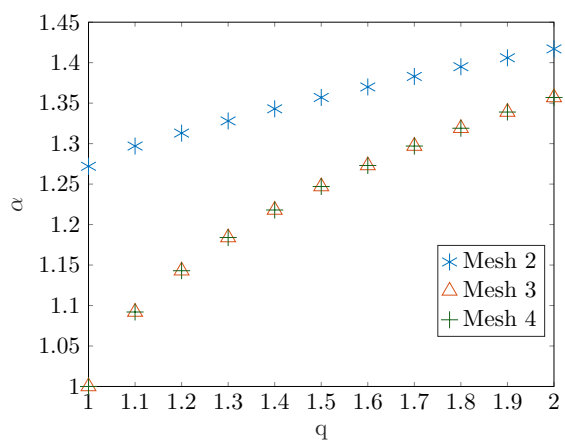
\begin{figure}
 \centering
%
%
\definecolor{mycolor1}{rgb}{0.00000,0.44700,0.74100}%
\definecolor{mycolor2}{rgb}{0.85000,0.32500,0.09800}%
\definecolor{mycolor3}{RGB}{21, 97, 24}
\begin{tikzpicture}[thick, scale=0.4, every node/.style={scale=2.2}]

\begin{axis}[%
width=6.028in,
height=4.754in,
at={(1.011in,0.642in)},
scale only axis,
xmin=1,
xmax=2,
xlabel style={font=\color{white!15!black}},
xlabel={q},
ymin=1,
ymax=1.45,
ylabel style={font=\color{white!15!black}},
ylabel={$\alpha$},
axis background/.style={fill=white},
legend style={at={(0.75,0.137)}, anchor=south west, legend cell align=left, align=left, draw=white!15!black, only marks}
]
\addplot [color=mycolor1, draw=none, mark=asterisk, mark options={solid, mycolor1, scale = 4.0}]
  table[row sep=crcr]{%
1	1.272\\
1.1	1.297\\
1.2	1.313\\
1.3	1.328\\
1.4	1.343\\
1.5	1.357\\
1.6	1.37\\
1.7	1.383\\
1.8	1.395\\
1.9	1.406\\
2	1.417\\
};
\addlegendentry{Mesh 2}

\addplot [color=mycolor2, draw=none, mark=triangle, mark options={solid, mycolor2, scale = 4.0}]
  table[row sep=crcr]{%
1	1\\
1.1	1.092\\
1.2	1.143\\
1.3	1.184\\
1.4	1.218\\
1.5	1.247\\
1.6	1.273\\
1.7	1.297\\
1.8	1.319\\
1.9	1.339\\
2	1.357\\
};
\addlegendentry{Mesh 3}

\addplot [color=mycolor3, draw=none, mark=+, mark options={solid, mycolor3, scale = 4.0}]
  table[row sep=crcr]{%
1	1\\
1.1	1.092\\
1.2	1.143\\
1.3	1.184\\
1.4	1.218\\
1.5	1.247\\
1.6	1.273\\
1.7	1.297\\
1.8	1.319\\
1.9	1.339\\
2	1.357\\
};
\addlegendentry{Mesh 4}

\end{axis}
\end{tikzpicture}%
 \caption{Values for $\alpha$ for different ranges of $q$.}\label{fig:alpha_2D}
\end{figure}

Fig.\ \ref{fig:alpha_2D} shows the maximum overshoot  for all three meshes for different values of $q$. The overshoot
for $q = 1$ is taken from the theoretically determined $L^1$-best approximations discussed in that
section. All remaining values have been determined numerically with an implementation
in FEniCS \citep{Alnes2015}. The plot shows that for the third and fourth mesh, the overshoot indeed
disappears as $q \rightarrow 1$, whereas for the second mesh it decreases but does not
vanish.

\subsubsection{Overshoot on Unstructured Meshes}
\begin{figure}[t!]
 \centering
 \begin{subfigure}{4.7cm}
    \centering
 \includegraphics[width = 4.7cm]{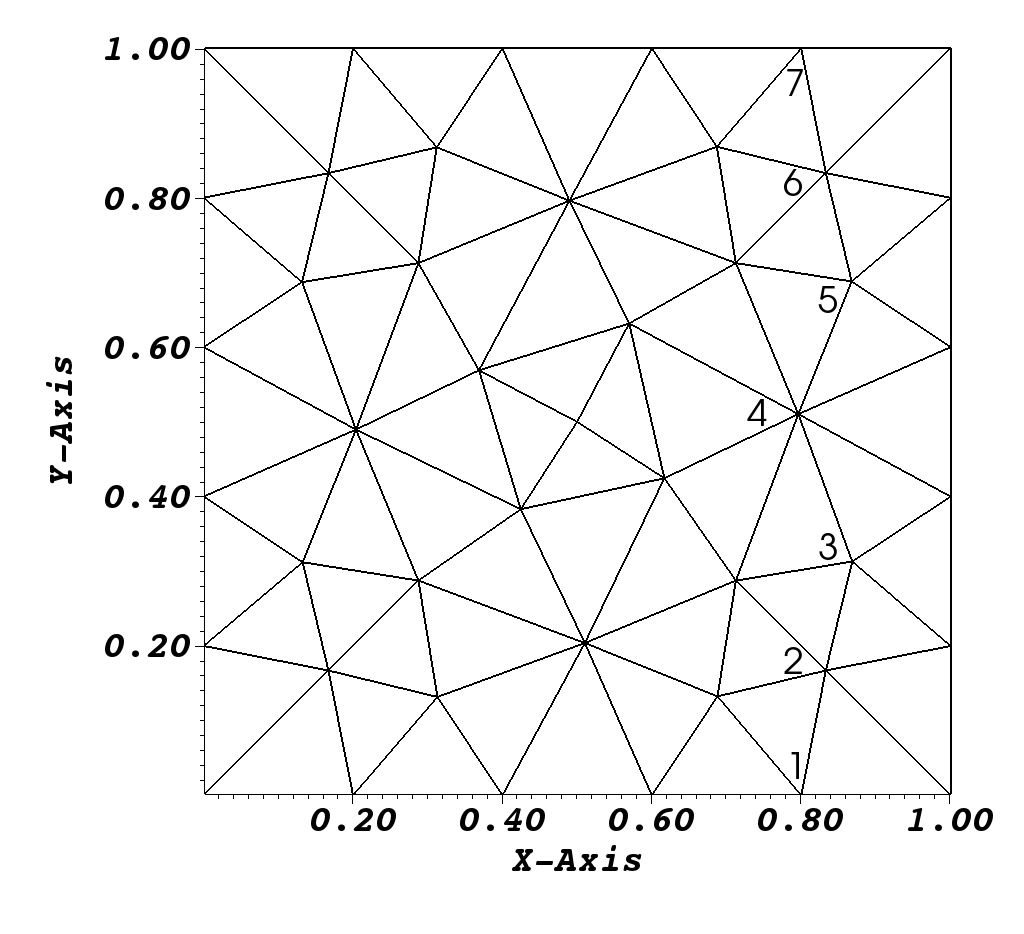}
 \caption{Mesh A}
 \end{subfigure}
\begin{subfigure}{4.7cm}
   \centering
\includegraphics[width = 4.7cm]{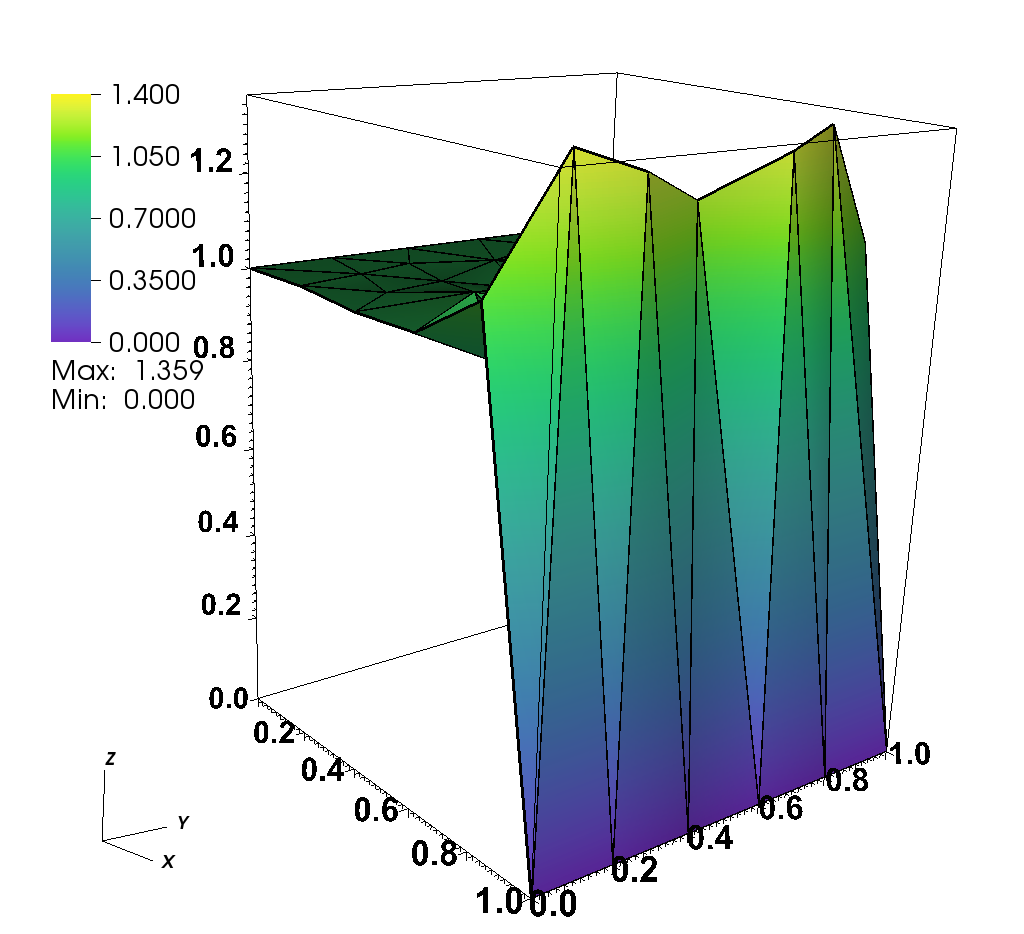}
\caption{ $q=2$}
\end{subfigure}
\begin{subfigure}{4.7cm}
   \centering
\includegraphics[width = 4.7cm]{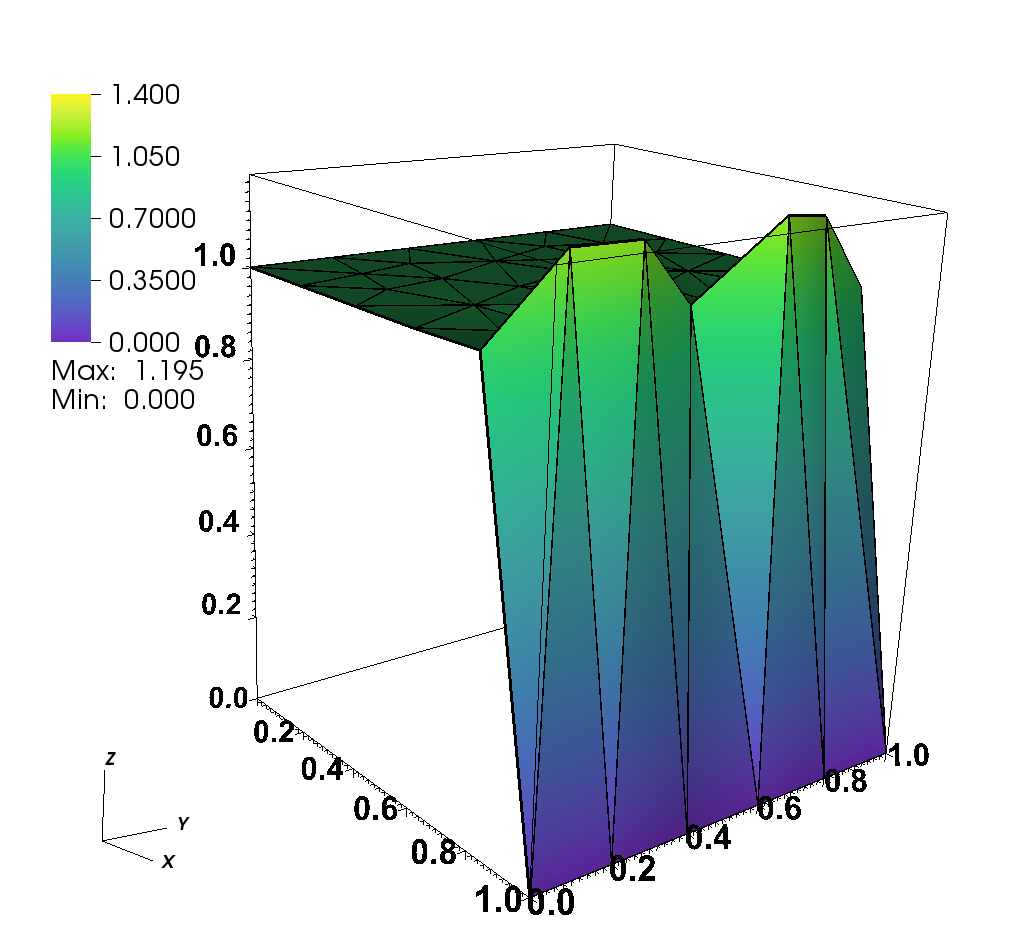}
\caption{ $q=1.2$}
\end{subfigure}
\caption{$L^q$-best approximation on an unstructured mesh.}
\label{fig:unstructerd_mesh}
\end{figure}

As a final example, we consider the unstructured mesh shown on the left in
Fig.\ \ref{fig:unstructerd_mesh}. From the computations for the
previous examples, we deduce that the $L^1$-best approximation
exhibits no overshoot if for every interior node $(x_i,y_i)$ that is connected
to the boundary $x = 1$ through one edge, the area of all triangles
whose boundaries contain the node $(x_i,y_i)$ \emph{and} at least one node on the boundary $x=1$ is
smaller than the area of all remaining triangles whose boundaries contain the node $(x_i,y_i)$.
Furthermore, the numerics for Mesh 2 (cf., Fig.\ \ref{fig:meshes_2d})
shown in Fig.\ \ref{fig:surface_plots} suggest that the overshoot
disappears for $q\rightarrow 1$ away from the nodes violating this
condition on the area of connected elements.

The interior nodes connected to the boundary are labelled $1,2, \dots 7$
in Fig.\ \ref{fig:unstructerd_mesh}. It is easy to see that for the
nodes $1,4$ and $7$ the total area of all triangles
touching both the node and the boundary $x=1$ is smaller than
the total area of all remaining triangles touching the node, whereas
this condition is violated for the nodes $2,3,5$ and $6$. Thus, we
expect that the overshoot vanishes at the nodes $1,4$ and $7$ as $q \rightarrow 1$,
while at the nodes $2,3,5$ and $6$ it reduces but does not disappear
entirely. Fig.\ \ref{fig:unstructerd_mesh} shows the
$L^q$-best approximation on the unstructured mesh for $q = 2$ and
$q=1.2$ in the center and on the right, respectively. Here, we clearly observe that the approximation for $q=2$ exhibits
overshoots at all nodes connected to the boundary $x=1$ with larger
overshoots at the nodes $2,3,5$ and $6$. At these nodes the
overshoot is reduced but still clearly visible for $q=1.2$. On the other hand at the nodes $1,4$ and $7$ the overshoot has nearly vanished for
$q=1.2$.

\begin{figure}
 \centering
 \begin{subfigure}{4cm}
    \centering
\includegraphics[width = 4cm]{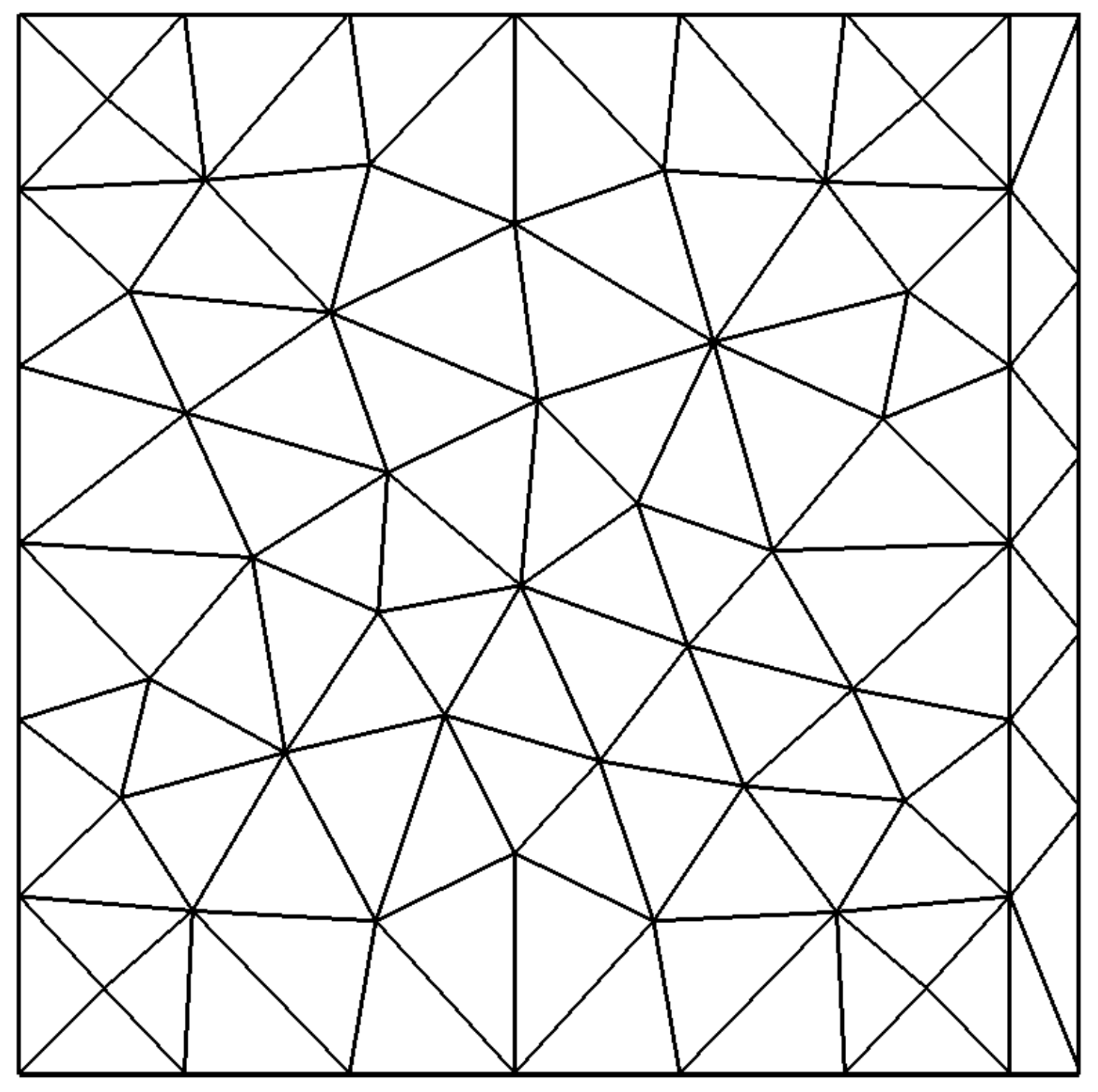}
 \caption{Mesh B}
 \end{subfigure}
 \begin{subfigure}{4cm}
    \centering
 \includegraphics[width = 4cm]{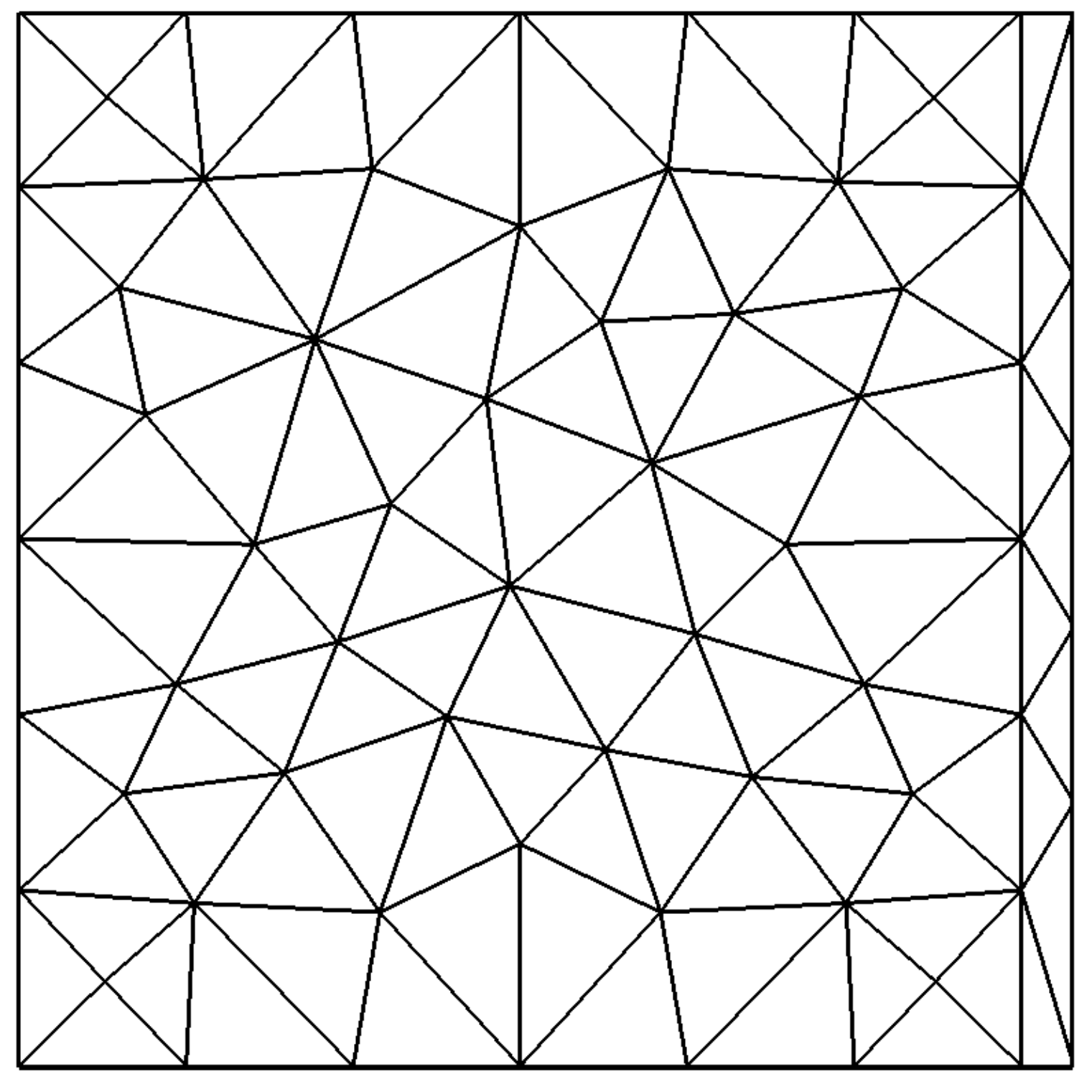}
 \caption{Mesh C}
 \end{subfigure}
\begin{subfigure}{6.1cm}
   \centering
%
%
\definecolor{mycolor1}{rgb}{0.00000,0.44700,0.74100}%
\definecolor{mycolor2}{rgb}{0.85000,0.32500,0.09800}%
\definecolor{mycolor3}{RGB}{21, 97, 24}
\begin{tikzpicture}[thick, scale=0.35, every node/.style={scale=1.2}]

\begin{axis}[%
width=6.028in,
height=4.754in,
at={(1.011in,0.642in)},
scale only axis,
xmin=0.98,
xmax=2.02,
xlabel style={font=\color{white!15!black}},
xlabel={q},
ymin=0.94,
ymax=1.37,
ylabel style={font=\color{white!15!black}},
ylabel={$\text{max(u}_\text{h}\text{)}$},
axis background/.style={fill=white},
legend style={at={(0.02,0.82)}, anchor=south west, legend cell align=left, align=left, draw=white!15!black, only marks}
]
\addplot [color=mycolor1, draw=none, mark=asterisk, mark options={solid, mycolor1, scale = 4.0}]
  table[row sep=crcr]{%
2	1.359\\
1.9	1.342\\
1.8	1.325\\
1.7	1.312\\
1.6	1.286\\
1.5	1.263\\
1.4	1.239\\
1.3	1.216\\
};
\addlegendentry{Mesh A}

\addplot [color=mycolor2, draw=none, mark=triangle, mark options={solid, mycolor2, scale = 4.0}]
  table[row sep=crcr]{%
2	1.274\\
1.9	1.253\\
1.8	1.231\\
1.7	1.207\\
1.6	1.182\\
1.5	1.155\\
1.4	1.124\\
1.3	1.089\\
};
\addlegendentry{Mesh B}

\addplot [color=mycolor3, draw=none, mark=+, mark options={solid, mycolor3, scale = 4.0}]
  table[row sep=crcr]{%
2	1.226\\
1.9	1.205\\
1.8	1.183\\
1.7	1.159\\
1.6	1.133\\
1.5	1.106\\
1.4	1.074\\
1.3	1.043\\
};
\addlegendentry{Mesh C}

\end{axis}
\end{tikzpicture}%
\caption{$\max(u_h)$}\label{fig:max_uh}
\end{subfigure}
\caption{$L^1$-best approximation unstructured mesh with vanishing overshoot}
\label{fig:unstructerd_mesh_mod}
\end{figure}
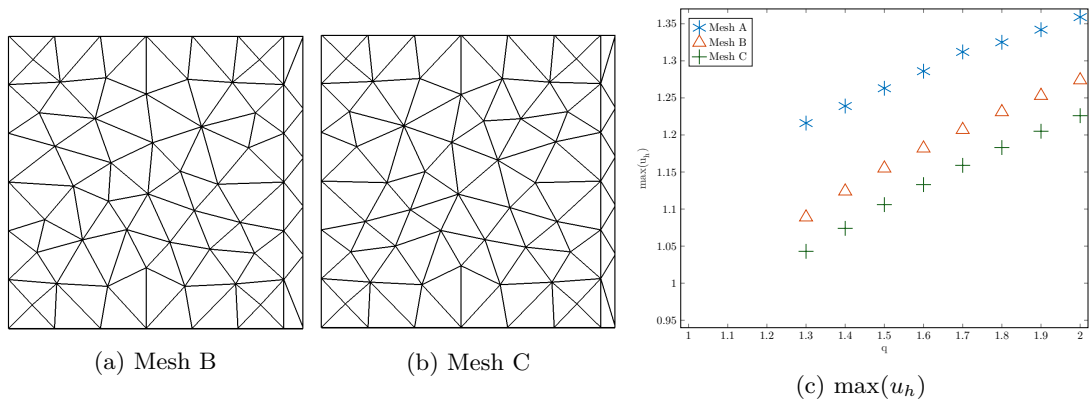

Fig.\ \ref{fig:unstructerd_mesh_mod} shows two further unstructured meshes which have been designed in such a way that for
every interior node $(x_i,y_i)$ that is connected
to the boundary $x = 1$ through one edge, the area of all triangles
whose boundaries contain the node $(x_i,y_i)$ \emph{and} at least one node on the boundary $x=1$ is
smaller than the area of all remaining triangles whose boundaries contain the node $(x_i,y_i)$. The difference between the two meshes is that the distance between the boundary $x = 1$ and the vertical line containing all nodes connected to this boundary is smaller in Mesh C than in Mesh B. Fig.\ \ref{fig:max_uh}
shows the maximum value of $u_h$
for different $q$ and Meshes A, B and C. This illustrates that the overshoot decreases on all three meshes as $q \rightarrow 1$. The overshoot on Mesh C is always smaller than the overshoot on the other two meshes and the overshoot on Mesh A is always larger than on the other two meshes. This illustrates that if the area  of the elements connected to the boundary is decreased in comparison the area of the remaining elements, then the overshoot is reduced for any $q$ and decreases more rapidly as $q \rightarrow 1$. This is consistent with the theoretical results in one dimension illustrated at the start of this article in Fig.\ \ref{fig:alpha_one_d}.

\section{Conclusions}
 In this article, we have investigated Gibbs phenomena in the $L^q$-best approximation of discontinuities within finite element spaces. Using selected examples, we have proven that the Gibbs phenomenon can be eliminated as $q \rightarrow 1$ on certain meshes. However, we have seen that  there exist non-uniform meshes in one dimension that lead to Gibbs phenomena even if $q = 1$. In two dimensions, even some uniform meshes lead to Gibbs phenomena if $q=1$. Nonetheless, the magnitude of the oscillations decreases as $q \rightarrow 1$ on all meshes.

 The computational examples presented in this article confirm the theoretical results. Moreover, we have seen that similar observations can be made for more general examples. Furthermore, we have demonstrated that Gibbs phenomenon cannot be eliminated on certain meshes under mesh refinement.
 For the final computational example, we have been able to establish a link between the structure of the mesh near the discontinuity and the magnitude of the overshoot at the nodes. This observation suggests that the oscillations can be eliminated in the limit as $q$ tends to $1$ if the mesh structure near the discontinuity is suitably adjusted. Indeed, this has been used to design meshes for the non-linear Petrov-Galerkin method for the convection-diffusion-reaction equation presented in \cite{Houston2019b}.

\bibliography{main}

\begin{thebibliography}{29}
\providecommand{\natexlab}[1]{#1}
\providecommand{\url}[1]{\texttt{#1}}
\expandafter\ifx\csname urlstyle\endcsname\relax
  \providecommand{\doi}[1]{doi: #1}\else
  \providecommand{\doi}{doi: \begingroup \urlstyle{rm}\Url}\fi

\bibitem[Aln{\ae}s et~al.(2015)Aln{\ae}s, Blechta, Hake, Johansson, Kehlet,
  Logg, Richardson, Ring, Rognes, and Wells]{Alnes2015}
M.~S. Aln{\ae}s, J.~Blechta, J.~Hake, A.~Johansson, B.~Kehlet, A.~Logg,
  C.~Richardson, J.~Ring, M.~E. Rognes, and G.~N. Wells.
\newblock The {FEniCS} project version 1.5.
\newblock \emph{Archive of Numerical Software}, 3\penalty0 (100), 2015.
\newblock \doi{10.11588/ans.2015.100.20553}.

\bibitem[Bank and Yserentant(2014)]{Bank2014}
R.~E. Bank and H.~Yserentant.
\newblock On the {$H^1$}-stability of the {$L_2$}-projection onto finite
  element spaces.
\newblock \emph{Numer. Math.}, 126\penalty0 (2):\penalty0 361--381, 2014.
\newblock \doi{10.1007/s00211-013-0562-4}.

\bibitem[Cioranescu(1990)]{Cioranescu1990}
I.~Cioranescu.
\newblock \emph{Geometry of {B}anach spaces, duality mappings and nonlinear
  problems}, volume~62 of \emph{Mathematics and its Applications}.
\newblock Kluwer Academic Publishers Group, Dordrecht, 1990.
\newblock \doi{10.1007/978-94-009-2121-4}.

\bibitem[Crouzeix and Thom\'{e}e(1987)]{Crouzeix1987}
M.~Crouzeix and V.~Thom\'{e}e.
\newblock The stability in {$L_p$} and {$W^1_p$} of the {$L_2$}-projection onto
  finite element function spaces.
\newblock \emph{Math. Comp.}, 48\penalty0 (178):\penalty0 521--532, 1987.
\newblock \doi{10.2307/2007825}.

\bibitem[Demkowicz and Gopalakrishnan(2014)]{Demkowicz2014}
L.~F. Demkowicz and J.~Gopalakrishnan.
\newblock An overview of the discontinuous {P}etrov {G}alerkin method.
\newblock In \emph{Recent developments in discontinuous {G}alerkin finite
  element methods for partial differential equations}, volume 157 of \emph{IMA
  Vol. Math. Appl.}, pages 149--180. Springer, Cham, 2014.
\newblock \doi{10.1007/978-3-319-01818-8_6}.

\bibitem[Douglas et~al.(1975)Douglas, Dupont, and Wahlbin]{Douglas1975}
J.~Douglas, Jr., T.~Dupont, and L.~Wahlbin.
\newblock Optimal {$L_{\infty }$} error estimates for {G}alerkin approximations
  to solutions of two-point boundary value problems.
\newblock \emph{Math. Comp.}, 29:\penalty0 475--483, 1975.
\newblock \doi{10.2307/2005565}.

\bibitem[Gibbs(1899)]{Gibbs1899}
J.~W. Gibbs.
\newblock Fourier's series.
\newblock \emph{Nature}, 59\penalty0 (1539):\penalty0 606, 1899.

\bibitem[Guermond(2004)]{Guermond2004}
J.-L. Guermond.
\newblock A finite element technique for solving first-order {PDEs} in {Lp}.
\newblock \emph{SIAM Journal on Numerical Analysis}, 42\penalty0 (2):\penalty0
  714--737, 2004.
\newblock \doi{10.1137/S0036142902417054}.

\bibitem[Guermond and Popov(2007)]{Guermond2007}
J.-L. Guermond and B.~Popov.
\newblock Linear advection with ill-posed boundary conditions via
  {$L^1$}-minimization.
\newblock \emph{Int. J. Numer. Anal. Model.}, 4\penalty0 (1):\penalty0 39--47,
  2007.

\bibitem[Guermond and Popov(2008/09)]{Guermond2008/09}
J.-L. Guermond and B.~Popov.
\newblock {$L^1$}-approximation of stationary {H}amilton-{J}acobi equations.
\newblock \emph{SIAM Journal on Numerical Analysis.}, 47\penalty0 (1):\penalty0
  339--362, 2008/09.
\newblock \doi{10.1137/070681922}.

\bibitem[Guermond and Popov(2009)]{Guermond2009}
J.-L. Guermond and B.~Popov.
\newblock An optimal {$L^1$}-minimization algorithm for stationary
  {H}amilton-{J}acobi equations.
\newblock \emph{Commun. Math. Sci.}, 7\penalty0 (1):\penalty0 211--238, 2009.
\newblock \doi{10.4310/CMS.2009.v7.n1.a11}.

\bibitem[Guermond et~al.(2008)Guermond, Marpeau, and Popov]{Guermond2008}
J.-L. Guermond, F.~Marpeau, and B.~Popov.
\newblock A fast algorithm for solving first-order {PDE}s by
  {$L^1$}-minimization.
\newblock \emph{Commun. Math. Sci.}, 6\penalty0 (1):\penalty0 199--216, 2008.
\newblock \doi{10.4310/CMS.2008.v6.n1.a10}.

\bibitem[Houston et~al.(2019)Houston, Roggendorf, and van~der
  Zee]{Houston2019b}
P.~Houston, S.~Roggendorf, and K.~G. van~der Zee.
\newblock {E}liminating {G}ibbs {P}henomena: A {N}on-linear {P}etrov-{G}alerkin
  {M}ethod for the {C}onvection-{D}iffusion-{R}eaction {E}quation.
\newblock \emph{arXiv preprint arXiv:1908.00996}, 2019.
\newblock URL \url{https://arxiv.org/pdf/1908.00996.pdf}.

\bibitem[Jiang(1993)]{Jiang1993}
B.-N. Jiang.
\newblock Non-oscillatory and non-diffusive solution of convection problems by
  the iteratively reweighted least-squares finite element method.
\newblock \emph{Journal of computational physics}, 105\penalty0 (1):\penalty0
  108--121, 1993.
\newblock \doi{10.1006/jcph.1993.1057}.

\bibitem[Jiang(1998)]{Jiang1998}
B.-N. Jiang.
\newblock \emph{The least-squares finite element method: theory and
  applications in computational fluid dynamics and electromagnetics}.
\newblock Springer Science \& Business Media, 1998.
\newblock \doi{10.1007/978-3-662-03740-9}.

\bibitem[John and Knobloch(2007{\natexlab{a}})]{John2007}
V.~John and P.~Knobloch.
\newblock On the performance of {SOLD} methods for convection-diffusion
  problems with interior layers.
\newblock \emph{Int. J. Comput. Sci. Math.}, 1\penalty0 (2-4):\penalty0
  245--258, 2007{\natexlab{a}}.
\newblock \doi{10.1504/IJCSM.2007.016534}.

\bibitem[John and Knobloch(2007{\natexlab{b}})]{John2007a}
V.~John and P.~Knobloch.
\newblock On spurious oscillations at layers diminishing ({SOLD}) methods for
  convection--diffusion equations: Part {I}--{A} review.
\newblock \emph{Computer Methods in Applied Mechanics and Engineering},
  196\penalty0 (17):\penalty0 2197--2215, 2007{\natexlab{b}}.
\newblock \doi{10.1016/j.cma.2006.11.013}.

\bibitem[John and Knobloch(2008)]{John2008}
V.~John and P.~Knobloch.
\newblock On spurious oscillations at layers diminishing ({SOLD}) methods for
  convection--diffusion equations: Part {II}--analysis for {P1} and {Q1} finite
  elements.
\newblock \emph{Computer Methods in Applied Mechanics and Engineering},
  197\penalty0 (21):\penalty0 1997--2014, 2008.
\newblock \doi{10.1016/j.cma.2007.12.019}.

\bibitem[Lavery(1988)]{Lavery1988}
J.~E. Lavery.
\newblock Nonoscillatory solution of the steady-state inviscid {B}urgers'
  equation by mathematical programming.
\newblock \emph{J. Comput. Phys.}, 79\penalty0 (2):\penalty0 436--448, 1988.
\newblock \doi{10.1016/0021-9991(88)90024-1}.

\bibitem[Lavery(1989)]{Lavery1989}
J.~E. Lavery.
\newblock Solution of steady-state one-dimensional conservation laws by
  mathematical programming.
\newblock \emph{SIAM J. Numer. Anal.}, 26\penalty0 (5):\penalty0 1081--1089,
  1989.
\newblock \doi{10.1137/0726060}.

\bibitem[Lavery(1991)]{Lavery1991}
J.~E. Lavery.
\newblock Solution of steady-state, two-dimensional conservation laws by
  mathematical programming.
\newblock \emph{SIAM J. Numer. Anal.}, 28\penalty0 (1):\penalty0 141--155,
  1991.
\newblock \doi{10.1137/0728007}.

\bibitem[Moskona et~al.(1995)Moskona, Petrushev, and Saff]{Moskona1995}
E.~Moskona, P.~Petrushev, and E.~Saff.
\newblock The gibbs phenomenon for best {L1}-trigonometric polynomial
  approximation.
\newblock \emph{Constructive Approximation}, 11\penalty0 (3):\penalty0
  391--416, 1995.
\newblock \doi{10.1007/BF01208562}.

\bibitem[Muga and {van der Zee}(2017)]{Muga2017}
I.~Muga and K.~G. {van der Zee}.
\newblock Discretization of {L}inear {P}roblems in {B}anach {S}paces:
  {R}esidual {M}inimization, {N}onlinear {P}etrov-{Ga}lerkin, and {M}onotone
  {M}ixed {M}ethods.
\newblock \emph{arXiv preprint arXiv:1511.04400v2}, 2017.
\newblock URL \url{https://arxiv.org/pdf/1511.04400v2.pdf}.

\bibitem[Muga et~al.(2019)Muga, Tyler, and van~der Zee]{Muga2019}
I.~Muga, M.~J. Tyler, and K.~G. van~der Zee.
\newblock The discrete-dual minimal-residual method ({DDMR}es) for weak
  advection-reaction problems in {B}anach spaces.
\newblock \emph{Computational Methods in Applied Mathematics}, 19\penalty0
  (3):\penalty0 557--579, 2019.
\newblock \doi{10.1515/cmam-2018-0199}.

\bibitem[Richards(1991)]{Richards1991}
F.~Richards.
\newblock A {G}ibbs phenomenon for spline functions.
\newblock \emph{Journal of approximation theory}, 66\penalty0 (3):\penalty0
  334--351, 1991.
\newblock \doi{10.1016/0021-9045(91)90034-8}.

\bibitem[Roos et~al.(2008)Roos, Stynes, and Tobiska]{Roos2008}
H.-G. Roos, M.~Stynes, and L.~Tobiska.
\newblock \emph{Robust numerical methods for singularly perturbed differential
  equations}, volume~24 of \emph{Springer Series in Computational Mathematics}.
\newblock Springer-Verlag, Berlin, second edition, 2008.
\newblock Con{\-}vec{\-}tion-diffusion-reaction and flow problems.

\bibitem[Saff and Tashev(1999)]{Saff1999}
E.~B. Saff and S.~Tashev.
\newblock Gibbs phenomenon for best {$L_p$} approximation by polygonal lines.
\newblock \emph{East J. Approx.}, 5\penalty0 (2):\penalty0 235--251, 1999.

\bibitem[Singer(1970)]{Singer1970}
I.~Singer.
\newblock \emph{Best Approximation in Normed Linear Spaces by Elements of
  Linear Subspaces}, volume 171 of \emph{Grundlehren der mathematischen
  Wissenschaften}.
\newblock Springer-Verlag Berlin Heidelberg, 1970.
\newblock \doi{10.1007/978-3-662-41583-2}.
\newblock Translated by R. Georgescu.

\bibitem[Wilbraham(1848)]{Wilbraham1848}
H.~Wilbraham.
\newblock On a certain periodic function.
\newblock \emph{Cambridge and Dublin Math. J}, 3\penalty0 (198):\penalty0 1848,
  1848.

\end{thebibliography}
\end{document}